\documentclass[10pt]{amsart}
\allowdisplaybreaks[4]
\usepackage{amssymb,color}

\usepackage{diagbox}
\usepackage[colorlinks=true, citecolor=blue, linkcolor=blue]{hyperref}
\usepackage{fancybox}
\newlength{\querylen}
\def\Mba{\eMM{M_\beta}}
\def\Mbb{{\eMM{M_{\beta_1}}}}

\definecolor{c20}{rgb}{0.,0.7,0.}
\definecolor{c30}{rgb}{0.,0.,1.}
\definecolor{c40}{rgb}{1,0.1,0.7}
\definecolor{c50}{rgb}{1,0,0}
\definecolor{c60}{rgb}{0,0.9,0.1}
\definecolor{c70}{rgb}{0,0,0}

\newcommand{\kb}[1]{\boldsymbol{#1}}
\newcommand{\vk}[1]{\kb{#1}}

\newcommand{\ve}{\varepsilon}

\newcommand{\abs}[1]{\lvert #1 \rvert}
\newcommand{\ABs}[1]{ \biggl \lvert #1 \biggr \rvert}

\newcommand{\EE}[1]{\mathbb{E}\left\{ #1 \right\}}

\newcommand{\pk}[1]{\mathbb{P} \left\{#1 \right\} }

\newcommand{\R}{\mathbb{R}}

\newcommand{\ldot}{,\ldots,}

\newcommand{\limit}[1]{\lim_{#1 \to   \infty}}

\newcommand{\BQN}{\begin{eqnarray}}
\newcommand{\EQN}{\end{eqnarray}}
\newcommand{\BQNY}{\begin{eqnarray*}}
\newcommand{\EQNY}{\end{eqnarray*}}

\def\eE#1{\textcolor{c20}{#1}}
\def\CE#1{\textcolor{c30}{#1}}

\newcommand{\BS}{\begin{sat}}
\newcommand{\ES}{\end{sat}}
\newcommand{\BT}{\begin{theo}}
\newcommand{\ET}{\end{theo}}
\newcommand{\BK}{\begin{korr}}
\newcommand{\EK}{\end{korr}}

\newcommand{\BD}{\begin{de}}
\newcommand{\ED}{\end{de}}

\newcommand{\BIT}{\begin{itemize}}
\newcommand{\EIT}{\end{itemize}}
\newcommand{\BDI}{\begin{description}}
\newcommand{\EDI}{\end{description}}

\newcommand{\BRM}{\begin{remarks}}
\newcommand{\ERM}{\end{remarks}}

\newcommand{\BEL}{\begin{lem}}
\newcommand{\EEL}{\end{lem}}

\newtheorem{theo}{Theorem}[section]
\newtheorem{sat}[theo]{Proposition}
\newtheorem{de}[theo]{Definition}
\newtheorem{lem}[theo]{Lemma}

\newtheorem{korr}[theo]{Corollary}
\newtheorem{remark}[theo]{Remark}
\newtheorem{remarks}[theo]{Remarks}

\newcommand{\nelem}[1]{{Lemma \ref{#1}}}

\newcommand{\netheo}[1]{{Theorem \ref{#1}}}

\newcommand{\prooftheo}[1]{ \textbf{Proof of Theorem} \ref{#1} }

\newcommand{\prooflem}[1]{\textbf{Proof of Lemma} \ref{#1}}

\newcommand{\COM}[1]{}

\newcommand{\QED}{\hfill $\Box$ \\}

\def\rw{\rightarrow}

\def\IF{\infty}

\def\ehb#1{\textcolor{c50}{#1}}
\def\ehc#1{\textcolor{c50}{#1}}
\def\ehb#1{#1}
\def\ehc#1{#1}

\def\Kd#1{{#1}}

\def\HH{\mathcal{H}}
\def\PP{\mathcal{P}}

\def\rw{\rightarrow}

\def\IF{\infty}

\date{}

\def\rw{\rightarrow}

\def\Var{\text{Var}}

\def\LL{\mathcal{\rho}}
\def\vv{v}

\def\KD#1{\textcolor{c70}{#1}}
\def\EHE#1{\textcolor{c50}{#1}}

\def\eM#1{\textcolor{c50}{#1}}
\def\eM#1{#1}
\def\eMM#1{{#1}}
\def\eMM#1{#1} 
\def\HEH#1{{#1}}
\def\NE#1{{#1}}

\begin{document}

\title[Gaussian Fields with regularly varying \Kd{dependence structure}] 
{Extremes of Gaussian Random Fields with regularly varying \Kd{dependence structure}}

\author{Krzysztof D\c{e}bicki}
\address{Krzysztof D\c{e}bicki, Mathematical Institute, University of Wroc\l aw, pl. Grunwaldzki 2/4, 50-384 Wroc\l aw, Poland}
\email{Krzysztof.Debicki@math.uni.wroc.pl}
\author{Enkelejd Hashorva}
\address{Enkelejd Hashorva, Department of Actuarial Science, University of Lausanne, UNIL-Dorigny 1015 Lausanne, Switzerland}
\email{enkelejd.hashorva.unil.ch}
\author{Peng Liu}
\address{Peng Liu, Mathematical Institute, University of Wroc\l aw, pl. Grunwaldzki 2/4, 50-384 Wroc\l aw, Poland
and Department of Actuarial Science, University of Lausanne, UNIL-Dorigny 1015 Lausanne, Switzerland}
\email{liupnankaimath@163.com}

\date{\today}

 \maketitle

\begin{quote}
	
{\bf Abstract}:
Let $X(t), t\in \mathcal{T}$ be a centered Gaussian random field with variance function $\sigma^2(\cdot)$ that attains its maximum at the unique point $t_0\in \mathcal{T}$, and let $M(\mathcal{T}):=\sup_{t\in \mathcal{T}} X(t)$.
For $\mathcal{T}$ a compact subset of $\R$, the current literature explains the \eMM{asymptotic} tail behaviour of $M(\mathcal{T})$ under some
\eMM{regularity}  conditions \eMM{including that} $1- \sigma(t)$ has a  polynomial decrease to 0 as $t \to  t_0$.
In this contribution we consider more general case that  $1- \sigma(t)$ is regularly varying at $t_0$.
We extend our analysis to random fields defined on some compact $\mathcal{T}\subset \R^2$, deriving the exact tail asymptotics of $M(\mathcal{T})$ for the class of Gaussian random fields
 with variance and correlation functions being regularly varying at $t_0$. A crucial novel element
is the analysis of families of Gaussian random fields
that do not possess locally additive dependence structures, which leads to qualitatively new types of asymptotics.\\

\end{quote}



\bigskip

\section{Introduction}
Let $X(t),t\ge 0$ be a centered stationary Gaussian processes with continuous trajectories,
unit variance and correlation function \KD{$r(\cdot)$} satisfying Pickands's condition
\BQN \label{PickandsC}
1-r(t) \sim  a\abs{t}^{\alpha} , \quad t \downarrow 0,  \quad a>0,  \quad \text{ and  } r(t)< 1, \forall t\not=0,
\EQN
with $\alpha \in (0,2]$;  in our notation  $\sim$ means asymptotic equivalence when the argument tends to 0 or $\IF$. \\
In the seminal contribution {\cite{PicandsA}} it was shown that
under \eqref{PickandsC}, for any $T$ positive
\begin{eqnarray}\label{eq1.2}
\pk{\sup_{t\in [0,T]} X(t)>u}\sim T\HH_{\alpha}a^{1/\alpha}u^{2/\alpha}\pk{X(0)> u}, \quad  \ u\rightarrow\infty,
\end{eqnarray}
with  the classical Pickands constant $\HH_{\alpha}$ defined by
$$\HH_{\alpha}=\lim_{T\rightarrow\infty} T^{-1} \EE{ \sup_{t\in[0,T]} e^{\sqrt{2}B_{\alpha}(t)-t^{\alpha}}},$$
where ${B_{\alpha}(t),t\ge 0}$  is a standard fractional Brownian motion with Hurst index $\alpha/2\in (0,1]$, see
 \cite{PicandsA,Pit72,Pit96,DE2002,DI2005,DE2014,DiekerY,DEJ14,Pit20, Tabis, DM, SBK} for various properties of $\HH_\alpha$.\\
\KD{The above finding was extended in various directions, including
$\alpha(t)$-locally-stationary Gaussian processes \eMM{(see \cite{MR2462286})},
and general  non-stationary Gaussian processes and random fields, \eMM{see e.g., \cite{Pit20}.}
A particularly important place in this theory is taken by the result of
Piterbarg \eMM{and Prisja{\v{z}}njuk} \cite{MR0494458}, where
the exact tail asymptotics of $\sup_{t\in [0,T]} X(t)$
is derived in the case that
the variance function $\sigma^2$ of a centered Gaussian process $X$ has a unique point of maximum in $[0,T]$, \eMM{say
$t_0$}.
More precisely, \eMM{for the correlation function it is assumed therein that} for some $\alpha \in (0,2]$
\BQN \label{countP}
1- \eM{r(s,t)} \sim a\abs{t-s}^{\alpha}, \quad s, t \downarrow 0, \quad a>0,
\EQN
whereas the behaviour of the variance function around the unique maximizer $t_0$ of
$\sigma^2(t)$ over $[0,T]$ \eMM{such that} $\sigma(t_0)=1$, is supposed  to satisfy
\BQN\label{eqA1}
1- \sigma(t)&\sim & b\abs{t}^{\beta}, \quad  t\downarrow 0, \quad b>0, \quad  \beta \in (0,\IF).
\EQN
}
{\eMM{Assume} further that the following H\"older continuity condition
\BQN\label{HOLDER}
\EE{(X(t)-X(t))^2} \le C \abs{t-s}^\nu
\EQN
is \eMM{valid} for all $s,t \in [0, \theta]$ with some $\theta\in (0,T]$ and $\nu\in (0,2]$},
\Kd{by \cite{MR0494458}, for $\alpha < \beta$}
\BQN\label{PLE:a}
\pk{ \sup_{t\in [0,T]} X(t)>u} \sim \HH_{{{\alpha}}}\frac{ a^{1/\alpha}}{b^{1/\beta}}\Gamma(1/{{\beta}}+1) u^{2/\alpha - 2/\beta}
\pk{X(0)> u},
\EQN
and \Kd{for} $\alpha = \beta$
\BQN\label{PLE:b}
\pk{ \sup_{t\in [0,T]} X(t)>u} \sim \PP_{\alpha}^{b/a}\pk{X(0)> u}, 
\EQN
where $\PP_{\alpha}^d,d>0$ is the Piterbarg constant defined by
$$ \PP_{\alpha }^d = \limit{S} \EE{ \sup_{t \in [0,S]} e^{\sqrt{2} B_\alpha(t)- (1+ d) t^\alpha}} \in (0,\IF).$$
\eMM{When $\alpha > \beta$, then} \eqref{PLE:b} holds with 1 instead of $\PP_{\alpha}^{b/a}$;
see \Kd{also} Theorem 2.1 in \cite{DEJ14} \eMM{for the case $T=\IF$}.\\
We note in passing that in fact {the H\"older continuity} (\ref{HOLDER})
is not needed to derive the \Kd{asymptotics} of (\ref{eq1.2}), which will be shown later in our main theorems;
{necessary and sufficient conditions that guarantee the
global  H\"older continuity of $X$ are presented in the deep contribution \cite{vitisari}.}\\

The original Pickands assumption \eqref{PickandsC}, and its counterpart \eqref{countP}
 can be relaxed to $1- r$ being regularly varying at 0 with index $\alpha \in \eMM{(0},2]$, see \cite{QuallsW, Berman92}. Specifically,
 in the case of \HEH{a} non-stationary $X$ we shall  assume
 for some \eMM{non-negative} $\LL\in \mathcal{R}_{\alpha/2}, \alpha \in (0,2]$
\BQN\label{eR}
1- r(s,t) &\sim& \LL^2(\abs{t-s}), \quad s,t\downarrow 0.
\EQN
Here  $f\in \mathcal{R}_\gamma$ means that the function $f$ is regularly varying at 0 with index $\gamma$, see \cite{Res, EKM97, Soulier} for details.\\
\Kd{The first goal of} this contribution is \Kd{an extension of} Piterbarg's \eMM{results} to a \KD{more general setup, that is to suppose that}
\BQN \label{eqsvv}
1- \sigma(t)\sim \vv^2(t), \quad t\downarrow 0,
\EQN
where $\vv\geq 0$ and $\vv \in \mathcal{R}_{\beta/2}, \beta>0$.
\Kd{In \netheo{mainT} we show}
that the \eMM{asymptotic} tail behaviour of
$\sup_{t\in [0,T]} X(t)$ can be determined under the assumption that $1- \sigma$ can be compared with $1- \rho$, \eMM{namely if
further}
\BQN\label{gamma}
\lim_{t\downarrow 0} \frac{\vv^2(t)}{\LL^2(t)}= \gamma \in [0,\IF].
\EQN
  \KD{Note that}, in Piterbarg's result mentioned above the limit $\gamma$ \eMM{is assumed to exist}.

\Kd{
Then we analyze tail distribution asymptotics of supremum of centered Gaussian random field $X(s,t),s\in[-T_1,T_1], t\in[\-T_2,T_2]$
with unique point that maximizes its variance function, say $(0,0)$.
Although extremes of Gaussian random fields with regularly varying correlation function are discussed in \cite{QuallsW},
see also \cite{ChengA,ChengC, YiminA,YiminB,StamatPopLet, Soja,MR3413855,Polnik} for new developments on extremes of Gaussian random fields,
most of the results in the existing literature are focused on the analysis of fields with {\it locally additive} dependence structure,
that is if
$$\Var(X(0,0))-\Var(X(s,t))\sim A_1|s|^{\beta_1} +A_1|t|^{\beta_2}$$
and
$$
1-{\rm Corr}(X(s,t),X(s_1,t_1))\sim B_1|s-s_1|^{\alpha_1}+B_2|t-t_1|^{\alpha_2}$$
as $s,s_1\to 0,$ $t,t_1\to 0.$
It appears that the investigation of fields that do not  satisfy this properties is considerably more delicate
and leads to qualitatively new results.
In Section 3 we derive several novel results \eM{concerned with the exact tail asymptotics of the
maximum of centered} Gaussian random fields when both variance and correlation functions
are regularly varying and \HEH{do not} possess a locally additive strucuture.}\\

Brief outline of the rest of the paper: Our main result for extremes of Gaussian processes is displayed in the Section 2, whereas Section 3 covers Gaussian random fields. The proofs of the theorems are presented in Section 4 and some technical results and their proofs are relegated to Appendix A and B.

\section{Gaussian Processes}
Before continuing with our investigation, we mention first that there are indeed  important cases of
Gaussian processes that satisfy our general setup in Section 1. \eMM{Indeed}, as remarked in \cite{HP99} and \cite{HP2004}, for any  function $\rho^2 \in \eMM{\mathcal{R}_\alpha}, \alpha\in (0,2]$ there exists a centered stationary Gaussian process $Y$ with continuous trajectories, unit variance and correlation function $r$ satisfying \eqref{eR}.
 Clearly, for any continuous function $\sigma(t),t\ge 0$ the process
 $X(t)=\sigma(t) Y(t),t\ge 0$ has continuous trajectories and variance function $\sigma^2$.

One \eMM{instance for the properties} of $\sigma$ is to assume that \eqref{eqsvv}
holds with $$\vv^2(t)= |\ln t|^c t^\beta, \quad A>0, c\in \mathbb{R}, \beta>0.$$ For such $\sigma$,
only the case $c=0$ can be dealt with using  Piterbarg's result mentioned in the Introduction. It is tempting to write
$$\vv^2(t)= (|\ln t|^{c/\beta} t) ^\beta.$$
\KD{Since in Piterbarg's result condition \eqref{eqA1}
explains the asymptotic expansion in \eqref{PLE:a}}
the  $u^{-2/\beta}$ term when $\alpha < \beta$, the above could imply that
\eqref{PLE:a} still holds if we replace $u^{-2/\beta}$ by $|\ln u|^{-2c/\beta^2} u^{-2/\beta} $.\\
 Detailed calculations show that \eMM{this intuition does not lead to the correct result, and in fact}
 the problem is much more complicated.
Indeed,
\KD{the tail asymptotics of the supremum is determined
in this case}
in terms of the (unique) asymptotic inverse of $v$, which is given by  (see Example 1.24 in \cite{Soulier} or Lemma 2 in \cite{Gira})
 \BQNY
\overleftarrow{v}(t)\sim \left(\frac{\beta}{2}\right)^{c/\beta}|\ln t|^{-c/\beta}t^{2/\beta}, \quad t\downarrow 0,
\EQNY
\COM{Consequently,
\BQNY
\int_0^{\overleftarrow{\vv}(u^{-1}g(u))}e^{-u^2v^2(t)}dt\sim \left(\frac{\beta}{2}\right)^{\frac{c}{\beta}}\Gamma(1/\beta+1)u^{-2/\beta}(\ln u)^{-c/\beta}.
\EQNY
}
\KD{where} $\overleftarrow{f}$ denotes the asymptotic (unique) inverse of $f\in \mathcal{R}_\gamma$.\\
Hereafter all regularly varying functions at 0 are assumed to be ultimately non-negative as $t\rw 0$.
\KD{Further $\Psi(u) \sim e^{-u^2/2} /(\sqrt{2 \pi }u)$, as $u \to \IF$, denotes the tail distribution of an $N(0,1)$ random variable, and
\eM{we set
	$$\PP_{\alpha}^\infty=:1, \quad \PP_{\alpha}^\infty[0,S]=:1, \quad  \alpha \in (0, 2], S>0.$$
}
}

We state next \eMM{the main result of  this section}.

\BT \label{mainT}
Let $X(t),t\ge 0$ be a centered Gaussian process with continuous trajectories and
variance {function} $\sigma^2$ having unique maximum at 0 with $\sigma(0)=1$.
 Suppose that $\sigma$ satisfies \eqref{eqsvv} and the correlation function $r$ of $X$  satisfies \eqref{eR}.
Assume further that
\KD{condition \eqref{gamma} is valid} for some $\gamma \in [0,\IF]$.  \\
i) If $\gamma=0$, then
\BQNY
\pk{ \sup_{t\in [0,T]} X(t)>u}\sim \ehc{\Gamma(1/\beta+1)\mathcal{H}_{\alpha}
\frac{\overleftarrow{v}(1/u)}{\overleftarrow{\LL}(1/u)}}
\Psi(u).
\EQNY
ii) If $\gamma\in (0,\IF]$, then
\BQNY
\pk{ \sup_{t\in [0,T]} X(t)>u}\sim \mathcal{P}_{\alpha}^\gamma\Psi(u) .
\EQNY
\ET

\BRM
i) If the maximum point of the variance is not 0, but an inner point, say $t_0\in (0,T)$
such that $\sigma(t_0)=1$, then
the results of \netheo{mainT} remain valid with  $\mathcal{H}_\alpha$
replaced by $\widehat{\mathcal{H}}_\alpha$
and $\mathcal{P}_\alpha^\gamma$ replaced by  $\widehat{\mathcal{P}}_\alpha^\gamma$,
where
$$ \widehat{\mathcal{H}}_\alpha= 2 \mathcal{H}_\alpha,  \quad
\widehat{ \PP}_{\alpha }^\gamma = \limit{S} \EE{ \sup_{t \in [-S,S]} e^{\sqrt{2} B_\alpha(t)- (1+ \gamma) t^\alpha}},
\quad \widehat \PP_{\alpha }^{\IF}=1.$$
ii) Since \netheo{mainT} remain valid if we substitute $v$ by an asymptotically equivalent $v^*$, we can assume that
$v^2(t)= \ell_\sigma(t) t^\beta$ with $\ell_\sigma$ a normalized slowly varying function (see e.g., \cite{BI1989, Soulier}).
  Similarly, let $\rho^2(t)=\ell_\rho(t) t^\alpha$ with $\ell_\rho$ another normalized slowly varying function.
Set next
$$
\ell_{\rho,\alpha}(x)= \sqrt{\ell_\rho(x^{1/\alpha})}, \quad \ell_{\sigma,\beta}(x)= \sqrt{\ell_\sigma(x^{1/\beta})}.$$
If further $\ell_{\sigma,\beta}^{\sharp}$ and
$ \ell_{\rho,\alpha}^{\sharp}$ denote the asymptotic inverses of $\ell_{\sigma, \beta}$ and $\ell_{\rho, \alpha}$, respectively then we have
$$v(x) = \ell_{\sigma,\beta}(x^\beta)x^{\beta/2}, \quad \rho(x) = \ell_{\rho,\alpha}(x^\alpha)x^{\alpha/2}$$
 and thus by Example 1.24 in \cite{Soulier} as $t\to 0$
 $$ \overleftarrow{v}(t) \sim  [\ell_{\sigma, \beta}^{\sharp}(t)]^{2/\beta}t^{2/\beta} ,  \quad \overleftarrow{\rho}(t) \sim  [\ell_{\rho,\beta}^{\sharp}(t)]^{2/\alpha}t^{2/\alpha}.$$
  Consequently,
$$ \frac{\overleftarrow{v}(1/u)}{\overleftarrow{\rho}(1/u)} \sim u^{2/\alpha- 2/\beta}
\frac{[\ell_{\sigma, \beta}^{\sharp}(1/u)]^{2/\beta}}
{[\ell_{\rho,\beta}^{\sharp}(1/u)]^{2/\alpha}}, \quad u\to \IF.$$

\ERM
 \netheo{mainT} is useful also when dealing with additive Gaussian random field. Specifically, assume that for \Kd{$T_1,T_2>0$}
$$ X(s,t)=  \eta_1(s) +  \eta_2(t), \quad s\in [-T_1,T_1], t\in [-T_2,T_2],$$
with $\eta_1,\eta_2$ two independent centered Gaussian random processes with continuous trajectories.
If both $\eta_1$ and $\eta_2$ are stationary satisfying \eqref{PickandsC}, or $\eta_1$ and $\eta_2$ satisfy the conditions of \netheo{mainT},  then
$$ \pk{\sup_{t\in \Kd{[-T_i,T_i]}} \eta_i(t) > u} \sim \mathcal{L}_i(u) u^{\tau_i} e^{- u^2/2} $$
 for some $\tau_i \ge -1$ with $\mathcal{L}_i(x)= C, x\ge 0$ if $\tau_i=-1$ and $\mathcal{L}_i$ is slowly varying at infinity if $\tau_i>-1$. Hence, since
$$ \sup_{s\in [-T_1,T_1], t\in [-T_2,T_2]} X(s,t)= \sup_{s\in [-T_1,T_1]} \eta_1(s) + \sup_{t\in [-T_2,T_2]} \eta_2(t),$$
then Lemma 2.3  in \cite{Farkas15} implies
\BQN\label{nuse}
\pk{ \sup_{s\in [-T_1,T_1], t\in [-T_2,T_2]} X(s,t)> u} \sim \sqrt{2 \pi} \mathcal{L}_1(u) \mathcal{L}_2(u) u^{\tau_1+ \tau_2-1} e^{-u^2/4}, \quad u\to \IF.
\EQN
In the particular case that $\eta_i$'s satisfy the conditions of \netheo{mainT} with $\rho_i, v_i, i=1,2$ instead of $\rho$ and $v$, where
\BQN \label{gammai}
\lim_{t\downarrow 0} \frac{\vv^2_i(t)}{\LL^2_i(t)}= \gamma_i \in [0,\IF], i=1,2,
\EQN
\eMM{then } \KD{ \eqref{nuse}} can be given more explicitly, see \netheo{th.T1} below. \\
\COM{
\BK \label{korrA} For the additive Gaussian random field $X$ above we have:\\
 i) If $\gamma_1=\gamma_2=0$, then
\BQNY
\pk{ \sup_{(s,t)\in [-T_1,T_1]\times[-T_2,T_2]} X(s,t)>u}\sim 4\prod_{i=1}^2 \Gamma(1/\beta_i+1)\mathcal{H}_{\alpha_i}
\frac{\overleftarrow{\vv}_i(1/u)}{\overleftarrow{\LL}_i(1/u)}\Psi(u).
\EQNY
ii) If $\gamma_1=0, \gamma_2 \in (0,\IF]$, then
\BQNY
\pk{ \sup_{(s,t)\in[-T_1,T_1]\times[-T_2,T_2]} X(s,t)>u}\sim2\Gamma(1/\beta_1+1)\mathcal{H}_{\alpha_1}\widehat{\mathcal{P}}^{\gamma_2}_{\alpha_2}\frac{\overleftarrow{\vv}_1(1/u)}{\overleftarrow{\LL}_1(1/u)}\Psi(u).
\EQNY

iii) If $\gamma_1, \gamma_2 \in (0,\IF]$, then
\BQNY
\pk{ \sup_{(s,t)\in [-T_1,T_1]\times[-T_2,T_2]} X(s,t)>u}\sim\widehat{\mathcal{P}}^{\gamma_1}_{\alpha_1}\widehat{\mathcal{P}}^{\gamma_2}_{\alpha_2}\Psi(u).
\EQNY
\EK
}
As we show in the next section, general Gaussian random fields are much more complex to deal with, and the results cannot be
derived from \netheo{mainT}.

\section{Gaussian Random Fields}
Extremes of \Kd{locally additive} Gaussian random fields with regularly varying correlation function are discussed in \cite{QuallsW}.
\KD{However there} are no results in the literature if the variance function is determined in terms of regularly varying functions
\Kd{and the dependence structure is non additive}.
In order to motivate our study, we consider first the additive Gaussian random field $X(s,t)= \eta_1(s)+ \KD{\eta_2}(t),
\Kd{ s\in [-T_1,T_1], t\in [-T_2,T_2]}$ \KD{introduced in Section 2}.
\KD{Thus, using that the variance function} $\sigma^2(s,t)$ of $X(s,t)$ is simply given by
$$\sigma^2(s,t)= \sigma_1^2(s)+ \sigma_2^2(t),$$
if $\eta_1,\eta_2$ satisfy the assumptions of \netheo{mainT}, then
$\sigma(s,t)$ achieves its unique  maximum at $(0,0)$.\\
In this section we shall discuss an extension of \netheo{mainT} to
\begin{eqnarray}
\pk{\sup_{(s,t) \in[-T_1,T_1]\times[-T_2,T_2]} X(s,t)>u}, \label{task.2}
\end{eqnarray}
as $u\to \IF$,
where $X(s,t),(s,t) \in[-T_1,T_1]\times[-T_2,T_2]$
is a centered  Gaussian random field, with variance function that is maximal on a unique point \Kd{but possess dependence structure that
is more complex than the additive one discussed above.}
\Kd{In particular, we suppose that}
\BQN\label{Cor2}
1- r(s,t,s_1,t_1) \sim \LL_{1}^2(\abs{a_{11}(s-s_1)+a_{12}(t-t_1)})+  \LL_{2}^2(\abs{a_{21}(s-s_1)+a_{22}(t-t_1)})
\EQN
as $s,s_1,t,t_1 \to 0$ with $\LL_{i}\geq 0$ and $\LL_{i}\in \mathcal{R}_{\alpha_i/2}, \alpha_i \in (0,2],i=1,2$.

For the variance function $\sigma^2(s,t)=Var(X(s,t))$  we shall assume that it attains \Kd{its} maximum at the unique point $(0,0)$ with $\sigma(0,0)=1$ and further
\BQN\label{Var2}
 1- \sigma(s,t) \sim \vv_{1}^2 (\abs{b_{11}(s-s_1)+b_{12}(t-t_1)})+ \vv_{2}^2(\abs{b_{21}(s-s_1)+b_{22}(t-t_1)}) , \quad s,t\downarrow0,
 \EQN
where $\vv_i\geq 0$ and $\vv_i\in \mathcal{R}_{\beta_i/2},\beta_i>0, i=1,2$.\\

\Kd{We note that recent results for the case that the variance function $\sigma^2$ is maximal on a line, which is the case for instance
if $\eta_1$ is stationary with unit variance 1 and $\eta_2$ satisfies the assumptions of \netheo{mainT},
are obtained in \cite{nonhomoANN}.
}

\KD{For further analysis it is useful to introduce the following matrices
 \BQNY
 A=\left(\begin{array}{cc}
 a_{11}& a_{12}\\
 a_{21}& a_{22}\\
 \end{array}
 \right), \ \
 B=\left(\begin{array}{cc}
 b_{11}& b_{12}\\
 b_{21}& b_{22}\\
 \end{array}
 \right).
 \EQNY
}
Let us observe that the assumption of uniqueness of
the maximizer of $\sigma(\cdot, \cdot)$ implies that ${\rm rank}(B)=2$.

We shall assume that \eqref{gammai} holds
and furthermore the following limits
\BQN\label{ette}
\lim_{t\downarrow 0}\frac{\rho_2^2(t)}{\rho_1^2(t)}= \eta \in [0,\IF], \ \ \lim_{t\downarrow 0}\frac{v_2^2(t)}{v_1^2(t)}=\theta \in [0,\IF]
\EQN
exist.

\KD{
It appears that the rank of matrix $A$ plays the key role
for the asymptotics of (\ref{task.2}), as $u\to\infty$.
Thus, in what follows, we shall distinguish
between two scenarios, when
${\rm rank}(A)=2$ and ${\rm rank}(A)=1$.
We exclude from further analysis the degenerated case of ${\rm rank}(A)=0$.
}
\subsection{\underline{Scenario I: ${\rm rank}(A)=2$}}

Suppose that $A$ is invertible and observe that
$Y(s,t):=X((A^{-1}(s,t)^\top )^\top )$ 
has, under  (\ref{Cor2}), (\ref{Var2}), correlation function such that
\BQN\label{Cor0}
1-r_{Y}(s,s_1,t,t_1)\sim \rho_1^2(|s-s_1|)+\rho^2_2(|t-t_1|), \ \ s,s_1,t,t_1\rw 0,
\EQN
and variance function $\sigma_Y^2$ satisfying
\BQN\label{Var0}
1-\sigma_Y(s,t)\sim v_1^2(|c_{11}s+c_{12}t|)+v_2^2(|c_{21}s+c_{22}t|), \ \ s,t\rw 0,
\EQN
with
$$ C=\left(\begin{array}{cc}
 c_{11}& c_{12}\\
 c_{21}& c_{22}\\
 \end{array}
 \right)
 =BA^{-1}.$$
Therefore, with no loss of generality,  in this section we tacitly assume that
$X$ satisfies (\ref{Cor2}) with
$$ A=\left(\begin{array}{cc}
 1& 0\\
 0& 1\\
 \end{array}
 \right)=:I.$$

Next, define an additive fractional Brownian field $W$ by
$$W(s,t)=\sqrt{2}B_{\alpha}(s)+\sqrt{2}\widetilde{B}_{\alpha}(t)- \abs{s}^\alpha- \abs{t}^\alpha, $$
where
$B_{\alpha}(t)$ and $\widetilde{B}_{\alpha}(t)$ are
independent standard fBm's with index $\alpha\in (0,2]$.
For a given matrix $D=(d_{ij})_{i,j=1,2}$, we define the generalized Piterbarg constant
$$\widehat{\mathcal{P}}_{\alpha}^{\gamma_1,\gamma_2, D}:=
\lim_{S\rw \IF}\EE{\sup_{(s,t)\in [-S,S]^2} e^{
W(s,t)-\gamma_1|d_{11}s+d_{12}t|^{\alpha}
-\gamma_2|d_{21}s+d_{22}t|^{\alpha}}},$$
where $\gamma_1, \gamma_2>0$. Note that if $det(\mathrm{D})\neq 0$, then there exists  $\gamma_3>0$ such that
$$\gamma_1|d_{11}s+d_{12}t|^{\alpha}
+\gamma_2|d_{21}s+d_{22}t|^{\alpha}\geq \gamma_3(|s|^{\alpha}+|t|^{\alpha}), \ \ s,t\in \mathbb{R},$$
which implies that
$\widehat{\mathcal{P}}_{\alpha}^{\gamma_1,\gamma_2, D}\leq \left(\widehat{\mathcal{P}}_{\alpha}^{\gamma_3}\right)^2<\IF.$
 Moreover,
\KD{for $D=I$} we have
$$
\widehat{\mathcal{P}}_{\alpha}^{\gamma_1,\gamma_2, I}=\widehat{ \mathcal{P}}_{\alpha}^{\gamma_1} \widehat{\mathcal{P}}_{\alpha}^{\gamma_2}.
$$

Let \eMM{for $S_1,S_2$ non-negative}
$$\mathcal{H}_\alpha^{\gamma_1,\gamma_2,b}(S_1,S_2):=
\EE{\sup_{(s+bt,t)\in [-S_1,S_1]\times[0,S_2]}  e^{ W(s,t) -\gamma_1|s+bt|^{\alpha}-\gamma_2|t|^{\alpha}
}} ,$$
and
$$\widehat{\mathcal{H}}_\alpha^{\gamma_1,\gamma_2,b}(S_1,S_2):=\EE{\sup_{(s+bt,t)\in [-S_1,S_1]\times[-S_2,S_2]}  e^{ W(s,t) -\gamma_1|s+bt|^{\alpha}-\gamma_2|t|^{\alpha}
}} .$$
In order to simplify the notation we set
$$\mathcal{H}_\alpha^{\gamma_1,\gamma_2,b}(S_1):=\mathcal{H}_\alpha^{\gamma_1,\gamma_2,b}(S_1,S_1), \ \ \widehat{\mathcal{H}}_\alpha^{\gamma_1,\gamma_2,b}(S_1)=\widetilde{\mathcal{H}}_\alpha^{\gamma_1,\gamma_2,b}(S_1,S_1), \ \  \mathcal{H}_\alpha^{\gamma_1,b}(S_1)=\mathcal{H}_\alpha^{\gamma_1,0,b}(S_1,S_1),$$
and
$$ \mathcal{H}_\alpha^{\gamma_1,b}:=\lim_{S\rw\IF}S^{-1}\mathcal{H}_\alpha^{\gamma_1,b}(S).$$

Now let us proceed to the analysis of (\ref{task.2}) for four special cases
whose proofs \eMM{are all different}, and
to which
one can reduce all other scenarios
(as will be advocated at the end of this section).
\\
Since below $A$ is taken to be the identity matrix, the cases discussed below are defined by the different choices of the matrix $B$.\\
$\diamond$  {\underline{Case 1.}}
We say that $X$ is {\it locally additive}, if both \eqref{Cor2} and \eqref{Var2} hold with $A=B=I$.
\COM{
(\ref{Cor2}) and
(\ref{Var2}) are satisfied with
$$ A=\left(\begin{array}{cc}
 1& 0\\
 0& 1\\
 \end{array}
 \right), \quad
B=\left(\begin{array}{cc}
1 & 0\\
0& 1\\
\end{array}
\right)
 .$$
\\
\\
}
\KD{The result below holds} for any $\theta, \eta \in [0,\IF]$ defined in \eqref{ette}.
\BT\label{th.T1}
Suppose that $X$ is \KD{a} {\it locally additive} Gaussian random field. \\
 i) If $\gamma_1=\gamma_2=0$, then
\BQNY \label{T1a}
\pk{ \sup_{(s,t)\in [-T_1,T_1]\times[-T_2,T_2]} X(s,t)>u}\sim 4\prod_{i=1}^2 \left(\Gamma(1/\beta_i+1)\mathcal{H}_{\alpha_i}
\frac{\overleftarrow{\vv}_i(1/u)}{\overleftarrow{\LL}_i(1/u)}\right)\Psi(u).
\EQNY
ii) If $\gamma_1=0, \gamma_2 \in (0,\IF]$, then
\BQNY
\pk{ \sup_{(s,t)\in[-T_1,T_1]\times[-T_2,T_2]} X(s,t)>u}\sim2\Gamma(1/\beta_1+1)\mathcal{H}_{\alpha_1}\widehat{\mathcal{P}}^{\gamma_2}_{\alpha_2}\frac{\overleftarrow{\vv}_1(1/u)}{\overleftarrow{\LL}_1(1/u)}\Psi(u).
\EQNY

iii) If $\gamma_1, \gamma_2 \in (0,\IF]$, then
\BQNY
\pk{ \sup_{(s,t)\in [-T_1,T_1]\times[-T_2,T_2]} X(s,t)>u}\sim\widehat{\mathcal{P}}^{\gamma_1}_{\alpha_1}\widehat{\mathcal{P}}^{\gamma_2}_{\alpha_2}\Psi(u).
\EQNY
\ET
\begin{remark}\label{remN}
We note that by the use of change of coordinates
Theorem \ref{th.T1} covers all the combinations of values of $\gamma_1, \gamma_2$.
\end{remark}
$\diamond$  {\underline{Case 2.}}  Here we shall assume that
(\ref{Cor2}) and
(\ref{Var2}) are satisfied with
\begin{eqnarray}
A=I, \quad B=\left(\begin{array}{cc}
1& b_{12}\\
0& 1\\
\end{array}
\right),\ \rm{with}\ b_{12}\neq 0. \label{case2}
\end{eqnarray}

\BT\label{th.T2}
\KD{Suppose that \eqref{ette} is satisfied with $\eta\in(0,\infty)$, $\theta=0$  and (\ref{case2}) holds}.

i) If $\gamma_1=\gamma_2=0$, then
\[
\pk{ \sup_{(s,t)\in[-T_1,T_1]\times[-T_2,T_2]} X(s,t)>u}
\sim
4 \prod_{i=1}^2 \left(\Gamma(1/\beta_i+1)\mathcal{H}_{\alpha_i}
\frac{\overleftarrow{\vv}_i(1/u)}{\overleftarrow{\LL}_i(1/u)}\right) \Psi(u).
\]
ii) If $\gamma_2=0, \gamma_1\in (0,\IF]$, then
\BQNY
\pk{ \sup_{(s,t)\in[-T_1,T_1]\times[-T_2,T_2]} X(s,t)>u}
\sim
2\Gamma(1/\beta_1+1) \mathcal{H}_{\alpha_1}^{\gamma_1, b_{12} \eta^{-1/\alpha_1}}
 \frac{\overleftarrow{\vv}_1(1/u)}{\overleftarrow{\LL}_1(1/u)}\Psi(u) .
\EQNY
iii) If $\gamma_2 \in (0,\IF], \gamma_1=\IF$, then
\[
\pk{ \sup_{(s,t)\in[-T_1,T_1]\times[-T_2,T_2]} X(s,t)>u}
\sim
\widehat{\mathcal{P}}_{\eMM{\alpha_1}}^{\gamma_2(|b_{12}|^{\alpha_1}\eta^{-1}+1)^{-1}}
\Psi(u) .
\]
\ET
\begin{remark}
The above theorem  \KD{covers} all the possible combinations of values
of $\gamma_1,\gamma_2$, since
the assumption that $\eta\in(0,\infty)$, $\theta=0$
excludes cases
$\gamma_1=0,\gamma_2 \in (0,\IF]$ and
$\gamma_1\in(0,\infty),\gamma_2 \in (0,\IF]$.

\KD{Although the same asymptotics \eMM{are imposed} in $i)$ of Theorem \ref{th.T1} and $i)$ of Theorem \ref{th.T2},
their proofs require \HEH{a} substantially different approach. Thus we did not combine
\HEH{those}  cases in one result.}
\end{remark}
$\diamond$  {\underline{Case 3.}}
The assumptions on $A$ and $B$ are the same as in Case 2 above, however
\COM{ case We say that
$X\in T_3$ if
(\ref{Cor2}) and
(\ref{Var2}) are satisfied with
$ A=\left(\begin{array}{cc}
 1& 0\\
 0& 1\\
 \end{array}
 \right)$,
$
B=\left(\begin{array}{cc}
1& b_{12}\\
0& 1\\
\end{array}
\right)
 $, $b_{12}\neq 0$, and
 }
 we shall suppose  that   $\eta=0$, $\theta\in(0,\infty)$. \\
Since $\theta\in (0,\IF)$, we set $\beta=\beta_1=\beta_2$. Let $\mu \in (-\IF,\IF)$ be the point at which
$|1+b_{12}t|^{\beta}+\theta |t|^{\beta}$
attains its minimum over $(-\IF,\IF)$. We have $\mu\in [-1/|b_{12}|, 1/|b_{12}|]$.
Further, set
\BQN\label{M} \Mba =\inf_{t\in (-\IF,\IF)}\left(|1+b_{12}t|^{\beta}+\theta |t|^{\beta}\right)
\EQN
and define the two-sided Piterbarg-type constant
$$\widehat{\PP}_{\beta }^{g_s} = \limit{S}\widehat{\PP}_{\beta }^{g_s}[-S,S], \ \ \text{with} \ \ \widehat{\PP}_{\beta }^{g_s}[-S,S]
=\EE{ \sup_{t \in \Kd{[-S,S]}} e^{\sqrt{2} B_\beta(t)-  t^\beta-g_s(t)}}, \quad  S>0, s\geq 0,$$
where
$$g_s(t)=\theta^{-1} \gamma_2\left(|s+b_{12}t|^{\beta}+\theta|t|^{\beta}-|(1+b_{12}\mu) s|^{\beta}-\theta|\mu s|^{\beta}\right), s\geq 0, t\in \mathbb{R}.$$
Further, set
$$\mathcal{I}_\beta :=  \int_{-\IF}^{\IF}\widehat{\mathcal{P}}_{\beta}^{g_{|s|}}
 e^{-\frac{\gamma_2  \Mba }{\theta}|s|^{\beta}}ds
\in (0,\IF).$$
The finiteness of $\eMM{\mathcal{I}_\beta}$ follows from the fact that for any $\epsilon>0$, there exists a positive constant $c_\epsilon>0$ such that
$$g_s(t)+\epsilon |s|^{\beta}\geq c_{\epsilon}|t|^{\beta}, \ \ s\geq 0, t\in \mathbb{R}
$$
implying that $\widehat{\mathcal{P}}_{\beta}^{g_{s}}\leq \widehat{\mathcal{P}}_{\beta}^{c_\epsilon}e^{\epsilon s^{\beta}}<\IF$, \eMM{and thus} for $\epsilon \in (0, \theta^{-1} \gamma_2 \Mba )$
$$\mathcal{I}_\beta\leq 2\int_0^\IF
\widehat{\mathcal{P}}_{\beta}^{g_{s}}e^{-\frac{\gamma_2 \Mba }{\theta}s^{\beta}}ds\leq 2\widehat{\mathcal{P}}_{\beta}^{c_\epsilon}\int_0^\IF
e^{-(\frac{\gamma_2 \Mba }{\theta}-\epsilon)s^{\beta}}ds<\IF.$$

\BT\label{th.T3}
\KD{Suppose that  \eqref{case2} holds and \eqref{ette} is satisfied}
with $\eta=0$, $\theta\in(0,\infty)$.\\
i) If $\gamma_1=\gamma_2=0$, then
\[
\pk{ \sup_{(s,t)\in [-T_1,T_1]\times[-T_2,T_2]} X(s,t)>u}
\sim
4\prod_{i=1}^2 \left(\Gamma(1/\KD{\beta}+1)\mathcal{H}_{\alpha_i}
\frac{\overleftarrow{\vv}_i(1/u)}{\overleftarrow{\LL}_i(1/u)}\right) \Psi(u) .
\]
ii) If $\gamma_1=0, \gamma_2 \in (0,\IF]$, then
\[
\pk{ \sup_{(s,t)\in[-T_1,T_1]\times[-T_2,T_2]} X(s,t)>u} \sim
\mathcal{H}_{\alpha_1} \left(\frac{\gamma_2}{\theta}\right)^{1/\beta}
\mathcal{I}_\beta \frac{\overleftarrow{\vv}_1(1/u)}{\overleftarrow{\LL}_1(1/u)}\Psi(u),
\]
iii)) If $\gamma_1=0, \gamma_2 =\IF$, then
\[
\pk{ \sup_{(s,t)\in[-T_1,T_1]\times[-T_2,T_2]} X(s,t)>u} \sim
2\Gamma(1/\beta+1)\left( \Mba \right)^{-1/\beta}\mathcal{H}_{\alpha_1}
\frac{\overleftarrow{\vv_1}(u^{-1})}{\overleftarrow{\LL_1}(u^{-1})}\Psi(u),
\]
iv) If $\gamma_1\in (0,\IF], \gamma_2=\IF$, then
\[
\pk{ \sup_{(s,t)\in [-T_1,T_1]\times[-T_2,T_2]} X(s,t)>u}
\sim
\widehat{\mathcal{P}}^{\gamma_1 \Mba }_{\alpha_1}
 \Psi(u).
\]
\ET
\begin{remark}
Analogously to the Case 2,
the assumption that $\eta=0$, $\theta\in(0,\infty)$
excludes case
$\gamma_1\in (0,\IF],\gamma_2\in[0,\infty)$.
\end{remark}

$\diamond$  {\underline{Case 4.}}
Here  we still assume that $A=I$ but there are no restrictions on the invertible $B$.
\COM{We say that
$X\in T_4$ if
(\ref{Cor2}) and
(\ref{Var2}) are satisfied with
$ A=\left(\begin{array}{cc}
 1& 0\\
 0& 1\\
 \end{array}
 \right)$
 and
$
B=\left(\begin{array}{cc}
b_{11}& b_{12}\\
b_{21}& b_{22}\\
\end{array}
\right)
 $,
 and $\eta,\theta\in(0,\infty)$.
\\
\\
This case leads to the following asymptotics.
}

\BT\label{th.T4}
Suppose that \eqref{Cor2} and \eqref{Var2} hold with $A=I$ and $B$ an invertible matrix,
and
\eqref{ette} is satisfied with  $\eta,\theta\in(0,\infty)$. \\
i)   If $\gamma_1=\gamma_2=0$, then
\BQN\label{gg0}
\pk{ \sup_{(s,t)\in [-T_1,T_1]\times[-T_2,T_2]} X(s,t)>u}
\sim
\frac{4}{\abs{  \mathrm{det}(B)}}
\prod_{i=1}^2 \left(\Gamma(1/\beta_i+1)\mathcal{H}_{\alpha_i}
\frac{\overleftarrow{\vv}_i(1/u)}{\overleftarrow{\LL}_i(1/u)}\right) \Psi(u) .
\EQN
ii) If $\gamma_1,\gamma_2 \in (0,\IF)$ or $\gamma_1=\gamma_2=\IF$, then
\BQN \label{f1}
\pk{ \sup_{(s,t)\in [-T_1,T_1]\times[-T_2,T_2]} X(s,t)>u}
\sim
\widehat{\mathcal{P}}_{\KD{\alpha_1}}^{\gamma_1,\gamma_1\theta, B_{\eta,\alpha_1}}
\Psi(u) ,
\EQN
where $\widehat{\mathcal{P}}_{\alpha_1}^{\gamma_1,\theta\gamma_1, B_{\eta,\alpha_1}}=1$ if $\gamma_1=\gamma_2= \IF$ and
$B_{\eta,\KD{\alpha_1}}=\left(\begin{array}{cc}
b_{11}& b_{12}\eta^{-1/\alpha_1}\\
b_{21}& b_{22}\eta^{-1/\alpha_1}\\
 \end{array}\right).$
\ET

\subsubsection{Discussion} As mentioned above, all other cases for
${\rm rank}(A)=2$ can be reduced to
the analysis of the field of one of types covered by Case 1-4. For the sake of transparency, let us first consider
$ A=I$
and $B$ such that exactly one element $b_{ij}$
equals to 0. With no loss of generality, by a change of variables, we can assume that
$$
B=\left(\begin{array}{cc}
1& b_{12}\\
0& 1\\
\end{array}
\right)
 , \quad b_{12}\neq 0.$$
Then the following holds:\\
\begin{itemize}
\item[$\diamond$] \underline{$\theta=\IF$:}
The asymptotics of (\ref{task.2})
in this case is covered by Case 1 above, since by Lemma \ref{simple} we obtain
$$v_1^2(|s+b_{12}t|)+v_2^2(|t|)\sim v_1^2(|s|)+v_2^2(|t|), \quad s, t\rw 0.$$
\item[$\diamond$] \underline{$\eta=\IF$:} Let $Z(s,t)=X(s-b_{12}t,t)$, which is a {\it locally additive} Gaussian random field.
Indeed, it follows from Lemma \ref{simple} that
$$1-r_Z(s,t,s_1,t_1)\sim \rho_1^2(|s-s_1-b_{12}(t-t_1)|)+\rho_2^2(|t-t_1|)\sim \rho_1^2(|s-s_1|)+\rho_2^2(|t-t_1|),\quad s,t,s_1,t_1\rw 0, $$
and $1-\sigma_Z(s,t)\sim  v_1^2(|s|)+v_2^2(|t|), s, t\rw 0.$
\item[$\diamond$] \underline{$\theta=0, \eta=0$:} Let $Z(s,t)=X(s,\frac{t-s}{b_{12}})$.
Then, again by Lemma \ref{simple},  $Z$ is a {\it locally additive} Gaussian random field
with
    $$1-r_Z(s,t,s_1,t_1)\sim\rho_1^2(|s-s_1|)+|b_{12}|^{-\alpha_2}\rho_2^2(|t-t_1|),\quad  s,t,s_1,t_1\rw 0, $$
and
$1-\sigma_Z(s,t)\sim  |b_{12}|^{-\beta_2}v_2^2(|s|)+v_1^2(|t|), s, t\rw 0.$
\item[$\diamond$]  \underline{$\theta=0$, $\eta\in (0,\infty)$:}
This is covered by Case 2 above.
\item[$\diamond$]  \underline{$\theta\in(0,\infty)$, $\eta=0$:}
This is covered by Case 3 above.
\item[$\diamond$]  \underline{$\theta\in(0,\infty)$, $\eta\in(0,\infty)$:}
This is covered by Case 4 above.
\\
\end{itemize}

Next, let $A=I$ and
\COM{Suppose now that
$ A=\left(\begin{array}{cc}
 1& 0\\
 0& 1\\
 \end{array}
 \right)$
and}
 $b_{ij}\neq 0$ for $i,j=1,2$.
With no loss of generality we can assume that
$$B=\left(\begin{array}{cc}
1& b_{12}\\
b_{21}& 1\\
\end{array}
\right)
 , \quad b_{12}b_{21}\neq 0.$$
 Let us observe that $det(B)=1-b_{12}b_{21}\neq 0$, which will be used in several places below.
Then the following holds:
\begin{itemize}
\item[$\diamond$] \underline{$\theta=0$, $\eta=0$:} Let $Z(s,t)=X(s, \frac{t-s}{b_{12}})$.
Again by Lemma \ref{simple}  $Z$ is a {\it locally additive} Gaussian random field
    $$1-r_Z(s,t,s_1,t_1)\sim  \rho_1^2(|s-s_1|)+|b_{12}|^{-\alpha_2}\rho_2^2(|t-t_1|), \quad s,t,s_1,t_1\rw 0, $$
and
$1-\sigma_Z(s,t)\sim  \left|\frac{detB}{b_{12}}\right|^{\beta_2}v_2^2(|s|)+v_1^2(|t|), \quad  s, t\rw 0.$

\item[$\diamond$] \underline{$\theta=0$, $\eta\in (0,\IF)$:}
This is Case 2 with $v_2^2$ replaced by $|det (B)|^{\beta_2}v_2^2$.
Indeed, by Lemma \ref{simple}, we have
\BQNY v_1^2(|s+b_{12}t|)+v_2^2(|b_{21}s+t|)&=&v_1^2(|s+b_{12}t|)+v_2^2(|b_{21}(s+b_{12}t)+(1-b_{12}b_{21})t|)\\
&\sim& v_1^2(|s+b_{12}t|)+|det (B)|^{\beta_2}v_2^2(|t|), \quad s, t\rw 0.
\EQNY
\item[$\diamond$] \underline{$\theta=0$, $\eta=\IF$:} Let $Z(s,t)=X(s-b_{12}t,t)$.
Again, by Lemma \ref{simple},  $Z$ is a {\it locally additive} Gaussian random field with
$$1-r_Z(s,t,s_1,t_1)\sim  \rho_1^2(|s-s_1|)+\rho_2^2(|t-t_1|),\quad s,t,s_1,t_1\rw 0, $$
and
$1-\sigma_Z(s,t)\sim v_1^2(|s|)+v_2^2(|b_{21}s+(1-b_{12}b_{21})t|)\sim v_1^2(|s|)+|det (B)|^{\beta_2}v_2^2(|t|), \quad s, t\rw 0.$
\item[$\diamond$] \underline{$\theta\in (0,\IF), \eta=0$:} Let $Z(s,t)=X(s,t-b_{21}s)$.
This is Case 3 with
$$1-r_Z(s,t,s_1,t_1)\sim \rho_1^2(|s-s_1|)+\rho_2^2(|t-t_1|), \quad s,t,s_1,t_1\rw 0, $$
and
$1-\sigma_Z(s,t)\sim  |det (B)|^{\beta_1}v_1^2(|s+b_{12}(det (B))^{-1}t|)
+v_2^2(|t|),  \quad s, t\rw 0.$
\item[$\diamond$]  \underline{$\theta\in(0,\infty)$, $\eta\in(0,\infty)$:}
This is covered by Case 4.
\item[$\diamond$]  \underline{$\theta\in(0,\infty)$, $\eta=\IF$:} Let $Z(s,t)=X(\frac{t-s}{b_{21}},s)$.
This is Case 3 with
$$1-r_Z(s,t,s_1,t_1)\sim \rho_2^2(|s-s_1|)+|b_{21}|^{-\alpha_1}\rho_1^2(|t-t_1|), \quad s,t,s_1,t_1\rw 0, $$
and
$1-\sigma_Z(s,t)\sim  \left|\frac{det (B)}{b_{21}}\right|^{\beta_1}v_1^2(|s+(-det (B))^{-1}t|)+v_2^2(|t|), s, t\rw 0.$
\item[$\diamond$]  \underline{$\theta=\IF$, $\eta=0$:} Let $Z(s,t)=X(s,t-b_{21}s)$.
This is a {\it locally additive} Gaussian random field with $v_1^2$ substituted by $|det (B)|^{\beta_1}v_1^2$.
\item[$\diamond$]  \underline{$\theta=\IF$, $\eta\in (0,\IF)$:} By Lemma \ref{simple}
we have that this is Case 2 with
$$1-r_X(s,t,s_1,t_1)\sim \rho_1^2(|s-s_1|)+\rho_2^2(|t-t_1|), \quad s,t,s_1,t_1\rw 0, $$
and
\BQNY
1-\sigma_X(s,t)&=&v_2^2(|b_{21}s+t|)+v_1^2(|b_{21}^{-1}(b_{21}s+t)+(b_{12}-b_{21}^{-1})t|)\\
&\sim&v_2^2(|b_{21}s+t|)+v_1^2((b_{12}-b_{21}^{-1})t|)\\
&\sim&  |b_{21}|^{\beta_2}v_2^2(|s+(b_{21})^{-1}t|)+\left|\frac{det (B)}{b_{21}}\right|^{\beta_1}v_1^2(|t|), \quad s, t\rw 0.
\EQNY
\item[$\diamond$]  \underline{$\theta=\IF$, $\eta=\IF$:} Let $Z(s,t)=X(\frac{s-t}{b_{21}},t)$.
We have that $Z$ is a {\it locally additive} Gaussian random field with
$$1-r_Z(s,t,s_1,t_1)\sim |b_{21}|^{-\alpha_1}\rho_1^2(|s-s_1|)+\rho_2^2(|t-t_1|),\quad  s,t,s_1,t_1\rw 0, $$
and
$1-\sigma_Z(s,t)\sim v_2^2(|s|)+ \left|\frac{det (B)}{b_{21}}\right|^{\beta_1}v_1^2(|t|),   s, t\rw 0.$
\end{itemize}

\subsection{\underline{Scenario II: ${\rm rank}(A)=1$}}
Suppose that ${\rm rank}(A)=1$. \Kd{Clearly} it suffices to consider Gaussian random fields with
covariance {function} that satisfies (\ref{Cor2}) with
$ A=\left(\begin{array}{cc}
 1& 0\\
 0& 0\\
 \end{array}
 \right)$
 and variance function satisfying (\ref{Var2}).
We begin with the analysis of
two  special cases, to
which all other structures of
field $X$ can be reduced.
\\
\\
$\diamond$ {\underline{Case 5.}}
Here we shall assume that (\ref{Cor2}) and
(\ref{Var2}) are satisfied with
\begin{eqnarray}
 A=\left(\begin{array}{cc}
 1& 0\\
 0& 0\\
 \end{array}
 \right)
\ {\rm{and}} \
B=I.
\label{case5}
\end{eqnarray}
\COM{
B=\left(\begin{array}{cc}
1& 0\\
0& 1\\
\end{array}
\right).
 $
\\}

\BT\label{th.T5}
\KD{Suppose that \eqref{case5} holds.}
\\
i) If $\gamma_1=0$, then
\BQNY
\pk{ \sup_{(s,t)\in[-T_1,T_1]\times[-T_2,T_2]} X(s,t)>u}
\sim 2\Gamma(1/\beta_1+1)\mathcal{H}_{\alpha_1}
 \frac{\overleftarrow{\vv}_1(1/u)}{\overleftarrow{\LL}_1(1/u)}\Psi(u) .
\EQNY
ii) If  $\gamma_1 \in (0,\IF]$, then
\BQN \label{f1}
\pk{ \sup_{(s,t)\in [-T_1,T_1]\times[-T_2,T_2]} X(s,t)>u}
\sim
\widehat{\mathcal{P}}^{\gamma_1}_{\alpha_1} \Psi(u).
\EQN
\ET
$\diamond$ {\underline{Case 6.}} Here we shall assume that
(\ref{Cor2}) and
(\ref{Var2}) are satisfied with
\begin{eqnarray}
A=\left(\begin{array}{cc}
 1& 0\\
 0& 0\\
 \end{array}
 \right), \quad
B=\left(\begin{array}{cc}
1& b_{12}\\
0& 1\\
\end{array}
\right),
\ \rm{and}\  b_{12}\neq 0. \label{case6}
\end{eqnarray}
\BT\label{th.T6}
\KD{Suppose that \eqref{case6} holds and \eqref{ette} is satisfied with $\theta \in (0,\IF)$.} \\
i) If $\gamma_1=0$, then
\BQN
\pk{ \sup_{(s,t)\in[-T_1,T_1]\times[-T_2,T_2]} X(s,t)>u}
\sim
2 (\Mbb )^{-1/\beta_1}  \Gamma(1/\beta_1+1)\mathcal{H}_{\alpha_1}
 \frac{\overleftarrow{\vv}_1(1/u)}{\overleftarrow{\LL}_1(1/u)}\Psi(u) .
\EQN

ii) If  $\gamma_1 \in (0,\IF]$, then
\BQN \label{f1}
\pk{ \sup_{(s,t)\in [-T_1,T_1]\times[-T_2,T_2]} X(s,t)>u}
\sim
\widehat{\mathcal{P}}^{\eMM{\gamma_1} \Mbb }_{\alpha_1}
\Psi(u)
\EQN
with $M_\beta$ defined in (\ref{M}).
\ET
\subsubsection{Discussion}
Having analyzed the above special cases, we are now in position to
give the asymptotics of (\ref{task.2}) for general structure of $X$. Suppose first, analogously to Scenario I, that $X$ satisfies
(\ref{Cor2}) and
(\ref{Var2}) with
$ A=\left(\begin{array}{cc}
 1& 0\\
 0& 0\\
 \end{array}
 \right)$
 and
 exactly one element of matrix $B$ equals 0.
With no loss of generality we can assume that
$
B=\left(\begin{array}{cc}
1& b_{12}\\
0 & 1\\
\end{array}
\right)
 $, $b_{12}\neq 0$.
Then the following holds.
\begin{itemize}
\item[$\diamond$]  \underline{$\theta=0$:}
Let $Z(s,t)=X(s,\frac{t-s}{b_{12}})$. Then, by Lemma \ref{simple}, this is Case 5
with
$$1-r_Z(s,t,s_1,t_1)\sim \rho_1^2(|s-s_1|), s,t,s_1,t_1\rw 0, \ \ 1-\sigma_Z(s,t)\sim |b_{12}|^{-\beta_2}v_2^2(|s|)+v_1^2(|t|),
\quad s,t\rw 0.$$
\item[$\diamond$]  \underline{$\theta\in(0,\infty)$:}
 This is case 6.
\item[$\diamond$]  \underline{$\theta=\infty$:}
The asymptotics of (\ref{task.2})
in this case is the same as the asymptotis derived in Case 5.
Indeed, by Lemma \ref{simple}, we have
$$v_1^2(|s+b_{12}t|)+v_2^2(|t|)\sim v_1^2(|s|)+v_2^2(|t|), \quad s,t\rw 0.$$
\end{itemize}

Finally we discuss the other case 
where the matrix $B$ is such that $b_{ij}\neq 0$ for $i,j=1,2$.
Again with no loss of generality we can assume that
$
B=\left(\begin{array}{cc}
1& b_{12}\\
b_{21} & 1\\
\end{array}
\right)
 $, $b_{12}, b_{21}\neq 0$.
Then the following holds with $det(B)=1-b_{12}b_{21}\neq 0$:
\begin{itemize}
\item[$\diamond$]  \underline{$\theta=0$:}
Let $Z(s,t)=X(s,\frac{t-s}{b_{12}})$.
This is covered by Case 5.
$$1-r_Z(s,t,s_1,t_1)\sim \rho_1^2(|s-s_1|), s,t,s_1,t_1\rw 0, \ \ 1-\sigma_Z(s,t)\sim \left|\frac{det (B)}{b_{12}}\right|^{\beta_2}v_2^2(|s|)+v_1^2(|t|), \quad s,t\rw 0.$$
\item[$\diamond$]  \underline{$\theta\in(0,\infty)$:}
Let $Z(s,t)=X(s,t-b_{21}s)$. Then, by
Lemma \ref{simple}, $Z$ is as in Case 6 with
$$1-r_Z(s,t,s_1,t_1)\sim \rho_1^2(|s-s_1|), \quad s,t,s_1,t_1\rw 0, \ \
1-\sigma_Z(s,t)\sim |det (B)|^{\beta_1}v_1^2(|s+b_{12}(det (B))^{-1}t|)
+v_2^2(|t|), \quad s,t\rw 0.$$
\item[$\diamond$]  \underline{$\theta=\infty$:}
Let $Z(s,t)=X(s, t-b_{21}s)$ . This is Case 5 with
$$1-r_Z(s,t,s_1,t_1)\sim \rho_1^2(|s-s_1|), s,t,s_1,t_1\rw 0, \ \ 1-\sigma_Z(s,t)\sim |det( B)|^{\beta_1}v_1^2(|s|)+v_2^2(|t|),
\quad s,t\rw 0.$$
\end{itemize}

\section{Proofs}
\KD{
In the rest of this section by
$\mathbb{Q}, \mathbb{Q}_i>0, i=1,2,...$
we denote constants that may differ from line to line.
}
\def\gu{\ehc{\xi(u)}}

\prooftheo{mainT}  
We set, for $u>1$ and $\gu:=u^{-1} \ln u$
\KD{
$$E_u=[0,\overleftarrow{\vv}(\gu )],
\quad I_k(u)=[kS\overleftarrow{\rho}(u^{-1}), (k+1)S\overleftarrow{\rho}(u^{-1})],
k\in \mathbb{N}\cup \{0\}$$
}
and, for given $\ve\in (0,1/2)$, define
 $$u_{k,\epsilon}^{-}=u(1+(1-\epsilon)\inf_{t\in I_{k}(u)}\vv^2(t)), \quad u_{k,\epsilon}^{+}=u(1+(1+\epsilon)\sup_{t\in I_{k}(u)}\vv^2(t)), \quad N(u)=\left[\frac{\overleftarrow{\vv}(\gu
 )}{\overleftarrow{\rho}(u^{-1})S}\right]+1.$$
  For $L>0$ sufficiently small
\BQN\label{new1}
\EE{(\overline{X}(t)-\overline{X}(t))^2}\leq 2(1-r(s,t))\leq 4\rho^2(|t-s|)\leq \mathbb{Q}|t-s|^{\alpha/2}, \ \ s,t\in [0,L],
\EQN
which ensures the H\"{o}lder condition in a neighborhood of $0$.
By the fact that
$\sup_{t\in [\overleftarrow{\vv}(\gu ),T]}\sigma(t)\leq 1-\mathbb{Q}(\gu)^2$
for $u$ sufficiently large,  (\ref{new1}), Theorem 8.1 in \cite{Pit96},  \Kd{we have}
\BQNY\label{proof2}
\pk{ \sup_{t\in [\overleftarrow{\vv}(\gu ),L]} X(t)>u}
\le  
 \mathbb{Q} T u^{4/\alpha}\Psi\left(\frac{u}{1-\mathbb{Q} (\gu)^2}\right).
\EQNY
Moreover, in light of  Borell inequality (see e.g., \cite{GennaBorell}) and the fact that  $\sup_{t\in [L,T]}\sigma(t)\leq 1-\delta$ with $\delta>0$,
$$\pk{ \sup_{t\in [\overleftarrow{\vv}(\gu ),L]} X(t)>u}
\le  
 e^{-\frac{(u-a)^2}{2(1-\delta)}}$$
 with $a=\mathbb{E}\left(\sup_{t\in [0,T]}X(t)\right)$.

Consequently,  for all large $u$ we have
\BQN\label{proof1}
\pi(u)\leq \pk{ \sup_{t\in [0,T]} X(t)>u}\leq \pi(u)+  \mathbb{Q} T u^{4/\alpha}\Psi\left(\frac{u}{1-\mathbb{Q} (\gu)^2}\right),
\EQN
where $\pi(u)=\pk{ \sup_{t\in [0, \overleftarrow{\vv}( \gu )]} X(t)>u}$.

\COM{
 For $L>0$ sufficiently small
\BQNY
\EE{(\overline{X}(t)-\overline{X}(t))^2}\leq 2(1-r(s,t))\leq 4\rho^2(|t-s|)\leq \mathbb{Q}|t-s|^{\alpha/2}, \ \ s,t\in [0,L].
\EQNY
By the fact that
$\sup_{t\in [\overleftarrow{\vv}(\gu ),T]}\sigma(t)\leq 1-\mathbb{Q}(\gu)^2$
for $u$ sufficiently large and
in view of (\ref{HOLDER}) and  Theorem 8.1 in \cite{Pit96}, we have
\BQN\label{proof2}
\pk{ \sup_{t\in [\overleftarrow{\vv}(\gu ),L]} X(t)>u}\leq \pk{ \sup_{t\in [\overleftarrow{\vv}(\gu ),L]} \overline{X}(t)>\frac{u}{1-\mathbb{Q} (\gu)^2}}\leq \mathbb{Q} T u^{4/\alpha}\Psi\left(\frac{u}{1-\mathbb{Q} (\gu)^2}\right)
\EQN
for some $\mathbb{Q}>0$.
Since \KD{$0$} is the unique point on $[0,T]$ such that $\sigma(\KD{0})=1$ is maximal,
there exists  $\delta\in(0,1)$ such that
$\sup_{t\in [L,T]}\sigma(t)\leq 1-\delta$,
which together with Borel theorem (see, e.g., Theorem D.1 in \cite{Pit96})
implies that there exists a constant $a>0$ such that for $u$ sufficiently large,
\BQN\label{proof0}
\pk{ \sup_{t\in [\theta,T]} X(t)>u}\leq 2\Psi\left(\frac{u-a}{1-\delta}\right).
\EQN
\BQN\label{proof2}
\pk{ \sup_{t\in [\overleftarrow{\vv}(\gu ),L]} X(t)>u}\leq \pk{ \sup_{t\in [\overleftarrow{\vv}(\gu ),L]} \overline{X}(t)>\frac{u}{1-\mathbb{Q} (\gu)^2}}\leq \mathbb{Q} T u^{4/\alpha}\Psi\left(\frac{u}{1-\mathbb{Q} (\gu)^2}\right)
\EQN
}
Next we give the exact asymptotics of $\pi(u)$ subject to three different scenarios.

{\underline{Case i) $\gamma=0$}}. For any $u$ positive we have
\BQN\label{proof3}
\sum_{k=0}^{N(u)-1}\pk{ \sup_{t\in I_{k}(u)} X(t)>u}-\sum_{i=1}^2\Lambda_i(u)\leq \pi(u)\leq \sum_{k=0}^{N(u)}\pk{ \sup_{t\in I_{k}(u)} X(t)>u},
\EQN
where
$$\Lambda_1(u)=\sum_{k=0}^{N(u)}\pk{ \sup_{t\in I_{k}(u)} X(t)>u, \sup_{t\in I_{k+1}(u)} X(t)>u},$$ and $$\Lambda_2(u)=\sum_{0\leq k,l\leq N(u), l\geq
k+2}\pk{ \sup_{t\in I_{k}(u)} X(t)>u, \sup_{t\in I_{l}(u)} X(t)>u}.$$

\KD{
The main difference in comparison with the proofs of the classical
cases considered in the literature, as e.g.,  in \cite{Pit96},
is contained in the approximation given below.
}
By uniform convergence theorem (UCT) for regularly varying functions, see e.g., \cite{BI1989}, we have
$$\sup_{s,t\in I_{k}(u), 1\leq k\leq N(u)}\left|\frac{\vv^2(s)}{\vv^2(t)}-\left(\frac{s}{t}\right)^{\alpha}\right|\rw 0, \ \ u\rw\IF,$$
which implies that for any $\epsilon>0$ and for $u$ sufficiently large,
$$\frac{\vv^2(s)}{\vv^2(t)}\geq\left(\frac{k}{k+1}\right)^{\alpha}-\epsilon/2, \ \ s,t \in I_k(u), 1\leq k\leq N(u).$$
 Thus for any $\epsilon>0$, there exists $k_\epsilon\in \mathbb{N}$ such that
$$\inf_{t\in I_{k}(u)}\vv^2(t)\geq (1-\epsilon)\sup_{t\in I_{k}(u)}\vv^2(t),  \ \ k_\epsilon\leq k\leq N(u).$$
Let $\eMM{X_{u,k}}(t)=\overline{X}(kS\overleftarrow{\rho}(u^{-1})+t), t\in I_{0}(u)$ with $k\in \mathcal{K}_u=\{k, 0\leq k\leq N(u)\}$ and $h_k(u)=u_{k,\epsilon}^{-}$. In light of  Lemma \ref{Pickands1}, we have
\BQN\label{callup}
\lim_{u\rw\IF}\sup_{0\leq k\leq N(u)}\left|(\Psi(u_{k,\epsilon}^{-}))^{-1}\pk{ \sup_{t\in I_{0}(u)}\eMM{X_{u,k}(t)}>u_{k,\epsilon}^{-}}
- \mathcal{H}_{\alpha}[0,S]\right|=0.
\EQN
\eMM{Consequently, as $u\to \IF$}
\BQNY\label{proof4}
\sum_{k=0}^{N(u)}\pk{ \sup_{t\in I_{k}(u)} X(t)>u}&\leq& \sum_{k=0}^{N(u)}\pk{ \sup_{t\in I_{k}(u)} \overline{X}(t)>u_{k,\epsilon}^{-}}\nonumber\\
&\leq& \sum_{k=0}^{N(u)}\pk{ \sup_{t\in I_{0}(u)}\eMM{X_{u,k}}(t)
>u_{k,\epsilon}^{-}}\nonumber\\
&\sim&  \sum_{k=0}^{N(u)}\mathcal{H}_{\alpha}[0,S]\Psi(u_{k,\epsilon}^{-})\nonumber\\
&\sim& \mathcal{H}_{\alpha}[0,S]\Psi(u)\sum_{k=0}^{N(u)}e^{-u^2(1-\epsilon)\inf_{t\in I_{k}(u)}\vv^2(t)}.
\label{q1}
\EQNY
\eMM{Further by Lemma \ref{integ}}
\BQNY
\sum_{k=0}^{N(u)}\pk{ \sup_{t\in I_{k}(u)} X(t)>u}&\leq&
\mathcal{H}_{\alpha}[0,S]\Psi(u)\left(k_\epsilon+\frac{1}{\overleftarrow{\LL}(u^{-1})S}\sum_{k=k_\epsilon}^{N(u)}\int_{t\in I_k(u)}e^{-(1-\epsilon)^2u^2\vv^2(t)}dt\right),\nonumber\\
&\sim&
\mathcal{H}_{\alpha}[0,S]
\left(k_\epsilon+\frac{1}{\overleftarrow{\LL}(u^{-1})S}\int_0^{\overleftarrow{\vv}(\gu)}e^{-(1-\epsilon)^2u^2v^2(t)}dt\right)\Psi(u)\nonumber\\
&\sim& \Gamma(1/\beta+1)\mathcal{H}_{\alpha}\frac{\overleftarrow{\vv}(u^{-1})}{\overleftarrow{\LL}(u^{-1})}\Psi(u), \ \ u\rw\IF, S\rw\IF, \epsilon\rw 0.
\label{q1}
\EQNY

Similarly, we obtain
\BQN\label{proof5}
\sum_{k=0}^{N(u)-1}\pk{ \sup_{t\in I_{k}(u)} X(t)>u}\geq
\Gamma(1/\beta+1)\mathcal{H}_{\alpha}\frac{\overleftarrow{\vv}(u^{-1})}{\overleftarrow{\LL}(u^{-1})}\Psi(u)
\KD{(1+o(1))}, \ \ u\rw\IF, S\rw\IF.
\label{q2}
\EQN
{Next we focus on $\Lambda_i(u), i=1,2$. Let $\hat{u}_{k,-\epsilon}=\min (u_{k,-\epsilon}^-, u_{k+1,-\epsilon}^-)$.  Then by (\ref{callup}) the fact that
\BQNY
&&\pk{ \sup_{t\in I_{k}(u)} X(t)>u, \sup_{t\in I_{k+1}(u)} X(t)>u}\\
&& \ \ \leq \pk{ \sup_{t\in I_{k}(u)}\overline{X}(t)>\hat{u}_{k,-\epsilon}}+\pk{\sup_{t\in
I_{k+1}(u)}\overline{X}(t)>\hat{u}_{k,-\epsilon}}-\pk{\sup_{t\in I_{k}(u)\cup I_{k+1}(u)} \overline{X}(t)>\hat{u}_{k,-\epsilon}},
\EQNY
 we have
$$\lim_{S\rw\IF}\lim_{u\rw\IF}\sup_{0\leq k\leq N(u)}\frac{\pk{ \sup_{t\in I_{k}(u)} X(t)>u, \sup_{t\in I_{k+1}(u)} X(t)>u}}{\mathcal{H}_{\alpha}[0,S]\Psi(\hat{u}_{k,-\epsilon})}\leq \lim_{S\rw\IF}\left(2-\frac{\mathcal{H}_{\alpha}[0,2S]}{\mathcal{H}_{\alpha}[0,S]}\right)=0.$$}
{Therefore,
\BQNY\label{proof6}
\Lambda_1(u) &=&
o(1)\sum_{k=0}^{N(u)}\mathcal{H}_{\alpha}[0,S]\Psi(\hat{u}_{k,-\epsilon})\leq o(1)\sum_{k=0}^{N(u)}2\mathcal{H}_{\alpha}[0,S]\Psi(u_{k,-\epsilon})\\
&=& o\left(\frac{\overleftarrow{\vv}(u^{-1})}{\overleftarrow{\LL}(u^{-1})}\Psi(u)\right), \ \ u\rw\IF, S\rw\IF.
\EQNY
}
By (\ref{eR}) and applying Lemma \ref{uni} in Appendix, we have \HEH{(note that below $k, l$ take values up to $N_u$, therefore an uniform
	upper bound for approximating the summuands derived in \nelem{uni} is essential)}
{\BQNY\label{proof7}
\Lambda_2(u)&\leq& \sum_{0\leq k,l\leq N(u), l\geq k+2} \pk{ \sup_{t\in I_{k}(u)} \overline{X}(t)>u_{k,-\epsilon}^-, \sup_{t\in I_{l}(u)} \overline{X}(t)>u_{l,-\epsilon}^-}\\
&\leq& \sum_{0\leq k,l\leq N(u), l\geq k+2} \mathbb{Q}S^2\Psi\left(\hat{u}_{k,l,-\epsilon}\right)
e^{-\mathbb{Q}_1|(l-k)S|^{\alpha/2}}\nonumber\\
&\leq&  \mathbb{Q}S^2\sum_{0\leq k\leq N(u)}\Psi\left(u_{k,-\epsilon}^{-}\right)\sum_{l=1}^\IF
e^{-\mathbb{Q}_1(lS)^{\alpha/2}}\nonumber\\
&\leq& \mathbb{Q}S^2 e^{-\mathbb{Q}_2S^{\alpha/2}}\sum_{k=0}^{N(u)}\Psi\left(u_{k,-\epsilon}^-\right)\nonumber\\
&=&o\left(\frac{\overleftarrow{\vv}(u^{-1})}{\overleftarrow{\LL}(u^{-1})}\Psi(u)\right), \ \ u\rw\IF, S\rw\IF, \nonumber
\EQNY
with $\hat{u}_{k,l,-\epsilon}=\min(u_{k,-\epsilon}^-, u_{l,-\epsilon}^-)$.}
 By the above \HEH{calculations} both $\Lambda_1(u)$ and $\Lambda_2(u)$  are negligible.
\eMM{Hence the results displayed in (\ref{proof2})-(\ref{proof7}) establish} the claim.\\
{\underline{Case ii) $\gamma \in (0,\IF]$ }.
\KD{
The proof of this case
is the same as the proof of the corresponding counterpart of Theorem D2 in \cite{Pit96},
with the exception that
\[
\pi(u)
\sim
\pk{ \sup_{t\in I_{0}(u)} X(t)>u}
\sim
\mathcal{P}_{\alpha}^{\gamma}[0, S]\Psi(u), \ \ u\rw\IF,
\]
where the last asymptotics follows by
Lemma \ref{Pickands1}.
}
This completes the proof.  \QED

\subsection{Proofs of Theorems \ref{th.T1}, \ref{th.T2}, \ref{th.T3}, \ref{th.T4}, \ref{th.T5} \text{and} \ref{th.T6}}
Define next for $S,u$ positive 
\begin{eqnarray*}
I_{k,l}(u)&=&[\overleftarrow{\LL_1}(u^{-1})kS, \overleftarrow{\LL_1}(u^{-1})(k+1)S]\times[\overleftarrow{\LL_2}(u^{-1})lS,
\overleftarrow{\LL_2}(u^{-1})(l+1)S], k, l \in \mathbb{N}\cup \{0\},\\
I_k(u)&=&[\overleftarrow{\LL_1}(u^{-1})kS, \overleftarrow{\LL_1}(u^{-1})(k+1)S], \quad
J_k(u)=[\overleftarrow{\LL_2}(u^{-1})kS,\overleftarrow{\LL_2}(u^{-1})(k+1)S],
\end{eqnarray*}
\KD{for $k\in\mathbb{N}\cup\{0\}$,}
and
$$ N_1(u)=\left[\frac{\overleftarrow{v_1}(u^{-1}\ln u)}{ \overleftarrow{\LL_1}(u^{-1})S}\right], \quad
N_2(u)=\left[\frac{ \overleftarrow{v_2}(u^{-1}\ln u)}{ \overleftarrow{\LL_2}(u^{-1})S}\right].$$
Additionally, let
   $$\mathbb{V}_1(u)=\{(k,l,k_1,l_1): -N_1(u)-2\leq k\leq k_1\leq N_1(u)+1, -N_2(u)-2\leq l,l_1\leq N_2(u)+1, I_{k,l}\cap I_{k_1,l_1}=\emptyset\},$$
    $$\mathbb{V}_2(u)=\{(k,l,k_1,l_1): -N_1(u)-2\leq k\leq k_1\leq N_1(u)+1, -N_2(u)-2\leq l,l_1\leq N_2(u)+1, (k,l)\neq (k_1,l_1), I_{k,l}\cap I_{k_1,l_1}\neq\emptyset\},$$
     $$u^-_{k,l,\epsilon}=u(1+(1-\epsilon)\inf_{(s,t)\in I_{k,l}(u)}\left(v_1^2(|b_{11}s+b_{12}t|)+v_2^2(|b_{21}s+b_{22}t|)\right)),$$
$$u^{+}_{k,l,\epsilon}=u(1+(1+\epsilon)\sup_{(s,t)\in I_{k,l}(u)}\left(v_1^2(|b_{11}s+b_{12}t|)+v_2^2(|b_{21}s+b_{22}t|)\right)),$$
$$u^{1,-}_{k,\epsilon}=u(1+(1-\epsilon)\inf_{s\in I_{k}(u)}v_1^2(|s|)), \ \ u^{1,+}_{k,\epsilon}=u(1+(1+\epsilon)\sup_{s\in I_{k}(u)}v_1^2(|s|)),$$
$$u^{2,-}_{l,\epsilon}=u(1+(1-\epsilon)\inf_{s\in J_{l}(u)}v_2^2(|s|)), \ \ u^{2,+}_{l,\epsilon}=u(1+(1+\epsilon)\sup_{s\in J_{l}(u)}v_2^2(|s|)), k,l\in \mathbb{Z},$$ where $u^{\pm}_{k,l,\epsilon}$ varies according to $B$.

\KD{
In what follows \eMM{for a given Gaussian random random field $Z$ we write $\overline{Z}$ for the standardised random field.}
}

\KD{\\
The general strategy of proofs of Theorems \ref{th.T1}, \ref{th.T2}, \ref{th.T3}, \ref{th.T4}, \ref{th.T5} \text{and} \ref{th.T6}
agrees from the double-sum technique developed for Gaussian random fields
in e.g.,  \cite{Pit96}.
However the variance-covariance structure of some cases substantially differs
from the families of Gaussian random fields analyzed in \cite{Pit96} and requires a case-specific approach,
on which we focus below.
\\
Observe that for all Cases 1-6
\BQN\label{F1}
\pi_1(u)\leq\pk{ \sup_{(s,t)\in [-T_1,T_1]\times[-T_2,T_2]} X(s,t)>u}\leq \pi_1(u) +\pk{ \sup_{(s,t)\in \left([-T_1,T_1]\times[-T_2,T_2]\right)
\setminus D_u} X(s,t)>u},
\EQN
where
$$\pi_1(u)=\pk{ \sup_{(s,t)\in D_u} X(s,t)>u},\text{with} \ \  D_u=[-\overleftarrow{v_1}\left(u^{-1}\ln u\right), \overleftarrow{v_1}\left(u^{-1}\ln
u\right)]\times[-\overleftarrow{v_2}\left(u^{-1}\ln u\right), \overleftarrow{v_2}\left(u^{-1}\ln u\right)].$$
For Case 1-Case 3 and Case 5-Case 6, by (\ref{Var2}) for $u$ sufficiently large we have
$$\sup_{(s,t)\in \left([-T_1,T_1]\times[-T_2,T_2]\right) \setminus D_u}\sigma(s,t)\leq 1-\mathbb{Q}u^{-2}\ln^2 u.$$
For Case 4, in light of (\ref{Var2}) and Lemma \ref{simple1}, we have
$$\sup_{(s,t)\in \left([-T_1,T_1]\times[-T_2,T_2]\right) \setminus D_u}\sigma(s,t)\leq 1-\mathbb{Q}\inf_{(s,t)\in \left([-T_1,T_1]\times[-T_2,T_2]\right) \setminus D_u}\left(v_1^2(|s|)+v_2^2(|t|)\right)\leq 1-\mathbb{Q}u^{-2}\ln^2 u.$$
\eMM{It follows by the fact that $(0,0)$ is the unique maximizer of $\sigma$, Theorem 8.1 in \cite{Pit96} and Borell theorem
that}
\BQN\label{neww}
\pk{ \sup_{(s,t)\in \left([-T_1,T_1]\times[-T_2,T_2]\right)
\setminus D_u} X(s,t)>u} \le \mathbb{Q}T_1 T_2u^{4/\alpha_1+4/\alpha_2}\Psi\left(\frac{u}{1-2u^{-2}\ln^2 u}\right).
\EQN
\COM{
point thatBy the assumption of the Similarly as in (\ref{eq3}), there exists $\delta>0$ sufficiently small such that
$$\EE{(\overline{X}(s,t)-\overline{X}(s_1,t_1))^2}\leq \mathbb{Q}_1(|s-s_1|^{\alpha_1/2}+|t-t_1|^{\alpha_2/2}), \ \ (s,t), (s_1,t_1)\in [-\delta,\delta]^2.$$
Thus in light of Theorem 8.1 in \cite{Pit96}, we have, for $u$ sufficiently large and $\delta$ sufficiently small
\BQN\label{F2}
\pk{ \sup_{(s,t)\in [-\delta,\delta]^2 \setminus D_u} X(s,t)>u}&\leq& \pk{ \sup_{(s,t)\in [-\delta,\delta]^2 \setminus D_u}
\overline{X}(s,t)>\frac{u}{1-2u^{-2}\ln^2 u}} \nonumber \\
&\leq& \mathbb{Q}T T_1u^{4/\alpha_1+4/\alpha_2}\Psi\left(\frac{u}{1-2u^{-2}\ln^2 u}\right).
\EQN
There exists a positive constant $0<\delta_1<1$ such that for all Case 1-Case 6
$$\sup_{(s,t)\in \left([-T_1,T_1]\times[-T_2,T_2]\right) \setminus [-\delta,\delta]^2}\sigma(s,t)\leq 1-\delta_1,$$
which together with Borell theorem gives that
\BQN\label{F3}
\pk{ \sup_{(s,t)\in  \left([-T_1,T_1]\times[-T_2,T_2]\right) \setminus[-\delta,\delta]^2} X(s,t)>u}\leq 2\Psi\left(\frac{u-a}{1-\delta_1}\right),
\EQN
with $a>0$.
\\
}
Therefore, for all Cases 1-6 we focus on the asymptotics of $\pi_1(u)$ as $u\rw\IF$,
proving that it delivers the asymptotics of (\ref{task.2}) as $u\to\infty$.
}
\\

\prooftheo{th.T1}\\
\underline{ Case i).}
\KD{Suppose that $\gamma_1=\gamma_2=0.$}
For any $0<\epsilon<1/2$ and $u$ large enough we have
\BQN\label{main}
\pi_{1,\epsilon}(u)-\sum_{i=1}^2\Lambda_i'(u)\leq \pi_1(u)\leq \pi_{1,-\epsilon}(u),
\EQN
with
\BQNY
\pi_{1,\pm \epsilon}(u):&=&\sum_{k=-N_1(u)\pm 2}^{N_1(u)\mp 1}\sum_{l=-N_2(u)\pm 2}^{N_2(u)\mp 1}\mathbb{P}\left(\sup_{(s,t)\in
I_{k,l}(u)}\overline{X}(s,t)>u_{k,l,\epsilon}^{\pm}\right),\\
 \Lambda_1'(u):&=&\sum_{(k,l,k_1,l_1)\in \mathbb{V}_1(u) }\mathbb{P}\left(\sup_{(s,t)\in I_{k,l}(u)}\overline{X}(s,t)>u_{k,l,\epsilon}^-, \sup_{(s,t)\in
 I_{k_1,l_1}(u)}\overline{X}(s,t)>u_{k_1,l_1,\epsilon}^-\right),\\
 \Lambda_2'(u):&=&\sum_{(k,l,k_1,l_1)\in \mathbb{V}_2(u) }\mathbb{P}\left(\sup_{(s,t)\in I_{k,l}(u)}\overline{X}(s,t)>u_{k,l,\epsilon}^-, \sup_{(s,t)\in
 I_{k_1,l_1}(u)}\overline{X}(s,t)>u_{k_1,l_1,\epsilon}^-\right).
 \EQNY
By  UCT, for any $\epsilon>0$, there exist two constants $k_\epsilon, l_\epsilon\in \mathbb{N}$ such that
\BQN\label{eq4}
&&\inf_{t\in I_{k}(u)}\vv_1^2(t)\geq (1-\epsilon)\sup_{t\in I_{k}(u)}\vv_1^2(t), \ \  \inf_{t\in J_{l}(u)}\vv_2^2(t)\geq (1-\epsilon)\sup_{t\in J_{l}(u)}\vv_2^2(t),
\EQN
hold for $k_\epsilon\leq |k|\leq N_1(u)+2, l_\epsilon\leq |l|\leq N_2(u)+2.$
Let
$$X_{u,k,l}(s,t)=\overline{X}(kS\overleftarrow{\rho}_1(u^{-1})+s,lS\overleftarrow{\rho}_2(u^{-1})+t), \mathcal{K}_u=\{(k,l), |k|\leq N_1(u)+2, |l|\leq N_2(u)+2\},$$
$$h_{k,l}(u)=u_{k,l,\epsilon}^-, \mathcal{E}_u=I_{0,0}(u), d_u=0.$$
One can easily check that conditions of Lemma \ref{PIPI} \eMM{are}  satisfied \eMM{implying that}
\BQN\label{uniform}
\lim_{u\rw\IF}\sup_{(k,l)\in \mathcal{K}_u}\left|(\Psi(u_{k,l,\epsilon}^-))^{-1}\mathbb{P}\left(\sup_{(s,t)\in
I_{0,0}(u)}\overline{X}(kS\overleftarrow{\rho}_1(u^{-1})+s,lS\overleftarrow{\rho}_2(u^{-1})+t)>u_{k,l,\epsilon}^-\right)-
\prod_{i=1}^2\mathcal{H}_{\alpha_i}[0,S]\right|=0.\EQN
Further, using Lemma \ref{integ} we have
\BQN\label{new}
\pi_{1,-\epsilon}(u)
&& \ = \sum_{k=-N_1(u)-2}^{N_1(u)+ 1}\sum_{l=-N_2(u)-2}^{N_2(u)+ 1}\mathbb{P}\left(\sup_{(s,t)\in
I_{0,0}(u)}\overline{X}(kS\overleftarrow{\rho}_1(u^{-1})+s,lS\overleftarrow{\rho}_2(u^{-1})+t)>u_{k,l,\epsilon}^-\right)\nonumber\\
&& \ \ \sim \sum_{k=-N_1(u)- 2}^{N_1(u)+ 1}\sum_{l=-N_2(u)- 2}^{N_2(u)+
1}\prod_{i=1}^2\mathcal{H}_{\alpha_i}[0,S]\Psi(u_{k,l,\epsilon}^-)\nonumber\\
&& \ \ \sim\prod_{i=1}^2\mathcal{H}_{\alpha_i}[0,S]\Psi(u)\left(R_1(u)+R_2(u)+\sum_{k_\epsilon\leq|k|\leq N_1(u)+2}
\frac{1}{\overleftarrow{\LL_1}(u^{-1})S}\int_{s\in I_k(u)}e^{-(1-\epsilon)^2u^2v_1^2(|t|)}dt\right.\nonumber\\
&& \ \ \ \left.\times\sum_{l_\epsilon\leq |l|\leq N_2(u)+2}\frac{1}{\overleftarrow{\LL_2}(u^{-1})S}\int_{t\in J_l(u)}e^{-(1-\epsilon)^2u^2v_2^2(t)}dt\right)\nonumber\\
&& \ \ \sim
4\prod_{i=1}^2\mathcal{H}_{\alpha_i}[0,S]\Psi(u)\frac{1}{\overleftarrow{\LL_1}(u^{-1})S}\int_0^{
\overleftarrow{v_1}(u^{-1}\ln u)}e^{-(1-\epsilon)^2u^2v^2_1(t)}dt\nonumber\\
&& \ \ \ \  \ \times \frac{1}{\overleftarrow{\LL_2}(u^{-1})S}\int_0^{
\overleftarrow{v_2}(u^{-1}\ln u)}e^{-(1-\epsilon)^2u^2v^2_2(t)}dt\nonumber\\
&& \ \ \sim
4(1-\epsilon)^{-1/\beta_1-1/\beta_2}\Gamma(1/\beta_1+1)\Gamma(1/\beta_2+1)S^{-2}\prod_{i=1}^2\mathcal{H}_{\alpha_i}[0,S]
\frac{\overleftarrow{\vv}_1(1/u)}{\overleftarrow{\LL}_1(1/u)}\frac{\overleftarrow{\vv}_2(1/u)}{\overleftarrow{\LL}_2(1/u)}\Psi(u)\nonumber\\
&& \ \ \sim 4\Gamma(1/\beta_1+1)\Gamma(1/\beta_2+1)\prod_{i=1}^2\mathcal{H}_{\alpha_i}
\frac{\overleftarrow{\vv}_1(1/u)}{\overleftarrow{\LL}_1(1/u)}\frac{\overleftarrow{\vv}_2(1/u)}{\overleftarrow{\LL}_2(1/u)}\Psi(u), \ \ u\rw \IF,
S\rw\IF, \epsilon\rw 0,\label{mainpart1}
\EQN
where
$$R_1(u)=\sum_{|k|\leq k_\epsilon}\sum_{|l|\leq N_2(u)+2}e^{-(1-\epsilon)u^2\inf_{s\in I_{k}(u)}v_1^2(|s|)-(1-\epsilon)u^2\inf_{t\in J_l(u)}v_2^2(|t|)},$$
$$R_2(u)=\sum_{|k|\leq N_1(u)+2}\sum_{|l|\leq l_\epsilon}e^{-(1-\epsilon)u^2\inf_{s\in I_{k}(u)}v_1^2(|s|)-(1-\epsilon)u^2\inf_{t\in J_l(u)}v_2^2(|t|)}.$$
Note that (\ref{new}) holds since  in light of Lemma \ref{integ} we have
\BQNY
R_1(u)&\leq& (2k_\epsilon+1) \left(2l_\epsilon+1+\frac{1}{\overleftarrow{\LL_2}(u^{-1})S}\int_0^{
\overleftarrow{v_2}(u^{-1}\ln u)}e^{-(1-\epsilon)^2u^2v^2_2(t)}dt\right)\\
&\sim&(2k_\epsilon+1)(1-\epsilon)^{-1/\beta_2} \frac{\overleftarrow{\vv}_2(1/u)}{\overleftarrow{\LL}_2(1/u)S}=
o\left(\frac{\overleftarrow{\vv}_1(1/u)}{\overleftarrow{\LL}_1(1/u)}\frac{\overleftarrow{\vv}_2(1/u)}{\overleftarrow{\LL}_2(1/u)}\right), \ \ u\rw\IF,
\EQNY
and
\BQNY
R_2(u)&\leq& (2l_\epsilon+1) \left(2k_\epsilon+1+\frac{1}{\overleftarrow{\LL_1}(u^{-1})S}\int_0^{
\overleftarrow{v_1}(u^{-1}\ln u)}e^{-(1-\epsilon)^2u^2v^2_1(t)}dt\right)\\
&\sim&(2l_\epsilon+1)(1-\epsilon)^{-1/\beta_1} \frac{\overleftarrow{\vv}_1(1/u)}{\overleftarrow{\LL}_1(1/u)S}\\
&=&
o\left(\frac{\overleftarrow{\vv}_1(1/u)}{\overleftarrow{\LL}_1(1/u)}\frac{\overleftarrow{\vv}_2(1/u)}{\overleftarrow{\LL}_2(1/u)}\right), \ \ u\rw\IF.
\EQNY
Similarly,
\BQN\label{Addition1}
\pi_{1,\epsilon}(u)&\geq& \sum_{k=-N_1(u)+1}^{N_1(u)-1}\sum_{l=-N_2(u)+1}^{N_2(u)- 1}\mathbb{P}\left(\sup_{(s,t)\in
I_{k,l}(u)}\overline{X}(s,t)>u_{k,l,\epsilon}^+\right)\nonumber\\ &\sim& 4\Gamma(1/\beta_1+1)\Gamma(1/\beta_2+1)\prod_{i=1}^2\mathcal{H}_{\alpha_i}
\frac{\overleftarrow{\vv}_1(1/u)}{\overleftarrow{\LL}_1(1/u)}\frac{\overleftarrow{\vv}_2(1/u)}{\overleftarrow{\LL}_2(1/u)}\Psi(u),
\EQN
as $u\rw\IF, S\rw\IF, \epsilon\rw 0$.

\KD{Next we prove that both $\Lambda'_1(u),\Lambda'_2(u)$ are asymptotically negligible.}
\eMM{From (\ref{Cor2})}, applying Lemma \ref{uni} in the Appendix,
with
$$\hat{u}_{k,l,k_1,l_1,\epsilon}=\min(u_{k,l,\epsilon}^-, u_{k_1,l_1,\epsilon}^-), \quad
\beta^*=\min(\alpha_1,\alpha_2),$$
we obtain
\BQN\label{lambda1}
\Lambda_1'(u)
&\leq& \mathbb{Q}S^{4} \sum_{(k,l,k_1,l_1)\in \mathbb{V}_1(u) }
\Psi\left(\hat{u}_{k,l,k_1,l_1,\epsilon}\right)e^{-\mathbb{Q}_1(|k-k_1|^2+|l-l_1|^2)^{\beta^*}S^{\beta^*/2}}\nonumber\\
&\leq& \mathbb{Q}S^{4}
\sum_{k=-N_1(u)-2}^{N_1(u)+1}\sum_{l=-N_2(u)-2}^{N_2(u)+1}\Psi(u_{k,l,\epsilon}^-)\sum_{m+n\geq 1, m,n \geq 0}
e^{-\mathbb{Q}_1(|m|^2+|n|^2)^{\beta^*}S^{\beta^*/2}}\nonumber\\
&\leq &\mathbb{Q}S^{4}
\sum_{k=-N_1(u)-2}^{N_1(u)+1}\sum_{l=-N_2(u)-2}^{N_2(u)+1}\Psi(u_{k,l,\epsilon}^-)e^{-\mathbb{Q}_2S^{\beta^*/2}}\nonumber\\
&=&o\left(\pi_{1,\epsilon}(u)\right), \ \ u\rw\IF, S\rw\IF.
\EQN

Now we focus on $\Lambda_2(u)$. Without loss of generality, we assume that $k_1=k+1$. Then let, for $k_1, l_1 \in \mathbb{N}\cup \{0\}$,
$$I^{(1)}_{k_1,l_1}=[\overleftarrow{\LL_1}(u^{-1})k_1S, \overleftarrow{\LL_1}(u^{-1})\left(k_1S+\sqrt{S}\right)]
\times[\overleftarrow{\LL_2}(u^{-1})l_1S, \overleftarrow{\LL_2}(u^{-1})(l_1+1)S],$$
$$I^{(2)}_{k_1,l_1}=[\overleftarrow{\LL_1}(u^{-1})\left(k_1S+\sqrt{S}\right),
\overleftarrow{\LL_1}(u^{-1})(k_1+1)S]\times[\overleftarrow{\LL_2}(u^{-1})l_1S, \overleftarrow{\LL_2}(u^{-1})(l_1+1)S].$$ For $(k,l,k_1,l_1)\in
\mathbb{V}_2(u), k_1=k+1$, we have
\BQNY
&&\mathbb{P}\left(\sup_{(s,t)\in I_{k,l}(u)}\overline{X}(s,t)>u_{k,l,\epsilon}^-, \sup_{(s,t)\in
I_{k_1,l_1}(u)}\overline{X}(s,t)>u_{k_1,l_1,\epsilon}^-\right)\\
&& \ \ \leq
\mathbb{P}\left(\sup_{(s,t)\in I_{k,l}(u)}\overline{X}(s,t)>u_{k,l,\epsilon}^-, \sup_{(s,t)\in
I^{(1)}_{k_1,l_1}(u)}\overline{X}(s,t)>u_{k_1,l_1,\epsilon}^-\right)\\
&& \ \ \ \ +
\mathbb{P}\left(\sup_{(s,t)\in I_{k,l}(u)}\overline{X}(s,t)>u_{k,l,\epsilon}^-, \sup_{(s,t)\in
I^{(2)}_{k_1,l_1}(u)}\overline{X}(s,t)>u_{k_1,l_1,\epsilon}^-\right)\\
&& \ \ :=p_{k,l,k_1,l_1}^{(1)}(u)+p_{k,l,k_1,l_1}^{(2)}(u).
\EQNY
It follows from Lemma \ref{PIPI} that
\BQNY
p_{k,l,k_1,l_1}^{(1)}(u)\leq \mathbb{P}\left( \sup_{(s,t)\in I^{(1)}_{k_1,l_1}(u)}\overline{X}(s,t)>u_{k_1,l_1,\epsilon}^-\right)\sim
\mathcal{H}_{\alpha_1}[0,\sqrt{S}]\mathcal{H}_{\alpha_2}[0,S]\Psi(u_{k,l,\epsilon}^-).
\EQNY
Further, since each $I_{k,l}(u)\times I_{k_1,l_1}(u)$ has at most $8$ neighbors, we have that
\BQNY
\lefteqn{\sum_{(k,l,k_1,l_1)\in \mathbb{V}_2(u)}p_{k,l,k_1,l_1}^{(1)}(u)
 \ \ \leq  8\sum_{k=-N_1(u)-2}^{N_1(u)+1}\sum_{l=-N_2(u)-2}^{N_2(u)+1}\Psi(u_{k,l,\epsilon}^-)}\\
&& \ \ \ \ \ \times \left(\mathcal{H}_{\alpha_1}[0,\sqrt{S}]
\mathcal{H}_{\alpha_2}[0,S]+
\mathcal{H}_{\alpha_1}
[0, S]\mathcal{H}_{\alpha_2}[0,\sqrt{S}]\right)\\
&& \ \ \leq 8\left(\frac{\mathcal{H}_{\alpha_1}[0,\sqrt{S}]}{\mathcal{H}_{\alpha_1}[0,S]}
+\frac{\mathcal{H}_{B_{\alpha_2}}
[0,\sqrt{S}]}{\mathcal{H}_{\alpha_2}[0,S]}\right)\nonumber\\
&& \ \ \ \ \ \times
\sum_{k=-N_1(u)-2}^{N_1(u)+1}\sum_{l=-N_2(u)-2}^{N_2(u)+1}\mathcal{H}_{\alpha_1}[0,S]
\mathcal{H}_{\alpha_2}[0,S]\Psi(u_{k,l,\epsilon}^-)\\
&& \ \ =o\left(\pi_{1,\epsilon}(u)\right), \ \ u\rw\IF, S\rw\IF.
\EQNY
In light of Lemma \ref{uni}, we have
\BQNY
\sum_{(k,l,k_1,l_1)\in \mathbb{V}_2(u)}p_{k,l,k_1,l_1}^{(2)}(u)
& \leq&  \mathbb{Q}S^{4}e^{-\mathbb{Q}_1S^{\beta^*/4}}
\sum_{(k,l,k_1,l_1)\in \mathbb{V}_2(u)}\Psi(\hat{u}_{k,l,k_1,l_1,\epsilon})\\
&\leq& \mathbb{Q}S^{4}e^{-\mathbb{Q}_1S^{\beta^*/4}}
\sum_{k=-N_1(u)-2}^{N_1(u)+1}\sum_{l=-N_2(u)-2}^{N_2(u)+1}\Psi(u_{k,l,\epsilon}^-)\\
& =& o\left(\pi_{1,\epsilon}(u)\right), \ \ u\rw\IF, S\rw\IF.
\EQNY
Consequently,
\BQN\label{lambda2}
\Lambda_2'(u)=o\left(\pi_{1,\epsilon}(u)\right), \ \ u\rw\IF, S\rw\IF.
\EQN
Combing (\ref{main}), (\ref{mainpart1}), (\ref{Addition1}) (\ref{lambda1}) with (\ref{lambda2}), we derive that
$$\pi_1(u)\sim 4\Gamma(1/\beta_1+1)\Gamma(1/\beta_2+1)\prod_{i=1}^2\mathcal{H}_{\alpha_i}
\frac{\overleftarrow{\vv}_1(1/u)}{\overleftarrow{\LL}_1(1/u)}\frac{\overleftarrow{\vv}_2(1/u)}{\overleftarrow{\LL}_2(1/u)}\Psi(u), \ \ u\rw \IF,$$
\eMM{hence the claim follows}. \\
\underline{ Case ii) $\gamma_1=0, \gamma_2\in (0,\IF)$.}  Let in the sequel
$$\widetilde{I}_{k,0}(u)=I_{k,0}(u)\cup I_{k,-1}(u), \quad \mathbb{V}_1^{(1)}(u)=\{(k,k_1):
-N_1(u)-2\leq k<k_1\leq N_1(u)+1, k_1-k\geq 2\},$$
and $$\mathbb{V}_2^{(1)}(u)=\{(k,k_1): -N_1(u)-2\leq k<k_1\leq N_1(u)+1, k_1=k+1\}.$$
\CE{For} any $0<\epsilon<1$ and all $u$ large enough
\BQN\label{main1}
\pi_{1,\epsilon}^{(1)}(u)-\sum_{i=1}^2\Lambda_i^{(1)}(u)\leq \pi_1(u)\leq \pi_{1,-\epsilon}^{(1)}(u)+\pi_{1,-\epsilon}^{(2)}(u),
\EQN
with
\BQNY
\pi_{1,\pm \epsilon}^{(1)}(u):&=&\sum_{k=-N_1(u)\pm 2}^{N_1(u)\mp1}\mathbb{P}\left(\sup_{(s,t)\in
\widetilde{I}_{k,0}(u)}\frac{\overline{X}(s,t)}{1+(1\pm\epsilon)v^2_2(t)}>u_{k,\epsilon}^{1,\pm}\right)\nonumber\\
\pi_{1,-\epsilon}^{(2)}(u):&=&\sum_{k=-N_1(u)-2}^{N_1(u)+1}\sum_{|l|=1, l\neq -1}^{N_2(u)+1}\mathbb{P}\left(\sup_{(s,t)\in
I_{k,l}(u)}\overline{X}(s,t)>u_{k,l,\epsilon}^-\right)\nonumber\\
 \HEH{\Lambda_i}^{(1)}(u):&=&\sum_{(k,k_1)\in \mathbb{V}_i^{(1)}(u) }\mathbb{P}\left(\sup_{(s,t)\in \widetilde{I}_{k,0}(u)}\overline{X}(s,t)>u_{k,\epsilon}^{1,-},
 \sup_{(s,t)\in \widetilde{I}_{k_1,0}(u)}\overline{X}(s,t)>u_{k_1,\epsilon}^{1,-}\right), \quad \HEH{i=1,2,}\nonumber\\
 \EQNY
 Set further  $X_{u,k}(s,t)=\overline{X}(k\overleftarrow{\rho}_1(u^{-1})S+s,t)$ and define
 $$ \mathcal{K}_u=\{k,  |k|\leq N_1(u)+2\}, \quad \mathcal{E}_u=\tilde{I}_{0,0}(u),  \quad h_k(u)=u_{k,\epsilon}^{1,-},  \quad d_u(s,t)=(1-\epsilon)v^2_2(t).$$
 Using Lemma {\ref{PIPI}}, we have
 $$\lim_{u\rw\IF}\sup_{k\in \mathcal{K}_u}\left|(\Psi(u_{k,\epsilon}^{1,-}))^{-1}\mathbb{P}\left(\sup_{(s,t)\in
\widetilde{I}_{k,0}(u)}\frac{\overline{X}(s,t)}{1+(1-\epsilon)v^2_2(t)}>u_{k,\epsilon}^{1,-}\right)-\mathcal{H}_{\alpha_1}[0,
S]\widehat{\mathcal{P}}^{\gamma_2(1-\epsilon)}_{\alpha_2}(S)\right|=0.$$

Further, by Lemma \ref{integ}, we have
\BQN\label{PG1}
\pi_{1,- \epsilon}^{(1)}(u)
&\sim&\sum_{k=-N_1(u)-2}^{N_1(u)+1}\mathcal{H}_{\alpha_1}[0,
S]\widehat{\mathcal{P}}^{\gamma_2(1-\epsilon)}_{\alpha_2}(S)\Psi(u_{k,\epsilon}^{1,-})\nonumber\\
&\leq& 2\mathcal{H}_{\alpha_1}[0,
S]\widehat{\mathcal{P}}^{\gamma_2(1-\epsilon)}_{\alpha_2}(S)
\left( 2k_{\epsilon}+1+\frac{\Psi(u)}{\overleftarrow{\LL}_1(u^{-1})S}\int_0^{\overleftarrow{v_1}(u^{-1}\ln u)}e^{-(1-\epsilon)^2u^2v^2_1(t)}dt\right)(1+o(1))\nonumber\\
&\sim&
2\Gamma(1/\beta_1+1)\mathcal{H}_{\alpha_1}\widehat{\mathcal{P}}^{\gamma_2}_{\alpha_2}\frac{\overleftarrow{\vv}_1(u^{-1})}{\overleftarrow{\LL}_1(u^{-1})
}\Psi(u), \ \ u\rw \IF,  \epsilon\rw 0, S\rw\IF.
\EQN
Similarly, as $u\rw \IF,  \eMM{\ve \to 0},  S\rw\IF,$
\BQN\label{Addition2}
\pi_{1, \epsilon}^{(1)}(u)
\sim
2\Gamma(1/\beta_1+1)\mathcal{H}_{\alpha_1}\widehat{\mathcal{P}}^{\gamma_2}_{\alpha_2}\frac{\overleftarrow{\vv}_1(u^{-1})}{\overleftarrow{\LL}_1(u^{-1})
}\Psi(u).
\EQN
Moreover,  by Lemma \ref{PIPI}
\BQN
\pi_{1,-\epsilon}^{(2)}(u)
&\sim& \sum_{k=-N_1(u)-2}^{N_1(u)+2}\sum_{|l|=1, l\neq -1}^{N_2(u)+1}\prod_{i=1}^2\mathcal{H}_{\alpha_i}[0,
S]\Psi(u_{k,l,\epsilon}^-)\nonumber\\
&\leq & 2\prod_{i=1}^2\mathcal{H}_{\alpha_i}[0, S]\Psi(u)
\sum_{k=-N_1(u)-2}^{N_1(u)+2}e^{-(1-\epsilon)u^2\inf_{s\in I_k(u)}v_1^2
(|s|)} \sum_{|l|=1, l\neq -1}^{N_2(u)+1}e^{-(1-\epsilon)u^2v_2^2(\overleftarrow{\LL_2}(u^{-1}))|lS|^{\beta_2'}}\nonumber\\
&\leq&2\prod_{i=1}^2\mathcal{H}_{\alpha_i}[0, S]\Psi(u)\sum_{k=-N_1(u)-2}^{N_1(u)+ 1}
e^{-(1-\epsilon)u^2\inf_{s\in I_k(u)}v_1^2
(|s|)}\sum_{|l|=1, l\neq -1}^{N_2(u)+1}e^{-(1-2\epsilon)\gamma_2|lS|^{\beta_2'}}\label{Addition3}\\
&\leq& 4\prod_{i=1}^2\mathcal{H}_{\alpha_i}\frac{\overleftarrow{\vv}_1(u^{-1})\Psi(u)}{(1-\epsilon)^{2/\beta_1}\overleftarrow{\LL}_1(u^{-1})
}S^2 e^{-\mathbb{Q}_4S^{\beta_2'}}\nonumber\\
&=&o\left(\pi_{1,\epsilon}^{(1)}(u)\right)
\label{PG2}
\EQN
 as $u\rw\IF, S\rw\IF$ and $\eMM{\ve \to 0}$,  with $0<\beta_2'<\beta_2$.
Note that in (\ref{Addition3}) we use  Potter's \HEH{bounds} (see e.g., \cite{Res})
 for regularly varying function $v_2(t)$ at zero to derive that, for $u$ and $S$ large enough,
\BQN\label{A0}
\frac{v_2^2(\overleftarrow{\LL_2}(u^{-1})lS)}{v_2^2(\overleftarrow{\LL_2}(u^{-1}))}\geq \left(lS\right)^{\beta_2'}
\EQN
holds for $1\leq l\leq N_2(u) $.
Using Lemma \ref{uni}, we have
\BQN\label{lambda11}
\Lambda_1^{(1)}(u)
&\leq& \mathbb{Q}S^{4}
\sum_{(k,k_1)\in \mathbb{V}_1^{(1)}(u) }\Psi
\left(u_{k,k_1,\epsilon}^1\right)e^{-\mathbb{Q}_1(|k-k_1|S)^{\beta^*/2}}\nonumber\\
&\leq& \mathbb{Q}S^{4}
\sum_{k=-N_1(u)-2}^{N_1(u)+1}\Psi(u_{k,\epsilon}^{1,-})\sum_{m\geq 1} e^{-\mathbb{Q}_1m^{\beta^*/2}S^{\beta^*/2}}\nonumber\\
&\leq & \mathbb{Q}S^{4}
\sum_{k=-N_1(u)-2}^{N_1(u)+1}\Psi(u_{k,\epsilon}^{1,-})e^{-\mathbb{Q}_2S^{\beta^*/2}}\nonumber\\
&=&o\left(\pi_{1,\epsilon}^{(1)}(u)\right),\ \ u\rw\IF, S\rw\IF,
\EQN
with $\hat{u}_{k,k_1,\epsilon}=\min(u_{k,\epsilon}^{1,-}, u_{k_1,\epsilon}^{1,-})$ and $\beta^*=\min (\alpha_1,\alpha_2)$. 
Using Lemma \ref{PIPI} yields that
\BQN\label{lambda21}
 \Lambda_2^{(1)}(u)&\leq&\sum_{k=-N_1(u)-2}^{N_1(u)+1}\Bigg[
 \pk{\sup_{(s,t)\in
 \widetilde{I}_{k,0}(u)}\frac{\overline{X}(s,t)}{1+(1-\epsilon)v^2_2(t)}>u_{k,\epsilon}^{1,-}} \\
&& +\pk{ \sup_{(s,t)\in
 \widetilde{I}_{k+1,0}(u)}\frac{\overline{X}(s,t)}{1+(1-\epsilon)v^2_2(t)}>u_{k+1,\epsilon}^{1,-}} \nonumber\\
 && -\pk{\sup_{(s,t)\in \widetilde{I}_{k,0}(u)\cup
 \widetilde{I}_{k+1,0}(u)}\frac{\overline{X}(s,t)}{1+(1+\epsilon)v^2_2(t)}>\check{u}_{k,k+1,\epsilon}} \Bigg]\nonumber\\
 &\leq&
 \left((1+\epsilon)2\mathcal{H}_{B_{\alpha_1}}[0,S]-(1-\epsilon)\mathcal{H}_{B_{\alpha_1}}[0,2S]\right)\nonumber\\
 &\ \ &\times\widehat{\mathcal{P}}^{\gamma_2(1-\epsilon)}_{\alpha_2}(S)
 \sum_{k=-N_1(u)-2}^{N_1(u)+1}\Psi(u_{k,\epsilon}^{1,-})\nonumber\\
& =&o\left(\pi_{1,\epsilon}^{(1)}(u)\right), \ \ u\rw\IF, S\rw\IF,
\EQN
 with $\check{u}_{k,k_1,\epsilon}=\max(u_{k,\epsilon}^{1,-}, u_{k_1,\epsilon}^{1,-})$.
 Combining (\ref{main1}), (\ref{PG1}), (\ref{Addition2}), (\ref{PG2}), (\ref{lambda11}) with (\ref{lambda21}) leads to
 $$\pi_1(u) \sim
 2\Gamma(1/\beta_1+1)\mathcal{H}_{\alpha_1}\widehat{\mathcal{P}}^{\gamma_2}_{\alpha_2}\frac{\overleftarrow{\vv}_1(u^{-1})}{\overleftarrow{\LL}_1(u^{-1})
}\Psi(u),  \ \ u\rw \IF.$$
The proof is \HEH{completed} by inserting the above asymptotic into (\ref{F1}).\\
\underline{Case ii) $\gamma_1=0, \gamma_2=\IF$.} In this case (\ref{main1}), (\ref{PG1}), (\ref{Addition2}), (\ref{PG2}), (\ref{lambda11}) and
(\ref{lambda21}) still hold  except for the fact that in light of Lemma {\ref{PIPI}}, we have to replace both
$\widehat{\mathcal{P}}^{\gamma_2(1-\epsilon)}_{\alpha_2}(S)$ and $\widehat{\mathcal{P}}^{\gamma_2}_{\alpha_2}$ by
1 in (\ref{PG1}) and (\ref{Addition2}). \\
\underline{Case iii) $\gamma_1, \gamma_2\in (0,\IF)$.} Let $\widehat{I}_{0,0}(u)=I_{0,0}(u)\cup I_{-1,0}(u)\cup I_{0,-1}(u)\cup I_{-1,-1}(u)$. It follows
straightforwardly that for any $0<\epsilon<1/2$ and $u$ large enough
\BQN\label{main2}
\pi_{1,\epsilon}^{(3)}(u)\leq \pi_1(u)\leq \pi_{1,-\epsilon}^{(3)}(u)+\pi_{1,-\epsilon}^{(4)}(u)
\EQN
with
$$
\pi_{1,\pm\epsilon}^{(3)}(u)=\mathbb{P}\left(\sup_{(s,t)\in
\widehat{I}_{0,0}(u)}\frac{\overline{X}(s,t)}{(1+(1\pm\epsilon)v^2_1(s))(1+(1\pm\epsilon)v^2_2(t))}>u\right),$$
and
$$\pi_{1,-\epsilon}^{(4)}(u)=\sum_{|k|=1, k\neq -1}^{N_1(u)+2}\sum_{|l|=1, l\neq -1}^{N_2(u)+2}\mathbb{P}\left(\sup_{(s,t)\in
I_{k,l}(u)}\overline{X}(s,t)>u_{k,l,\epsilon}^-\right).
$$
By Lemma {\ref{PIPI}}, it follows that
\BQN\label{PG3}
\pi_{1,\pm\epsilon}^{(3)}(u)\sim {\prod_{i=1}^2} \widehat{\mathcal{P}}^{(1\pm\epsilon)\gamma_i}_{\alpha_i}(S)\Psi(u), \ \ u\rw\IF, S\rw\IF.
\EQN
In addition,
using Lemma \ref{PIPI} and (\ref{A0}),
the same argument as given in the derivation of
the upper bound for $\pi_{1,-\epsilon}^{(2)}(u)$ yields
\BQN\label{PG4}
\pi_{1,-\epsilon}^{(4)}(u)
=o(\pi_{1, \pm\epsilon}^{(3)}(u))
\EQN

as $u\rw\IF$ and $S\rw\IF$.
Combination of (\ref{main2}) and (\ref{PG3}) with (\ref{PG4})
leads to $$\pi_1(u)\sim \prod_{i=1}^2 \widehat{\mathcal{P}}^{\gamma_i}_{\alpha_i}\Psi(u), \ \ u\rw\IF, S\rw\IF,$$
\eMM{hence the proof of this case is complete}.\\
\underline{Case iii) $\gamma_1\in (0,\IF), \gamma_2=\IF$.}
\KD{The proof follows the same lines} as given in previous case, 
with \KD{$\widehat{\mathcal{P}}^{\gamma_2}_{\alpha_2}$} replaced by 1.\\
\underline{Case iii) $\gamma_1, \gamma_2=\IF$.} Similarly,  (\ref{main2}), (\ref{PG3}) and (\ref{PG4}) hold with
\KD{$\widehat{\mathcal{P}}^{\gamma_1}_{\alpha_1},\widehat{\mathcal{P}}^{\gamma_2}_{\alpha_2}$ replaced by 1.}
\QED

\prooftheo{th.T4}
Similarly as in (\ref{main})
\BQN\label{E0}
\pi_{1}^+(u)-\Lambda^{(1)}(u)\leq \pi_1(u)\leq \pi_{1}^-(u),
\EQN
\KD{with}
\BQNY
\pi_{1}^{\pm}(u):&=&\sum_{k=-N_1(u)\pm 2}^{N_1(u)\mp 1}\sum_{l=-N_2(u)\pm 2}^{N_2(u)\mp 1}\mathbb{P}\left(\sup_{(s,t)\in
I_{k,l}(u)}\overline{X}(s,t)>u_{k,l,\epsilon}^{\pm}\right),\\
 \Lambda^{(1)}(u):&=&\sum_{(k,l,k_1,l_1)\in \mathbb{V}_1(u)\cup \mathbb{V}_2(u) }\mathbb{P}\left(\sup_{(s,t)\in I_{k,l}(u)}X(s,t)>u, \sup_{(s,t)\in
 I_{k_1,l_1}(u)}X(s,t)>u\right).
 \EQNY

Since \KD{$B$} is \HEH{non-singular matrix},
then there exists a positive constant $\mu>0$ such that for any $s, t$,
$$|b_{11}s+b_{12}t|+|b_{21}s+b_{22}t|\geq \mu \left(|s|+|t|\right).$$
Thus, for $(s,t)\in I_{k,l}(u)$ with $|k|\geq k_0\geq 2$ and $|l|\geq l_0\geq 2$
$$|b_{11}s+b_{12}t|+|b_{21}s+b_{22}t|\geq \mu S\left((k_0-1)\overleftarrow{\LL_1}(u^{-1})+(l_0-1)\overleftarrow{\LL_2}(u^{-1})\right).$$
\def\AS{\HEH{a(s,t)}}
By UCT,
for any $(s,t), (s',t')\in I_{k,l}(u)$ with $k_0\leq |k|\leq N_1(u)+2, l_0\leq |l|\leq N_2(u)+2$ and $u$ large enough \HEH{set $\AS:= v_1^2(|b_{11}s+b_{12}t|)+v_2^2(|b_{21}s+b_{22}t|)$)}
$$\frac{\AS }{v_1^2(|b_{11}s+b_{12}t|+|b_{21}s+b_{22}t|)}\geq
(1-\epsilon/3)\left(\nu^{\beta_1}+\theta(1-\nu)^{\beta_1}\right) $$
and
\BQNY\frac{ \NE{a(s',t') } 
	}{v_1^2(|b_{11}s+b_{12}t|+|b_{21}s+b_{22}t|)}\leq
\frac{(1-\epsilon/3)^2}{1-\epsilon}\left((\nu+\delta)^{\beta_1}+\theta(1-\nu+\delta)^{\beta_1}\right) ,
\EQNY
where
\BQNY \NE{\nu}=\frac{|b_{11}s+b_{12}t|}{|b_{11}s+b_{12}t|+|b_{21}s+b_{22}t|}\in [0,1],
\EQNY
and
$$0\leq \delta\leq
\frac{\left(|b_{11}|+|b_{12}|+|b_{21}|+|b_{22}|\right)\left(\overleftarrow{\LL_1}(u^{-1})+\overleftarrow{\LL_2}(u^{-1})\right)S}
{|b_{11}s+b_{12}t|+|b_{21}s+b_{22}t|}\leq
\frac{2\left(\sum_{i,j=1}^2|b_{ij}|\right)(1+\eta^{-1/\alpha_1})}{\mu\left((k_0-1)+(l_0-1)\eta^{-1/\alpha_1}\right)}\rw 0, $$
as $k_0, l_0\rw\IF.$
\KD{Using that for any $0<\epsilon<1/2$, when $k_0$ and $l_0$ are large enough
\BQNY
\left(\nu^{\beta_1}+\theta(1-\nu)^{\beta_1}\right)\geq (1-\epsilon/3)\left((\nu+\delta)^{\beta_1}+\theta(1-\nu+\delta)^{\beta_1}\right), \ \ \nu\in [0,1]
\EQNY
}
 for any
 $0<\epsilon<1/2$, there exists $k_\epsilon, l_\epsilon$ such that for any $k_\epsilon\leq |k|\leq N_1(u)+2, l_\epsilon\leq |l|\leq N_2(u)+2$  and $u$
 large enough
 \NE{$$a(s,t)\geq (1-\epsilon)a(s',t'), \quad (s,t), (s',t')\in I_{k,l}(u),$$}
 which is equivalent to
\BQN\label{Basic1}
\inf_{(s,t)\in I_{k,l}(u)}
\AS 
\geq (1-\epsilon)\sup_{(s,t)\in
I_{k,l}(u)}\AS 
.
\EQN
\underline{ Case i).}
Using (\ref{uniform}) and by (\ref{Basic1}), we have
\BQNY
\pi_1^{-}(u)
&\sim& \sum_{k=-N_1(u)- 2}^{N_1(u)+ 1}\sum_{l=-N_2(u)- 2}^{N_2(u)+
1}\prod_{i=1}^2\mathcal{H}_{\alpha_i}[0,S]\Psi(u_{k,l,\epsilon}^-)\nonumber\\
&\sim&\prod_{i=1}^2\mathcal{H}_{\alpha_i}[0,S]\Psi(u)\sum_{k=-N_1(u)-2}^{N_1(u)+ 1}\sum_{l=-N_2(u)-2}^{N_2(u)+ 1}
e^{-(1-\epsilon)u^2\inf_{(s,t)\in I_{k,l}(u)} \AS}\\
&\leq& \prod_{i=1}^2\mathcal{H}_{\alpha_i}[0,S]\Psi(u)
\left(R_3(u)+R_4(u)+\right.\nonumber\\
&& \ +\left.\frac{1}{\overleftarrow{\LL_1}(u^{-1})\overleftarrow{\LL_2}(u^{-1})S^2}\sum_{|k|=k_\epsilon}^{N_1(u)+2}\sum_{|l|=l_\epsilon}^{N_2(u)+ 2}
\int_{(s,t)\in I_{k,l}(u)}e^{-(1-\epsilon)^2u^2 \AS 
	}dsdt\right),
\EQNY
where
$$R_3(u)=\sum_{|k|\leq k_\epsilon}\sum_{l=-N_2(u)-2}^{N_2(u)+ 1}
e^{-(1-\epsilon)u^2\inf_{(s,t)\in I_{k,l}(u)} \AS 
	},$$
and
$$R_4(u)=\sum_{|l|\leq l_\epsilon}\sum_{k=-N_1(u)-2}^{N_1(u)+ 1}
e^{-(1-\epsilon)u^2\inf_{(s,t)\in I_{k,l}(u)}\AS 
	}.$$
\KD{By linear transformation
$(s',t')^{\top}=B(s,t)^\top$
and Lemma \ref{integ}, we have}
\BQNY
&&\frac{1}{\overleftarrow{\LL_1}(u^{-1})\overleftarrow{\LL_2}(u^{-1})}\sum_{|k|=k_\epsilon}^{N_1(u)+2}\sum_{|l|=l_\epsilon}^{N_2(u)+ 2}
\int_{(s,t)\in I_{k,l}(u)}e^{-(1-\epsilon)^2u^2 \AS
	}dsdt\\
&& \ \ \ \leq \frac{1}{\overleftarrow{\LL_1}(u^{-1})\overleftarrow{\LL_2}(u^{-1})}\int_{-2\overleftarrow{v_1}(u^{-1}\ln
u)}^{2\overleftarrow{v_1}(u^{-1}\ln u)}\int_{-2\overleftarrow{v_2}(u^{-1}\ln u)}^{2\overleftarrow{v_2}(u^{-1}\ln
u)}e^{-(1-\epsilon)^2u^2 \AS 
}dsdt\\
&& \ \ \ \leq\frac{1}{\overleftarrow{\LL_1}(u^{-1})\overleftarrow{\LL_2}(u^{-1})}\frac{1}{\abs{ det(B)}}\int_{-\mathbb{Q}\overleftarrow{v_1}(u^{-1}\ln
u)}^{\mathbb{Q}\overleftarrow{v_1}(u^{-1}\ln u)}\int_{-\mathbb{Q}\overleftarrow{v_2}(u^{-1}\ln u)}^{\mathbb{Q}\overleftarrow{v_2}(u^{-1}\ln
u)}e^{-(1-\epsilon)^2u^2v_1^2(|s'|)}e^{-(1-\epsilon)^2u^2v_2^2(|t'|)}ds'dt'\\
&& \ \ \ =\frac{1}{\overleftarrow{\LL_1}(u^{-1})\overleftarrow{\LL_2}(u^{-1})}\frac{4}{\abs{ det(B)}}\int_{0}^{\mathbb{Q}\overleftarrow{v_1}(u^{-1}\ln
u)}\int_{0}^{\mathbb{Q}\overleftarrow{v_2}(u^{-1}\ln u)}e^{-(1-\epsilon)^2u^2v_1^2(|s'|)}e^{-(1-\epsilon)^2u^2v_2^2(|t'|)}ds'dt'\\
&& \ \ \ \sim (1-\epsilon)^{-1/\beta_1-1/\beta_2}\frac{4}{\abs{
det(B)}}\Gamma(1/\beta_1+1)\Gamma(1/\beta_2+1)\frac{\overleftarrow{\vv_1}(u^{-1})\overleftarrow{\vv_2}(u^{-1})}
{\overleftarrow{\LL_1}(u^{-1})\overleftarrow{\LL_2}(u^{-1})}\rw \IF, \ \ u\rw\IF.
\EQNY
Moreover, in light of Lemma \ref{simple1}, there exists a constant $\kappa_1>0$ such that
\BQN\label{eq5}
\kappa_1v^2_1(|s|)+\kappa_1v^2_2(|t|)\leq  \AS 
, \ \ s,t\in \mathbb{R}.
\EQN
Thus we have
\BQNY
R_3(u)&\leq& \sum_{|k|\leq k_\epsilon}\sum_{l=-N_2(u)-2}^{N_2(u)+ 1}
e^{-(1-\epsilon)u^2\inf_{(s,t)\in I_{k,l}(u)}\kappa_1\left(v_1^2(|s|)+v_2^2(|t|)\right)}\\
&\leq& (2k_\epsilon+1)\sum_{l=-N_2(u)-2}^{N_2(u)+ 1}
e^{-(1-\epsilon)u^2\inf_{t\in J_l(u)}\kappa_1v_2^2(|t|)}\\
&\leq& \mathbb{Q}\frac{\overleftarrow{\vv_2}(u^{-1})}
{\overleftarrow{\LL_2}(u^{-1})}=o\left(\frac{\overleftarrow{\vv_1}(u^{-1})\overleftarrow{\vv_2}(u^{-1})}
{\overleftarrow{\LL_1}(u^{-1})\overleftarrow{\LL_2}(u^{-1})}\right), \ \ u\rw\IF.
\EQNY
Similarly,
$$R_4(u)\leq \mathbb{Q}_1\frac{\overleftarrow{\vv_1}(u^{-1})}
{\overleftarrow{\LL_1}(u^{-1})}=o\left(\frac{\overleftarrow{\vv_1}(u^{-1})\overleftarrow{\vv_2}(u^{-1})}
{\overleftarrow{\LL_1}(u^{-1})\overleftarrow{\LL_2}(u^{-1})}\right), \ \ u\rw\IF.$$
Therefore,
\BQN\label{E1}
\pi_1^-(u)\leq \frac{4}{\abs{ det(B)}}\prod_{i=1}^2 \Bigl[\Gamma(1/\beta_i+1)\mathcal{H}_{\alpha_i}
\frac{\overleftarrow{\vv}_i(1/u)}{\overleftarrow{\LL}_i(1/u)}\Bigr]\Psi(u)(1+o(1)), \ \ u\rw\IF, S\rw\IF,\epsilon\rw 0.
\EQN
\KD{In the same way we obtain that}
\BQN\label{E2}
\pi_1^{+}(u)\geq \frac{4}{\abs{ det(B)}}\prod_{i=1}^2 \Bigl[\Gamma(1/\beta_i+1)\mathcal{H}_{\alpha_i}
\frac{\overleftarrow{\vv}_i(1/u)}{\overleftarrow{\LL}_i(1/u)}\Bigr]\Psi(u)(1+o(1)), \ \ u\rw\IF, S\rw\IF.
\EQN
\KD{Due to} (\ref{eq5}), letting
\BQN\label{Y}
Y(s,t)=\frac{\overline{X}(s,t)}{1+\frac{\kappa_1}{2}v_1^2(|s|)+\frac{\kappa_1}{2}v_2^2(|t|)}, \ \ (s,t)\in \mathbb{R}^2,
\EQN we have
\BQNY
 \Lambda^{(1)}(u)&\leq &\sum_{(k,l,k_1,l_1)\in \mathbb{V}_1(u)\cup \mathbb{V}_2(u) }\mathbb{P}\left(\sup_{(s,t)\in I_{k,l}(u)}Y(s,t)>u, \sup_{(s,t)\in
 I_{k_1,l_1}(u)}Y(s,t)>u\right).
\EQNY
\KD{The same argument as given in the proof of Theorem \ref{th.T1} leads to}
\BQN\label{E3}
\Lambda^{(1)}(u)=o(\pi_1^-(u)), \ \ u\rw\IF, S\rw\IF.
\EQN
Inserting (\ref{E1})-(\ref{E3}) into (\ref{E0}) yields
$$\pi_1(u)\sim \frac{4}{\abs{ det(B)}}\prod_{i=1}^2 \Bigl[\Gamma(1/\beta_i+1)\mathcal{H}_{\alpha_i}
\frac{\overleftarrow{\vv}_i(1/u)}{\overleftarrow{\LL}_i(1/u)}\Bigr]\Psi(u),$$
which together with (\ref{F1})  completes the proof.\\
\underline{Case ii) $\gamma_1,\gamma_2\in (0,\IF)$.} Using the same notation for $\widehat{I}_{0,0}(u)$ as that in the proof of \KD{Theorem \ref{th.T1}} for case iii)
$\gamma_1, \gamma_2\in (0,\IF)$, (\ref{main2}) holds with
$$
\pi_{1,\pm\epsilon}^{(3)}(u)=\mathbb{P}\left(\sup_{(s,t)\in
\widehat{I}_{0,0}(u)}\frac{\overline{X}(s,t)}{1+(1\pm\epsilon) \AS 
}>u\right),$$
and
$$\pi_{1,-\epsilon}^{(4)}(u)=\sum_{|k|=1, k\neq -1}^{N_1(u)+2}\sum_{|l|=1, l\neq -1}^{N_2(u)+2}\mathbb{P}\left(\sup_{(s,t)\in I_{k,l}(u)}X(s,t)>u\right).
$$
Noting that
\BQNY
&&u^2\left(v^2_1(|b_{11}\overleftarrow{\LL}_1(1/u)s+b_{12}\overleftarrow{\LL}_2(1/u)t|)
+v^2_2(|b_{21}\overleftarrow{\LL}_1(1/u)s+b_{22}\overleftarrow{\LL}_2(1/u)t|)\right)\\ \ \ &&\rw
\gamma_1|b_{11}s+b_{12}\eta^{-1/\alpha_1}t|^{\alpha_1}+\theta\gamma_1|b_{21}s+b_{22}\eta^{-1/\alpha_1}t|^{\alpha_1}, \ \ u\rw\IF
\EQNY
uniformly with respect to $s,t\in [-S,S]^2$,
it follows from Lemma \ref{PIPI} that
\BQN
\pi_{1,\pm\epsilon}^{(3)}(u)\sim\widehat{\mathcal{P}}_\alpha^{(1\pm \epsilon)\gamma_1, (1\pm\epsilon)\theta\gamma_1, B_{\eta,\alpha}}(S)\Psi(u)\sim
\widehat{\mathcal{P}}_\alpha^{\gamma_1,\theta\gamma_1, B_{\eta,\alpha}}\Psi(u), \ \ u\rw\IF, S\rw\IF, \epsilon\rw 0.
\EQN
Moreover, by Lemma \ref{simple1} and (\ref{PG4}), with $Y$ defined by (\ref{Y}),
\BQNY
\pi_{1,-\epsilon}^{(4)}(u)\leq \sum_{|k|=1, k\neq -1}^{N_1(u)+2}\sum_{|l|=1, l\neq -1}^{N_2(u)+2}\mathbb{P}\left(\sup_{(s,t)\in
I_{k,l}(u)}Y(s,t)>u\right)\KD{=o\left(\Psi(u)\right)},
\EQNY
\KD{as $u\to\infty, S\rw\IF$}.
Thus
$\pi_1(u)\sim\widehat{\mathcal{P}}_\alpha^{\gamma_1,\theta\gamma_1, B_{\eta,\alpha}}\Psi(u),$
which completes the proof.\\
\underline{Case ii) $\gamma_1=\gamma_2=\IF$.}
\KD{The proof follows by the same argument as the proof of Case ii) $\gamma_1,\gamma_2\in (0,\IF)$,
with  $\widehat{\mathcal{P}}_\alpha^{\gamma_1,\theta\gamma_1, B_{\eta,\alpha}}$
replaced by 1}.
\QED

\prooftheo{th.T2}
\KD{This scenario requires a modification of set $D_u$.}
Let
$D_u^{(1)}=\{(s,t), |s+b_{12}t|\leq \overleftarrow{\vv}_1(u^{-1}\ln u), |t|\leq \overleftarrow{\vv}_2(u^{-1}\ln u) \}$.
\eMM{It follows that (\ref{F1}})-(\ref{neww}) also hold with $D_u$ replaced by
$D_u^{(1)}$. In this scenario, denote $\alpha=\alpha_1=\alpha_2$.\\
 \underline{Case i) $\gamma_1=\gamma_2=0$.}
 \KD{Let }
$$E_{l}^+(u)=\{k: I_{k,l}(u)\subset D_u^{(1)}\}, E_{l}^-(u)=\{k:
 I_{k,l}(u)\cap D_u^{(1)}\neq \emptyset\},$$ $$E^{(1)}(u)=\{(k,l, k_1,l_1), k\leq k_1, I_{k,l}(u)\cap D_u^{(1)}\neq \emptyset, I_{k_1,l_1}(u)\cap D_u^{(1)}\neq
 \emptyset \ \ \text{and} \ \ I_{k,l}(u)\cap I_{k',l'}(u)= \emptyset\},$$
 $$E^{(2)}(u)=\{(k,l, k_1,l_1), k\leq k_1, I_{k,l}(u)\cap D_u^{(1)}\neq \emptyset, I_{k_1,l_1}(u)\cap D_u^{(1)}\neq \emptyset, (k,l)\neq (k_1,l_1)\ \  \text{and} \ \ I_{k,l}(u)\cap
 I_{k',l'}(u) \neq\emptyset\}.$$
It follows that
\BQN\label{E4}
\pi_{2}^+(u)-\sum_{i=1}^2\Lambda_i^{(2)}(u)\leq \pi_1(u)\leq \pi_{2}^-(u),
\EQN
where
\BQNY
\pi_{2}^{\pm}(u):&=&\sum_{l=-N_2(u)\pm 2}^{N_2(u)\mp 1}\sum_{k\in E^{\pm}_{l}(u)}\mathbb{P}\left(\sup_{(s,t)\in
I_{k,l}(u)}\overline{X}(s,t)>u_{k,l,\epsilon}^{\pm}\right),\\
 \Lambda_i^{(2)}(u):&=&\sum_{(k,l,k_1,l_1)\in E^{(i)}(u)}\mathbb{P}\left(\sup_{(s,t)\in I_{k,l}(u)}\overline{X}(s,t)>u_{k,l,\epsilon}^- , \sup_{(s,t)\in
 I_{k_1,l_1}(u)}\overline{X}(s,t)>u_{k_1,l_1,\epsilon}^-\right).
 \EQNY
Using (\ref{uniform}), we have
\BQNY
\pi_2^-(u)
&=&\sum_{l=-N_2(u)-2}^{N_2(u)+ 1}\sum_{k\in E^{-}_{0}(u)}\mathbb{P}\left(\sup_{(s,t)\in I_{0,0}(u)}\overline{X}(\overleftarrow{\LL_1}(u^{-1})kS+s,\overleftarrow{\LL_2}(u^{-1})lS+t)>u_{k,l,\epsilon}^-\right)\nonumber\\
&\sim& \sum_{l=-N_2(u)-2}^{N_2(u)+ 1}\sum_{k\in E^{-}_{l}(u)}\prod_{i=1}^2\mathcal{H}_{\alpha_i}[0,S]\Psi(u_{k,l,\epsilon}^-)\nonumber\\
&\sim&\prod_{i=1}^2\mathcal{H}_{\alpha_i}[0,S]\Psi(u)\sum_{l=-N_2(u)-2}^{N_2(u)+ 1}\sum_{k\in E^{-}_{l}(u)}
e^{-(1-\epsilon)u^2\inf_{(s,t)\in I_{k,l}(u)}\left(v_1^2(|s+b_{12}t|)+v_2^2(|t|)\right)}
\EQNY
We observe that, for $u$ sufficiently large and all $|l|\leq N_2(u)+2$,
$$E_{l}^{-}\subset  \Biggl\{k\in \mathbb{N}, |k-\left[b_{12}l\frac{\overleftarrow{\LL}_2(1/u)}{\overleftarrow{\LL}_1(1/u)}\right]|\leq
N_1(u)+2(2+\left[|b_{12}|\eta^{-1/\alpha}\right]) \Biggr\},$$
and
$$E_l^+\supset\Biggl \{k\in \mathbb{N}, |k-\left[b_{12}l\frac{\overleftarrow{\LL}_2(1/u)}{\overleftarrow{\LL}_1(1/u)}\right]|\leq
N_1(u)-2(2+\left[|b_{12}|\eta^{-1/\alpha}\right])\Biggr\}.$$
By  UCT, we have \KD{that} for any $\epsilon>0$ there exists $l_\epsilon>0$ such that
\BQNY
\inf _{t\in J_l(u)}v_2^2(|t|)\geq (1-\epsilon) \sup_{t\in J_l(u)}v_2^2(|t|)
\EQNY
holds for $l_\epsilon\leq |l|\leq N_2(u)+2$.
Moreover, for any $\epsilon>0$ there exists $k_\epsilon>0$ such that
\BQNY
\inf _{t\in I_{k,l}(u)}v_1^2(|s+b_{12}t|)\geq (1-\epsilon) \sup_{s\in I_k(u)}v_1^2(|s+b_{12}l\overleftarrow{\LL}_2(u^{-1})S|)
\EQNY
hold for $|l|\leq N_2(u)+2$ and
$$k\in E_{l,\epsilon}^-(u):= \Bigg\{k, k_\epsilon\leq
|k-\left[b_{12}l\frac{\overleftarrow{\LL}_2(1/u)}{\overleftarrow{\LL}_1(1/u)}\right]|\leq N_1(u)+2(1+\left[|b_{12}|\eta^{-1/\alpha}\right]) \Bigg\}.$$ Therefore, in light of Lemma \ref{integ}, we have
\BQNY
&&\sum_{l=-N_2(u)-2}^{N_2(u)+ 1}\sum_{k\in E^{-}_{l}(u)}
e^{-(1-\epsilon)u^2\inf_{(s,t)\in I_{k,l}(u)}\left(v_1^2(|s+b_{12}t|)+v_2^2(|t|)\right)}\\
&& \ \ \  \leq \sum_{l=-N_2(u)-2}^{N_2(u)+ 1}e^{-(1-\epsilon)u^2\inf_{t\in J_l(u)}v_2^2(|t|)}\sum_{k\in E^{-}_{l}(u)}
e^{-(1-\epsilon)u^2\inf_{(s,t)\in I_{k,l}(u)}v_1^2(|s+b_{12}t|)}\\
&& \ \ \ \leq \frac{1}{\overleftarrow{\LL}_2(u^{-1})S}\sum_{|l|\geq l_\epsilon}^{N_2(u)+ 2}\int_{t\in
J_l(u)}e^{-(1-\epsilon)^2u^2v_2^2(|t|)}dt \\
&& \ \ \times \left(2k_\epsilon+1+\frac{1}{\overleftarrow{\LL}_1(1/u)S} \sum_{k\in E_{l,\epsilon}(u)}
\int_{s\in I_k(u)}e^{-(1-\epsilon)^2u^2v_1^2(|s+b_{12}l\overleftarrow{\LL}_2(1/u)S|)}ds\right)\\
&&\ \  + \sum_{|l|=0}^{l_\epsilon}\Bigg(2k_\epsilon+1+\frac{1}{\overleftarrow{\LL}_1(1/u)S}\sum_{k\in E_{l,\epsilon}(u)}
\int_{s\in I_k(u)}e^{-(1-\epsilon)^2u^2v_1^2(|s+b_{12}l\overleftarrow{\LL}_2(1/u)S|)}ds\Bigg)\\
&& \ \ \ \leq \frac{2(1+o(1))}{\overleftarrow{\LL}_2(u^{-1})S} \int_{0}^{\mathbb{Q}\overleftarrow{v_2}(u^{-1}\ln
u)}e^{-(1-\epsilon)^2u^2v_2^2(|t|)}dt\Bigg (2k_\epsilon+1+\frac{2}{\overleftarrow{\LL}_1(1/u)S}\int_{0}^{\mathbb{Q}\overleftarrow{v_1}(u^{-1}\ln u)}
e^{-(1-\epsilon)^2u^2v_1^2(|s|)}ds\Bigg)\\
&& \ \ \ \sim (1-\epsilon)^{-1/\beta_1-1/\beta_2}\frac{4}{S^2}\prod_{i=1}^2
\Gamma(1/\beta_i+1)\frac{\overleftarrow{\vv}_i(1/u)}{\overleftarrow{\LL}_i(1/u)}, \ \ u\rw\IF.
\EQNY
Consequently,
\KD{
\BQN\label{E5}
\pi_2^-(u)&\leq&
4\prod_{i=1}^2\left[\Gamma(1/\beta_i+1)\mathcal{H}_{\alpha_i}\frac{\overleftarrow{\vv}_i(1/u)}{\overleftarrow{\LL}_i(1/u)}\right]
\Psi(u)(1+o(1)), \ \ u\rw\IF, S\rw\IF, \epsilon\rw 0.
\EQN
Let $E_{l,\epsilon}^+(u):=\{k, k_\epsilon\leq |k-\left[bl\frac{\overleftarrow{\LL}_2(1/u)}{\overleftarrow{\LL}_1(1/u)}\right]|\leq
N_1(u)-2(1+\left[|b|\eta^{-1/\alpha}\right])\}$.
}
Similarly,
\BQN\label{E6}
\pi_2^+(u)&\geq& \prod_{i=1}^2\mathcal{H}_{\alpha_i}[0,S]\Psi(u)\sum_{l=-N_2(u)+2}^{N_2(u)-1}\sum_{k\in E^{+}_{l}(u)}
e^{-(1-\epsilon)u^2\sup_{(s,t)\in I_{k,l}(u)}\left(v_1^2(|s+bt|)+v_2^2(|t|)\right)}\nonumber\\
&\geq& \prod_{i=1}^2\mathcal{H}_{\alpha_i}[0,S]\Psi(u)\sum_{|l|=k_\epsilon}^{N_2(u)-2}e^{-(1-\epsilon)u^2\sup_{t\in
J_l(u)}v_2^2(|t|)}\sum_{k\in E^{+}_{l,\epsilon}(u)}
e^{-(1-\epsilon)u^2\sup_{(s,t)\in I_{k,l}(u)}v_1^2(|s+bt|)}\nonumber\\
&\sim& 4\prod_{i=1}^2\left[\Gamma(1/\beta_i+1)\mathcal{H}_{\alpha_i}\frac{\overleftarrow{\vv}_i(1/u)}{\overleftarrow{\LL}_i(1/u)}\right]
\Psi(u), \ \ u\rw\IF, S\rw\IF, \epsilon\rw 0.
\EQN
\KD{
Following the same argumentation as given
in
(\ref{lambda1}) and (\ref{lambda2}),
we get that
$\Lambda^{(2)}_i(u)
=o\left(\pi_2^-(u)\right), \ \  i=1,2, u\rw\IF,$
which together with (\ref{E5}) and (\ref{E6}) completes the proof.}\\
\underline{Case ii) $\gamma_2=0, \gamma_1\in (0,\IF)$.}
 We first
\KD{introduce}
\begin{eqnarray*}
L_{0,l}^*(u)&=&\{(s,t), |s+b_{12}t|\leq \overleftarrow{\LL}_1(1/u)S, t\in [l\overleftarrow{\LL}_2(1/u)S, (l+1)\overleftarrow{\LL}_2(1/u)S]\},\\
L_{k,l}(u)&=&\{(s,t), k\overleftarrow{\LL}_1(1/u)S\leq s+b_{12}t\leq (k+1)\overleftarrow{\LL}_1(1/u)S, t\in [l\overleftarrow{\LL}_2(1/u)S,
(l+1)\overleftarrow{\LL}_2(1/u)S]\},
\end{eqnarray*}
$$u_{k,l,\epsilon,*}^-=u\left(1+(1-\epsilon)\inf_{(s,t)\in L_{k,l}(u)}\left(v_1^2(|s+b_{12}t|)+v_2^2(|t|)\right)\right)$$ 
with $k,l \in \mathbb{Z}$.
Then we have
\BQN
\pi_3^+(u)+\sum_{i=1}^2\Lambda_i^{(3)}(u)\leq \pi_1(u)\leq \pi_3^-(u)+\pi_4(u),
\EQN
where
\BQNY
\pi_3^{\pm}(u)&:=&\sum_{l=-N_2(u)\pm 2}^{N_2(u)\mp 1}\mathbb{P}\left(\sup_{(s,t)\in L_{0,l}^*(u)}\frac{\overline{X}(s,t)}{1+(1\pm
\epsilon)v_1^2(|s+b_{12}t|)}>u_{l,\epsilon}^{2,\pm}\right),\\
\pi_4(u)&:=& \sum_{|k|\leq N_1(u)+2, k\neq 0, -1} \sum_{l=-N_2(u)\pm 2}^{N_2(u)\mp 1}\mathbb{P}\left(\sup_{(s,t)\in
L_{k,l}(u)}\overline{X}(s,t)>u_{k,l,\epsilon,*}^-\right),\\
 \Lambda_1^{(3)}(u)&:=&\sum_{-N_2(u)-2\leq l+1<l_1\leq N_2(u)+2}\mathbb{P}\left(\sup_{(s,t)\in L_{0,l}(u)}\overline{X}(s,t)>u_{l,\epsilon}^{2,-}, \sup_{(s,t)\in
L_{0,l_1}(u)}\overline{X}(s,t)>u_{l_1,\epsilon}^{2,-}\right)\\
 \Lambda_2^{(3)}(u)&:=&\sum_{l=-N_2(u)-2}^{N_2(u)+2}\mathbb{P}\left(\sup_{(s,t)\in L_{0,l}(u)}\frac{\overline{X}(s,t)}{1+(1-
 \epsilon)v_1^2(|s+b_{12}t|)}>u_{l,\epsilon}^{2,-}, \sup_{(s,t)\in L_{0,l+1}(u)}\frac{\overline{X}(s,t)}{1+(1- \epsilon)v_1^2(|s+b_{12}t|)}>u_{l+1,\epsilon}^{2,-}\right).
 \EQNY
 Let
 $$X_l(s,t)=\overline{X}(-b_{12}l\overleftarrow{\rho}_2(u^{-1})S+s,l\overleftarrow{\rho}_2(u^{-1})S+t),\quad  \mathcal{K}_u=\{l, |l|\leq N_2(u)+2\}, \quad  \mathcal{E}_u=L_{0,0}^*(u), $$ $$ h_l(u)=u_{l,\epsilon}^{2,-}, \quad d_u(s,t)=(1-\epsilon)v_1^2(|s+b_{12}t|)$$
 Since
 $$\lim_{u\rw\IF}\sup_{l\in \mathcal{K}_u}\left|(u_{l,\epsilon}^{2,-})^2v_1^2(|\overleftarrow{\rho}_1(1/u)s+b_{12}\overleftarrow{\rho}_2(1/u)t|)-
 \gamma_1|s+b_{12}\eta^{-1/\alpha}t|\right|=0, $$
 uniformly over any compact set,
by  Lemma \ref{PIPI}, we have
$$\lim_{u\rw\IF}\sup_{l\in\mathcal{K}_u}\left|(\Psi(u_{l,\epsilon}^{2,-}))^{-1}\mathbb{P}\left(\sup_{(s,t)\in L_{0,0}^*(u)}\frac{\overline{X}(-b_{12}l\overleftarrow{\LL}_2(1/u)S+s,l\overleftarrow{\LL}_2(1/u)S+t)}{1+(1-
 \epsilon)v_1^2(|s+b_{12}t|)}>u_{l,\epsilon}^{2,-}\right)-\mathcal{H}_\alpha^{(1-\epsilon)\gamma_1, b_{12}\eta^{-1/\alpha}}(S)\right|=0.$$
Thus we have,
 \BQN\label{eq2}
 \pi_3^{-}(u)
 &=& \sum_{l=-N_2(u)- 2}^{N_2(u)+ 1}\mathbb{P}\left(\sup_{(s,t)\in L_{0,0}^*(u)}\frac{\overline{X}(-b_{12}l\overleftarrow{\LL}_2(1/u)S+s,l\overleftarrow{\LL}_2(1/u)S+t)}{1+(1-
 \epsilon)v_1^2(|s+b_{12}t|)}>u_{l,\epsilon}^{2,-}\right)\nonumber\\
&\sim& \sum_{l=-N_2(u)- 2}^{N_2(u)+ 1}\Psi(u_{l,\epsilon}^{2,-})\mathcal{H}_\alpha^{(1-\epsilon)\gamma_1, b_{12}\eta^{-1/\alpha}}(S)\nonumber\\
 &\sim&\mathcal{H}_\alpha^{(1-\epsilon)\gamma_1, b_{12}\eta^{-1/\alpha}}(S) \Psi(u)\sum_{l=-N_2(u)- 2}^{N_2(u)+
1}e^{-(1-\epsilon)^2u^2\inf_{t\in J_l(u)}v_2^2(|t|)}\nonumber\\
 &\sim& \frac{\mathcal{H}_\alpha^{\gamma_1, b_{12}\eta^{-1/\alpha}}(S)}{S}2\Gamma(1/\beta_2+1)\Psi(u)\frac{\overleftarrow{v}_2(1/u)}{\overleftarrow{\rho}_2(1/u)}(1+o(1)), \
 \ u\rw\IF, \upsilon, \epsilon\rw 0.
 \EQN
Moreover, in light of \cite{nonhomoANN},
$$\lim_{S\rw\IF} \frac{\mathcal{H}_\alpha^{\gamma_1, b_{12}\eta^{-1/\alpha}}(S)}{S}=\mathcal{H}_\alpha^{\gamma_1, b_{12}\eta^{-1/\alpha}}\in (0,\IF).$$
Thus we have
\BQN\label{E7}
 \pi_3^{-}(u)\leq \mathcal{H}_\alpha^{\gamma_1, b_{12}\eta^{-1/\alpha}}2\Gamma(1/\beta_2+1)\Psi(u)\frac{\overleftarrow{v}_2(1/u)}{\overleftarrow{\rho}_2(1/u)}(1+o(1)), \ \ u\rw\IF, S\rw\IF, \eMM{\epsilon\rw 0} .
\EQN
Similarly,
\BQN\label{E12}
 \pi_3^{+}(u)\geq \mathcal{H}_\alpha^{\gamma_1, b_{12}\eta^{-1/\alpha}}2\Gamma(1/\beta_2+1)\Psi(u)\frac{\overleftarrow{v}_2(1/u)}{\overleftarrow{\rho}_2(1/u)}(1+o(1)), \ \ u\rw\IF,  S\rw\IF, \eMM{\epsilon\rw 0}.
\EQN
Note that for $u$ sufficiently large
\BQN\label{eq6}
L_{0,0}(u)\subset [-(1+2|b_{12}|\eta^{-1/\alpha})\overleftarrow{\rho}_1(1/u)S, (1+2|b_{12}|\eta^{-1/\alpha})\overleftarrow{\rho}_1(1/u)S]\times [0,
 \overleftarrow{\rho}_2(1/u)S]=:J_{0,0}(u).
 \EQN
Thus,
with $S_2=(1+2|b|\eta^{-1/\alpha})S$, by (\ref{uniform}) with $u_{k,l,\epsilon,*}^-$ instead of $u_{k,l,\epsilon}^-$ , we obtain
 \BQN\label{E8}
 \pi_4(u)&=&\sum_{|k|\leq N_1(u)+2, k\neq 0, -1} \sum_{l=-N_2(u)- 2}^{N_2(u)+ 1}\mathbb{P}\left(\sup_{(s,t)\in
 L_{0,0}(u)}\overline{X}(k\overleftarrow{\LL}_1(1/u)S-b_{12}l\overleftarrow{\LL}_2(1/u)S+s,l\overleftarrow{\LL}_2(1/u)S+t)>u_{k,l,\epsilon,*}^-\right)\nonumber\\
 &\leq& \sum_{|k|\leq N_1(u)+2, k\neq 0, -1} \sum_{l=-N_2(u)- 2}^{N_2(u)+ 1}\mathbb{P}\left(\sup_{(s,t)\in
 J_{0,0}(u)}\overline{X}(k\overleftarrow{\LL}_1(1/u)S-b_{12}l\overleftarrow{\LL}_2(1/u)S+s,l\overleftarrow{\LL}_2(1/u)S+t)>u_{k,l,\epsilon,*}^-\right)\nonumber\\
 &\sim& \sum_{|k|\leq N_1(u)+2, k\neq 0, -1} \sum_{l=-N_2(u)- 2}^{N_2(u)+
 1}\widehat{\mathcal{H}}_\alpha(S_2)\mathcal{H}_{\alpha}(S)\Psi(u)e^{-(1-\epsilon)u^2\inf_{(s,t)\in L_{k,l}(u)}\left(v_1^2(|s+bt|)+v_2^2(|t|)\right)}\nonumber\\
 &\leq& \widehat{\mathcal{H}}_\alpha(S_2)\mathcal{H}_{\alpha}(S)\Psi(u)\sum_{1\leq|k|\leq
 N_1(u)+2}e^{-\mathbb{Q}u^2v_1^2(\overleftarrow{\rho}_1(1/u)|k|S)} \sum_{l=-N_2(u)- 2}^{N_2(u)+ 1}e^{-(1-\epsilon)u^2\inf_{t\in J_l(u)}v_2^2(|t|)}\nonumber\\
 &\leq& 2\Gamma(1/\beta_2+1)(1-\epsilon)^{-2/\beta_2}\frac{\widehat{\mathcal{H}}_\alpha(S_2)}{S}\frac{\mathcal{H}_{\alpha}(S)}{S}
 \frac{\overleftarrow{v}_2(1/u)}{\overleftarrow{\rho}_2(1/u)}\Psi(u)S\sum_{1\leq|k|\leq N_1(u)+2}e^{-\mathbb{Q}_1|kS|^{\beta_1/2}}\nonumber\\
 &\leq&\mathbb{Q}_2
 \frac{\overleftarrow{v}_2(1/u)}{\overleftarrow{\rho}_2(1/u)}\Psi(u)e^{-\mathbb{Q}_3S^{\beta_1/2}}=o\left(\pi_3^-(u)\right), \ \ u\rw\IF, S\rw\IF.
 \EQN
\KD{
Following the same idea as given in
the proof of Theorem \ref{th.T1}, we get that
}
$
\Lambda_1^{(3)}(u)+\Lambda_2^{(3)}(u)
=
o\left(\pi_3^-(u)\right)
$, as $\to\infty$,
which completes the proof of this case.
\\
\underline{Case ii) $\gamma_2=0, \gamma_1=\IF$.} It follows straightforwardly that, for any $x>0$ and $u$ sufficiently large,
\BQN
\mathbb{P}\left(\sup_{|t|\leq \overleftarrow{\vv}_2(u^{-1}\ln u)}X(-b_{12}t,t)>u\right)\leq \pi_1(u)\leq \mathbb{P}\left(\sup_{(s,t)\in
D_u^{(1)}}\frac{\overline{X}(s,t)}{1+x\rho_1^2(|s+b_{12}t|)+v_2^2(|t|)}>u\right).
\EQN
Using that the Gaussian random field on the right hand side of the above
satisfies case $\gamma_2=0, \gamma_1=x\in (0,\IF)$, by (\ref{eq2}) and (\ref{E8}) we get for $S$ sufficiently large
\BQNY
\mathbb{P}\left(\sup_{(s,t)\in D_u^{(1)}}\frac{\overline{X}(s,t)}
       {1+x\rho_1^2(|s+b_{12}t|)+v_2^2(|t|)}>u\right)
&\leq&
       \frac{\mathcal{H}_\alpha^{x,b_{12}\eta^{-1/\alpha}}(S)2\Gamma(1/\beta_2+1)}{S}\Psi(u)\frac{\overleftarrow{v}_2(1/u)}{\overleftarrow{\rho}_2(1/u)}(1+o(1)).
\EQNY
It follows that for any $S$ positive
\BQNY\lim_{x\rw\IF}\mathcal{H}_\alpha^{x, b_{12}\eta^{-1/\alpha}}(S)&=&\lim_{x\rw\IF}\EE{\sup_{(s+ b_{12}\eta^{-1/\alpha}t,t)\in [-S,S]\times[0,S]}  e^{ W(s,t) -x|s+ b_{12}\eta^{-1/\alpha}t|^{\alpha}
}}\\
&=&\EE{\sup_{(s+ b_{12}\eta^{-1/\alpha}t,t)\in \{0\}\times[0,S]}  e^{ W(s,t)
}}\\&=&
\mathcal{H}_\alpha(\left(|b_{12}|^{\alpha}\eta^{-1}+1\right)^{1/\alpha}S).
\EQNY
Hence, as $u\rw\IF, x\rw\IF, S\rw\IF$
$$\mathbb{P}\left(\sup_{(s,t)\in D_u^{(1)}}
\frac{\overline{X}(s,t)}{1+x\rho_1^2(|s+b_{12}t|)+v_2^2(|t|)}>u\right)\leq
2\left(|b_{12}|^{\alpha}\eta^{-1}+1\right)^{1/\alpha}\mathcal{H}_\alpha\Gamma(1/\beta_2+1)\Psi(u)\frac{\overleftarrow{v}_2(1/u)}{\overleftarrow{\rho}_2(1/u)}(1+o(1)).$$
Further,  for the process $X(-b_{12}t,t)$, we have
\BQN\label{E9}
1-\sqrt{Var\left(X(-b_{12}t,t)\right)}&\sim& v_2^2(|t|), \ \ t\rw 0,\nonumber\\
 1-Corr\left(X(-b_{12}t,t), X(-b_{12}s,s)\right)&\sim& \left(|b_{12}|^{\alpha}\eta^{-1}+1\right)\rho_2^2(|t-s|), \ \ s,t\rw 0.
\EQN
Thus in light of Theorem \ref{mainT}, we have
$$
\KD{\mathbb{P}\left(\sup_{|t|\leq \overleftarrow{\vv}_2(u^{-1}\ln u)}X(-b_{12}t,t)>u\right)}
\sim
2\left(|b_{12}|^{\alpha}\eta^{-1}+1\right)^{1/\alpha}\mathcal{H}_\alpha\Gamma(1/\beta_2+1)\Psi(u)\frac{\overleftarrow{v}_2(1/u)}{\overleftarrow{\rho}_2(1/u)}.$$
Consequently,
$$\pi_1(u)\sim
2\left(|b_{12}|^{\alpha}\eta^{-1}+1\right)^{1/\alpha}\mathcal{H}_\alpha\Gamma(1/\beta_2+1)\Psi(u)\frac{\overleftarrow{v}_2(1/u)}{\overleftarrow{\rho}_2(1/u)},$$
which completes the proof.\\
\underline{Case iii) $\gamma_2\in (0,\IF), \gamma_1=\IF$.}
Let
$\widehat{I}_{0,0}^*(u)=\{(s,t), |s+b_{12}t|\leq \overleftarrow{\LL}_1(1/u)S, |t|\leq
\overleftarrow{\LL}_2(1/u)S\}.$
Then for $u$ sufficiently large, we have
\BQN\label{E10}
\mathbb{P}\left(\sup_{(s,t)\in D_u^{(1)}}X(-b_{12}t,t)>u\right)\leq \pi_1(u)\leq \pi_{1,-\epsilon}^{(5)}(u)+\pi_{1,-\epsilon}^{(6)}(u)
\EQN
with
$$
\pi_{1,-\epsilon}^{(5)}(u)=\mathbb{P}\left(\sup_{(s,t)\in
\widehat{I}_{0,0}^*(u)}\frac{\overline{X}(s,t)}{1+x\rho^2_1(|s+b_{12}t|)+(1-\epsilon)v^2_2(|t|)}>u\right),$$
and
$$\pi_{1,-\epsilon}^{(6)}(u)=\sum_{|k|=1, k\neq -1}^{N_1(u)+2}\sum_{|l|=1, l\neq -1}^{N_2(u)+2}\mathbb{P}\left(\sup_{(s,t)\in
L_{k,l}(u)}\overline{X}(s,t)>u_{k,l,\epsilon}^-\right).
$$
Since
$$u^2(x\rho^2_1(|\overleftarrow{\rho}_1(1/u)s+b_{12}\overleftarrow{\rho}_2(1/u)t|)+(1-\epsilon)v^2_2(|\overleftarrow{\rho}_2(1/u)t|))\rw x|s+b_{12}\eta^{-1/\alpha}t|^{\alpha}+(1-\epsilon)\gamma_2|t|^{\alpha}, \ \ u\rw\IF$$
uniformly on any compact set,
then, by Lemma \ref{PIPI},
$$\pi_{1,-\epsilon}^{(5)}(u)\sim \widehat{\mathcal{H}}_{\alpha}^{x,\gamma_2,b_{12}\eta^{-1/\alpha}}(S)\Psi(u), \ \ u\rw\IF, \epsilon\rw 0.
$$
\KD{Moreover, by the same argument as given in case ii), we have
$$\lim_{x\rw\IF}\widehat{\mathcal{H}}_{\alpha}^{x,\gamma_2,b_{12}\eta^{-1/\alpha}}(S)=\widehat{\mathcal{P}}_\alpha^
{\gamma_2\left(|b_{12}|^\alpha\eta^{-1}+1\right)^{-1}}\left(\left(|b|^\alpha\eta^{-1}+1\right)^{1/\alpha}S\right).$$
Then}
\BQNY
\pi_{1,-\epsilon}^{(5)}(u)\sim \widehat{\mathcal{P}}_\alpha^
{\gamma_2\left(|b_{12}|^\alpha\eta^{-1}+1\right)^{-1}}\Psi(u), \ \ u\rw\IF, x\rw\IF, \epsilon\rw 0, S\rw\IF.
\EQNY
\KD{Using that
$L_{0,0}(u)\subset J_{0,0}(u),$ with $J_{0,0}(u)$ defined by (\ref{eq6}),
and following the same steps as in (\ref{E8}), we get}
\BQNY
\pi_{1,-\epsilon}^{(6)}(u)
=o(\Psi(u)), \ \
 u\rw\IF, S\eMM{\rw \IF}.
\EQNY
Hence, from Theorem \ref{mainT} and (\ref{E9})
$$\mathbb{P}\left(\sup_{(s,t)\in D_u^{(1)}}X(-b_{12}t,t)>u\right)\sim \widehat{\mathcal{P}}_\alpha^
{\gamma_2\left(|b_{12}|^\alpha\eta^{-1}+1\right)^{-1}}\Psi(u), \ \ u\rw\IF,$$
which establishes the claim.\\
\underline{Case iv) $\gamma_2=\IF, \gamma_1=\IF$.}
Clearly, (\ref{E10}) holds with
$$
\pi_{1,-\epsilon}^{(5)}(u):=\mathbb{P}\left(\sup_{(s,t)\in \widehat{I}_{0,0}(u)}
\frac{\overline{X}(s,t)}{1+x\rho^2_1(|s+b_{12}t|)+y\rho^2_2(|t|)}>u\right), x,
y>0.$$
Moreover,
$$\pi_{1,-\epsilon}^{(5)}(u)\sim \widehat{\mathcal{H}}_{\alpha}^{x,y,b_{12}\eta^{-1/\alpha}}(S)\Psi(u)$$
and
$$\lim_{y\rw\IF}\lim_{x\rw\IF}\widehat{\mathcal{H}}_{\alpha}^{x,y,b_{12}\eta^{-1/\alpha}}(S)=\lim_{y\rw\IF}\widehat{\mathcal{P}}_\alpha^
{y\left(|b_{12}|^\alpha\eta^{-1}+1\right)^{-1}}\left(\left(|b_{12}|^\alpha\eta^{-1}+1\right)^{1/\alpha}S\right)=1.$$
Hence
$\pi_{1,-\epsilon}^{(5)}(u)\sim \Psi(u), u\rw\IF, x\rw\IF, y\rw\IF.$
\KD{The rest of the proof is the same as the case
$\gamma_2\in (0,\IF), \gamma_1=\IF$.}
 \QED
\prooftheo{th.T3}.
We focus on $\pi_1(u)$ as $u\rw\IF$. \\
\underline{ Case i)}
\KD{
The proof of this case follows line by line the same arguments as given in the proof
of Case i) of Theorem \ref{th.T6}.
}
\\
\underline{Case ii) $\gamma_1=0, \gamma_2\in (0,\IF)$.}
\KD{
First we introduce some new notation.}
Let
$$u_{k,\epsilon}^{{* -}}=1+(1 {-} 3\epsilon)\inf_{t\in I_k(u)}\left(v_1^2(|(1+b_{12}\mu)s|)+v_2^2(|\mu s|)\right),$$
%
and
$$\widehat{I}_{k,0}(u)=I_k(u)\times \left(J_{-1}(u)\cup J_0(u)\right), \ \
v(s,t)=v_1^2(|s+b_{12}t|)+v_2^2(|t|)-v_1^2(|(1+b_{12}\mu)s|)-v_2^2(|\mu s|), \ \
(s,t)\in D_u,
$$
where $\mu$ is defined right before Theorem \ref{th.T3}.
 For any $0<x<y<\frac{S}{2|b_{12}|}$ and $0<\epsilon<1/4$, we have
\BQN\label{EQ3}
\pi_5^+(u)+\Lambda(u)\leq\pi_1(u)\leq \pi_5^-(u)+\pi_6(u)+\pi_7(u)+\pi_8(u),
\EQN
where
\BQNY
\pi_5^\pm(u):&=&\sum_{k\in E_{x,y}^\pm(u)}\mathbb{P}\left(\sup_{(s,t)\in
\widehat{I}_{k,0}(u)}\frac{\overline{X}(s,t)}{1+(1\pm\epsilon)\left(v_1^2(|s+b_{12}t|)+v_2^2(|t|)\right)}>u\right),\\
\pi_6(u):&=&\sum_{k\in E_{0,x}(u)}\mathbb{P}\left(\sup_{(s,t)\in
\widehat{I}_{k,0}(u)}\frac{\overline{X}(s,t)}{1+(1-\epsilon)\left(v_1^2(|s+b_{12}t|)+v_2^2(|t|)\right)}>u\right),\\
\pi_7(u):&=&\sum_{k\in E_{y,\IF}(u)}\mathbb{P}\left(\sup_{(s,t)\in
\widehat{I}_{k,0}(u)}\frac{\overline{X}(s,t)}{1+(1-\epsilon)\left(v_1^2(|s+b_{12}t|)+v_2^2(|t|)\right)}>u\right),\\
\pi_8(u):&=&\sum_{|k|=0}^{N_1(u)+2}\sum_{|l|=1, l\neq -1}^{N_2(u)+2}\mathbb{P}\left(\sup_{(s,t)\in
I_{k,l}(u)}\frac{\overline{X}(s,t)}{1+(1-\epsilon)\left(v_1^2(|s+b_{12}t|)+v_2^2(|t|)\right)}>u\right),\\
\Lambda(u):&=&\sum_{k<k_1\in E_{x,y}^-(u)}\mathbb{P}\left(\sup_{(s,t)\in \widehat{I}_{k,0}(u)}X(s,t)>u, \sup_{(s,t)\in
\widehat{I}_{k_1,0}(u)}X(s,t)>u\right),
\EQNY
with
\begin{eqnarray*}
E_{0,x}
&=&
\{k, |k|\leq N_1(u)+2, I_{k}(u)\cap [-\overleftarrow{\LL_2}(u^{-1})x, \overleftarrow{\LL_2}(u^{-1})x]\neq \KD{\emptyset}\},\\
E_{x,y}^-
&=&
\{k, |k|\leq N_1(u)+2, I_{k}(u)\cap \left([-\overleftarrow{\LL_2}(u^{-1})y, -\overleftarrow{\LL_2}(u^{-1})x
]\cup[\overleftarrow{\LL_1}(u^{-1})x, \overleftarrow{\LL_1}(u^{-1})y ]\right)\neq \KD{\emptyset}\}, \\
E_{x,y}^+
&=&
\{k, |k|\leq N_1(u)+2, I_{k}(u)\subset \left([-\overleftarrow{\LL_2}(u^{-1})y, -\overleftarrow{\LL_2}(u^{-1})x
]\cup[\overleftarrow{\LL_1}(u^{-1})x, \overleftarrow{\LL_1}(u^{-1})y ]\right)\},\\
E_{y,\IF}^-
&=&
\{k, |k|\leq N_1(u)+2, I_{k}(u)\cap \left([-\IF, -\overleftarrow{\LL_2}(u^{-1})y]
\cup[\overleftarrow{\LL_1}(u^{-1})y, \IF]\right)\neq \KD{\emptyset}\}.
\end{eqnarray*}
We observe that for $|s|\in [\frac{i-1}{n}\overleftarrow{\LL_2}(u^{-1}), \frac{i+2}{n}\overleftarrow{\LL_2}(u^{-1})]$ with $x/2\leq \frac{i}{n}\leq 2y$ and  and
$|t|\in [0, \overleftarrow{\vv}_2(\ln u/u)]$
\BQN\label{EQ2}
 &&1+(1-\epsilon)\left(v_1^2(|s+b_{12}t|)+v_2^2(|t|)\right)\nonumber\\
 && \ \ \geq  \left[1+(1-3\epsilon)\left(v_1^2(|(1+b_{12}\mu)s|)+v_2^2(|\mu s|)\right)\right]\left[1+(1-3\epsilon)v(i\overleftarrow{\LL_2}(u^{-1})/n,t)\right],
\EQN
whose proofs is postponed in the Appendix.
Let
$$X_{u,k}(s,t)=\overline{X}(k\overleftarrow{\LL}_1(1/u)S+s,t), \quad \mathcal{K}_u=E_{i/n,(i+1)/n }^-, \quad \mathcal{E}_u=\widehat{I}_{0,0}(u),$$
$$d_u(s,t)= (1-3\epsilon)u^2v(i\overleftarrow{\LL_2}(u^{-1})/n, t), \quad h_k(u)=u_{k,\epsilon}^{*-}.$$
We note that
$$\lim_{u\rw\IF}\sup_{k\in \mathcal{K}_u, t\in[-S,S]}\left|(u_{k,\epsilon}^{*-})^2v(i\overleftarrow{\LL_2}(u^{-1})/n, \overleftarrow{\LL_2}(u^{-1})t)-g_{i/n}(t)\right|=0.$$
Thus in light of Lemma \ref{PIPI}, we have
$$\lim_{u\rw\IF}\sup_{x/2\leq \frac{i}{n}\leq 2y}\sup_{k\in E_{i/n,(i+1)/n }^-}\left|(\Psi(u_{k,\epsilon}^{*-}))^{-1}\mathbb{P}\left(\sup_{(s,t)\in
\widehat{I}_{0,0}(u)}\frac{\overline{X}(k\overleftarrow{\LL}_1(1/u)S+s,t)}{1+(1-3\epsilon)
v(i\overleftarrow{\LL_2}(u^{-1})/n,t)}>u_{k,\epsilon}^{*-}\right)-\mathcal{H}_{\alpha_1}(S)\widehat{\mathcal{P}}_{\beta}^{(1-3\epsilon)g_{i/n}}
(S)\right|=0.$$

Thus for $[nx]-1\leq i\leq [ny]$,
it follows that
\BQNY
&&\sum_{k\in E_{i/n, (i+1)/n}^-}\mathbb{P}\left(\sup_{(s,t)\in
\widehat{I}_{k,0}(u)}\frac{\overline{X}(s,t)}{1+(1-\epsilon)\left(v_1^2(|s+b_{12}t|)+v_2^2(|t|)\right)}>u\right)\\
&& \ \ \ \leq\sum_{k\in E_{i/n, (i+1)/n}^-}\mathbb{P}\left(\sup_{(s,t)\in
\widehat{I}_{k,0}(u)}\frac{\overline{X}(s,t)}{\left[1+(1-3\epsilon)\left(v_1^2(|(1+b_{12}\mu)s|)+v_2^2(|\mu
s|)\right)\right]\left[1+(1-3\epsilon)v(i\overleftarrow{\LL_2}(u^{-1})/n,t)\right]}>u\right)\nonumber\\
&& \ \ \ \leq\sum_{k\in E_{i/n, (i+1)/n}^-}\mathbb{P}\left(\sup_{(s,t)\in
\widehat{I}_{k,0}(u)}\frac{\overline{X}(s,t)}{1+(1-3\epsilon)v(i\overleftarrow{\LL_2}(u^{-1})/n,t)}>u_{k,\epsilon}^{*-}\right)\nonumber\\
&&\ \ \ =\sum_{k\in E_{i/n, (i+1)/n}^-}\mathbb{P}\left(\sup_{(s,t)\in
\widehat{I}_{0,0}(u)}\frac{\overline{X}(k\overleftarrow{\LL}_1(1/u)S+s,t)}{1+(1-3\epsilon)v(i\overleftarrow{\LL_2}(u^{-1})/n,t)}>u_{k,\epsilon}^{*-}\right)\nonumber\\
&& \ \ \ \sim \sum_{k\in E_{i/n, (i+1)/n}^-}\mathcal{H}_{\alpha_1}(S)\widehat{\mathcal{P}}_{\beta}^{(1-3\epsilon)g_{i/n}}
(S)\Psi(u_{k,\epsilon}^{*-})\nonumber\\
&& \ \ \ \sim\mathcal{H}_{\alpha_1}(S)\widehat{\mathcal{P}}_{\beta}^{(1-3\epsilon)g_{i/n}}
(S)\Psi(u)\sum_{k\in E_{i/n, (i+1)/n}^-}e^{-u^2(1-3\epsilon)\inf_{t\in I_k(u)}\left(v_1^2(|(1+b_{12}\mu)s|)+v_2^2(|\mu
s|)\right)}\nonumber\\
&& \ \ \ \leq\frac{\mathcal{H}_{\alpha_1}(S)}{S}\widehat{\mathcal{P}}_{\beta}^{(1-3\epsilon)g_{i/n}}
(S)\frac{\Psi(u)}{\overleftarrow{\LL_1}(u^{-1})}2\int_{i\overleftarrow{\LL_2}(u^{-1})/n}^{(i+1)
\overleftarrow{\LL_2}(u^{-1})/n}e^{-(1-4\epsilon)\frac{\gamma_2 \Mba }{\theta} u^2 \rho_2^2(|s|)}ds(1+o(1)),
\EQNY
as $u\rw\IF$, with $ \Mba $ defined right before Theorem {\ref{th.T3}}.
\KD{Using the same arguments as in} the proof of Lemma \ref{integ}, we have
\BQNY
\int_{i\overleftarrow{\LL_2}(u^{-1})/n}^{(i+1)
\overleftarrow{\LL_2}(u^{-1})/n}e^{-(1-4\epsilon)\frac{\gamma_2 \Mba }{\theta} u^2 \rho_2^2(|s|)}ds &\sim&
\frac{2}{\beta}\overleftarrow{\LL_2}(u^{-1})\int_{(i/n)^{\beta/2}}^{\left((i+1)/n\right)^{\beta/2}}t^{2/\beta-1}
e^{-(1-4\epsilon)\frac{\gamma_2 \Mba }{\theta}t^2}dt\\
&\sim& \overleftarrow{\LL_2}(u^{-1})\int_{i/n}^{(i+1)/n}
e^{-(1-4\epsilon)\frac{\gamma_2 \Mba }{\theta}t^{\beta}}dt.
\EQNY
Hence
\BQNY
&&\sum_{k\in E_{i/n, (i+1)/n}^-}\mathbb{P}\left(\sup_{(s,t)\in
\widehat{I}_{k,0}(u)}\frac{\overline{X}(s,t)}{1+(1-\epsilon)\left(v_1^2(|s+b_{12}t|)+v_2^2(|t|)\right)}>u\right)\nonumber\\
&& \ \ \ \leq 2\frac{\mathcal{H}_{\alpha_1}(S)}{S}\widehat{\mathcal{P}}_{\beta}^{(1-3\epsilon)g_{i/n}}
(S)\frac{\overleftarrow{\LL_2}(u^{-1})}{\overleftarrow{\LL_1}(u^{-1})}\Psi(u)\int_{i/n}^{(i+1)/n}
e^{-(1-4\epsilon)\frac{\gamma_2 \Mba }{\theta}t^{\beta}}dt(1+o(1))\nonumber\\
&& \ \ \ \leq 2\mathcal{H}_{\alpha_1}\widehat{\mathcal{P}}_{\beta}^{g_{i/n}}
\frac{\overleftarrow{\LL_2}(u^{-1})}{\overleftarrow{\LL_1}(u^{-1})}\Psi(u)\int_{i/n}^{(i+1)/n}
e^{-\frac{\gamma_2 \Mba }{\theta}t^{\beta}}dt(1+o(1)), \ \ u\rw\IF,  S\rw\IF, \eMM{\epsilon\rw 0}.
\EQNY
Further, by the continuity of $\widehat{\mathcal{P}}_{\beta}^{g_{s}} $ over $s\in [x/2, 2y]$, we have
\BQN\label{EQ4}
\pi_5^-(u)
&\leq&   2\mathcal{H}_{\alpha_1}
\frac{\overleftarrow{\LL_2}(u^{-1})}{\overleftarrow{\LL_1}(u^{-1})}\Psi(u)\sum_{i=[nx]-1}^{[ny]+1}\int_{i/n}^{(i+1)/n}
\widehat{\mathcal{P}}_{\beta}^{g_{i/n}}e^{-\frac{\gamma_2 \Mba }{\theta}t^{\beta}}dt(1+o(1))\nonumber\\
&\leq&2\mathcal{H}_{\alpha_1}
\frac{\overleftarrow{\LL_2}(u^{-1})}{\overleftarrow{\LL_1}(u^{-1})}\Psi(u)\int_x^y
\widehat{\mathcal{P}}_{\beta}^{g_{t}}e^{-\frac{\gamma_2 \Mba }{\theta}t^{\beta}}dt(1+o(1)), \ \ u\rw\IF,
 S\rw\IF, \eMM{\epsilon\rw 0}, n\rw\IF.
\EQN
Similarly,
\BQN\label{EQ5}
\pi_5^+(u)
\geq 2\mathcal{H}_{\alpha_1}
\frac{\overleftarrow{\LL_2}(u^{-1})}{\overleftarrow{\LL_1}(u^{-1})}\Psi(u)\int_x^y
\widehat{\mathcal{P}}_{\beta}^{g_{t}}e^{-\frac{\gamma_2 \Mba }{\theta}t^{\beta}}dt(1+o(1)), \ \ u\rw\IF,
S\rw\IF, \eMM{\epsilon\rw 0}, n\rw\IF.
\EQN
Next we focus on $\pi_6(u)$.
In light of (\ref{eq5}) and (\ref{Y}), we have
\BQNY
\pi_6(u)\leq\sum_{k\in E_{0,x}(u)}\mathbb{P}\left(\sup_{(s,t)\in \widehat{I}_{k,0}(u)}Y(s,t)>u\right).
\EQNY
Hence, following case ii) $\gamma_1=0, \gamma_2\in (0,\IF)$ in Theorem \ref{th.T1}, we have
\BQN\label{EQ6}
\pi_6(u)&\leq& 2\mathcal{H}_{\alpha_1}(S)\widehat{\mathcal{P}}^{\gamma_2\kappa_1/2}_{\beta}(S)
\frac{\Psi(u)}{\overleftarrow{\LL}_1(u^{-1})S}\int_0^{x\overleftarrow{\LL_2}(u^{-1})}e^{-\frac{\kappa_1}{2}u^2v_1^2(t)}dt(1+o(1))\nonumber\\
&\leq& 2\mathcal{H}_{\alpha_1}\widehat{\mathcal{P}}^{\gamma_2\kappa_1/2}_{\beta}
\frac{\overleftarrow{\LL}_2(u^{-1})}{\overleftarrow{\LL}_1(u^{-1})}\Psi(u)\int_0^{x}e^{-\frac{\kappa_1}{2}\frac{\gamma_2}{\theta}t^{\beta}}dt(1+o(1)), \ \
u\rw\IF, S\rw\IF.
\EQN
Similarly,
\BQN\label{EQ7}
\pi_7(u)
\leq 2\mathcal{H}_{\alpha_1}\widehat{\mathcal{P}}^{\gamma_2\kappa_1/2}_{\beta}
\frac{\overleftarrow{\LL}_2(u^{-1})}{\overleftarrow{\LL}_1(u^{-1})}\Psi(u)\int_y^{\IF}e^{-\frac{\kappa_1}{2}\frac{\gamma_2}{\theta}t^{\beta}}dt(1+o(1)), \
\ u\rw\IF, S\rw\IF.
\EQN
Moreover, by Lemma {\ref{simple1}}, (\ref{lambda11}) and (\ref{lambda21}),
\BQN\label{EQ8}
\Lambda(u)\leq \sum_{k<k_1\in E_{x,y}^-(u)}\mathbb{P}\left(\sup_{(s,t)\in \widehat{I}_{k,0}(u)}Y(s,t)>u, \sup_{(s,t)\in
\widehat{I}_{k_1,0}(u)}Y(s,t)>u\right)=o(\pi_5^+(u)), \ \ u\rw\IF, S\rw\IF.
\EQN
Moreover, it follows from Lemma \ref{simple1} and (\ref{PG2}) that
\BQN\label{EQ9}
\pi_8(u)\leq \sum_{|k|=0}^{N_1(u)+2}\sum_{|l|=1, l\neq -1}^{N_2(u)+2}\mathbb{P}\left(\sup_{(s,t)\in I_{k,l}(u)}Y(s,t)>u\right)=o(\pi_5^+(u)), \ \ u\rw\IF,
S\rw\IF.
\EQN
Inserting (\ref{EQ4})--(\ref{EQ9}) into (\ref{EQ3}), we have
$$\pi_1(u)\geq 2\mathcal{H}_{\alpha_1}
\frac{\overleftarrow{\LL_2}(u^{-1})}{\overleftarrow{\LL_1}(u^{-1})}\Psi(u)\int_x^y
\widehat{\mathcal{P}}_{\beta}^{g_{t}}e^{-\frac{\gamma_2 \Mba }{\theta}t^{\beta}}dt(1+o(1), \ \ u\rw\IF, S\rw\IF,$$
and
\BQN\label{EQ10}
\pi_1(u)&\leq& 2\mathcal{H}_{\alpha_1}
\frac{\overleftarrow{\LL_2}(u^{-1})}{\overleftarrow{\LL_1}(u^{-1})}\Psi(u)\left(\int_x^y
\widehat{\mathcal{P}}_{\alpha_2}^{g_{t}}e^{-\frac{\gamma_2 \Mba }{\theta}t^{\beta}}dt
+\widehat{\mathcal{P}}^{\gamma_2\kappa_1/2}_{\beta}\int_0^{x}e^{-\frac{\kappa_1}{2}\frac{\gamma_2}{\theta}t^{\beta}}dt\right.\nonumber\\
&& \ \ \ \ \left.+
\widehat{\mathcal{P}}^{\gamma_2\kappa_1/2}_{\beta}\int_y^{\IF}e^{-\frac{\kappa_1}{2}\frac{\gamma_2}{\theta}t^{\beta}}dt\right)(1+o(1)),
\EQN
as $u\rw\IF, S\rw\IF$.
Letting $x\rw 0$ and $y\rw\IF$ leads to
$$\pi_1(u)\sim 2\mathcal{H}_{\alpha_1}
\frac{\overleftarrow{\LL_2}(u^{-1})}{\overleftarrow{\LL_1}(u^{-1})}\Psi(u)\int_0^\IF
\widehat{\mathcal{P}}_{\beta}^{g_{t}}e^{-\frac{\gamma_2 \Mba }{\theta}t^{\beta}}dt(1+o(1), \ \ u\rw\IF, S\rw\IF.$$
which, together with the fact that
$\overleftarrow{\LL_2}(u^{-1})\sim \left(\frac{\gamma_2}{\theta}\right)^{1/\beta}\overleftarrow{\vv_1}(u^{-1}),$
derives the claim.
This completes the proof.\\
\underline{Case iii) $\gamma_1=0, \gamma_2=\IF$}
Let $X_z^\epsilon(s,t), (s,t)\in \mathbb{R}^2, z>0, \epsilon>0$ be homogeneous Gaussian random fields with
correlation function
$$1-Corr\left(X_z^\epsilon(s,t), X_z^\epsilon(S,t_1) \right)\sim (1+\epsilon)\rho_1^2(|s-s_1|)(1+o(1))+\frac{1}{z}v_2^2(|t-t_1|)(1+o(1)), \ \ |s-s_1|,
|t-t_1|\rw\IF.$$
Thus, by Slepian inequality
\BQN\label{EQ12}
\pi_9^+(u)\leq\pi_1(u)\leq \pi_9^-(u),
\EQN
where
\BQNY
\pi_9^+(u):&=&\mathbb{P}\left(\sup_{|s|\leq \overleftarrow{\vv_1}(u^{-1})}X(s,\mu s)>u\right),\\
\pi_9^-(u):&=&\mathbb{P}\left(\sup_{(s,t)\in D_u}\frac{X_z^{\epsilon}(s,t)}{1+(1-\epsilon)\left(v_1^2(|s+b_{12}t|)+v_2^2(|t|)\right)}>u\right).
\EQNY
\KD{It is straightforward to check that
$\frac{X_z^{\epsilon}(s,t)}{1+(1-\epsilon)\left(v_1^2(|s+b_{12}t|)+v_2^2(|t|)\right)}$
satisfies assumptions of}
Case ii) $\gamma_1=0, \gamma_2=(1-\epsilon)z\in (0,\IF)$. Thus
\BQNY
\pi_9^-(u)&\leq& 2\left(\frac{z}{\theta}\right)^{1/\beta}\mathcal{H}_{\alpha_1}
\frac{\overleftarrow{\vv_1}(u^{-1})}{\overleftarrow{\LL_1}(u^{-1})}\Psi(u)\left(\int_0^{\IF}
\widehat{\mathcal{P}}_{\beta}^{g_{t,z}}e^{-\frac{z \Mba }{\theta}t^{\beta}}dt
 \right)(1+o(1)), \ \ u\rw\IF,  \eMM{ S\rw\IF, \ve \to 0},
\EQNY
with $$g_{s,z}(t)=\frac{z}{\theta}\left(|s+bt|^{\beta}+\theta|t|^{\beta}-|(1+b\mu) s|^{\beta}-\theta|\mu s|^{\beta}\right), s\geq 0, t\in \mathbb{R}.$$
Replacing $t$ by $z^{-1/\beta} s$ in the above integral yields
$$z^{1/\beta} \int_0^{\IF} \widehat{\mathcal{P}}_{\beta}^{g_{t,z}}e^{-\frac{z \Mba }{\theta}t^{\beta}}dt=\int_0^{\IF}
\widehat{\mathcal{P}}_{\beta}^{g_{z^{-1/\beta} s,z}}e^{-\frac{ \Mba }{\theta}s^{\beta}}ds.$$
Note that for any $\epsilon>0$, there exists a positive constant
$M_\epsilon>0$ such that  \KD{for} $z$ sufficiently large
$$g_{z^{-1/\beta}s,z}(t)+\epsilon|s|^\beta=\frac{1}{\theta}\left(|s+b_{12}tz^{1/\beta}|^{\beta}+\theta|tz^{1/\beta}|^{\beta}-|(1+b_{12}\mu) s|^{\beta}-\theta|\mu
s|^{\beta}\right)+\epsilon|s|^\beta\geq M_\epsilon z|t|^{\beta}, \ \ t\in \mathbb{R},$$
which implies that
$$\widehat{\mathcal{P}}_{\beta}^{g_{z^{-1/\beta} s,z}}\leq e^{\epsilon |s|^\beta}\widehat{\mathcal{P}}_{\beta}^{M_{\epsilon}z}.$$
Since
$$\lim_{z\rw\IF}\widehat{\mathcal{P}}_{\beta}^{M_{\epsilon}z}=1,$$  then using dominated convergence theorem
\BQNY
\limsup_{z\rw\IF}\int_0^{\IF} \widehat{\mathcal{P}}_{\beta}^{g_{z^{-1/\beta} s,z}}e^{-\frac{ \Mba }{\theta}s^{\beta}}ds&\leq&
\limsup_{z\rw\IF}\int_0^{\IF} \widehat{\mathcal{P}}_{\beta}^{M_{\epsilon}z}e^{-\left(\frac{ \Mba }{\theta}-\epsilon\right)s^{\beta}}ds\\
&\rw&
\left(\frac{ \Mba }{\theta}\right)^{-1/\beta}\Gamma(1/\beta+1), \ \ \epsilon\rw 0.
\EQNY
Thus we conclude that
\BQN\label{EQ13}
\pi_9^-(u)&\leq& 2\Gamma(1/\beta+1)\left( \Mba \right)^{-1/\beta}\mathcal{H}_{\alpha_1}
\frac{\overleftarrow{\vv_1}(u^{-1})}{\overleftarrow{\LL_1}(u^{-1})}\Psi(u)(1+o(1)), \ \ u\rw\IF, \eMM{ S\rw\IF, \ve \to 0}.
 \EQN
 Next we focus on $\pi_9^+(u)$. One can easily check that
 the variance and correlation \NE{functions of $X(s,\mu s)$  satisfy}
 $$1-Var\left(X(s, \mu s)\right)\sim v_1^2(|(1+b_{12}\mu) s|)+v_2^2(|\mu s|)\sim  \Mba v_1^2(|s|), \ \ s\rw 0,$$ and
 $$1-Corr\left(X(s, \mu s), X(s_1, \mu s_1)\right)\sim \rho_1^2(|s-s_1|)+\rho_2^2(|\mu(s-s_1)|)\sim \rho_1^2(|s-s_1|), \ \ s, s_1\rw 0.$$
 In light of Theorem \ref{mainT}, we have
 $$\pi_9^+(u)\sim 2\Gamma(1/\beta+1)\left( \Mba \right)^{-1/\beta}\mathcal{H}_{\alpha_1}
\frac{\overleftarrow{\vv_1}(u^{-1})}{\overleftarrow{\LL_1}(u^{-1})}\Psi(u), \ \ u\rw\IF,$$
which combined with
(\ref{EQ12}) and (\ref{EQ13}) establishes the proof.\\
\underline{Case iv) $\gamma_1\in (0,\IF), \gamma_2=\IF$.} \ \  Let $Z(s,t)$ be a homogeneous Gaussian random field with variance $1$ and correlation
function satisfying
$$1-Corr\left(Z(s,t), Z(s_1,t_1)\right)\sim 2\rho_1^2(|s-s_1|)+\rho_1^2(|t-t_1|), \ \ |s-s_1|\rw 0, |t-t_1|\rw 0,$$
and
$\widehat{I}_{0,0}(u)=[-\overleftarrow{\LL_1}(u^{-1})S, \overleftarrow{\LL_1}(u^{-1})S]\times [-\overleftarrow{\LL_1}(u^{-1})S_1,
\overleftarrow{\LL_1}(u^{-1})S_1]. $

Then, by Slepian's inequality and Lemma \ref{simple1},
\BQN
\pi_{10}^+(u)\leq \pi_1(u)\leq \pi_{10}^-(u)+\pi_{11}(u),
\EQN
where
\BQNY
\pi_{10}^{\pm}(u)&=&\mathbb{P}\left(\sup_{(s,t)\in \widehat{I}_{0,0}(u)}\frac{\overline{X}(s,t)}{1+(1\pm
\epsilon)\left(v_1^2(|s+b_{12}t|)+v_2^2(|t|)\right)}>u\right),\\
\pi_{11}(u)&=&\mathbb{P}\left(\sup_{(s,t)\in
\left(D_u-\widehat{I}_{0,0}(u)\right)}\frac{Z(s,t)}{1+\frac{\kappa_1}{2}\left(v_1^2(|s|)+v_2^2(|t|)\right)}>u\right).
\EQNY
Note that $\rho_2^2(t)=o(\rho_1^2(t))$ as $t\rw 0$ and
$$(1\pm
\epsilon)u^2\left(v_1^2(|\overleftarrow{\LL_1}(u^{-1})s+b_{12}\overleftarrow{\LL_1}(u^{-1})t|)+v_2^2(|\overleftarrow{\LL_1}(u^{-1})t|)\right)\rw (1\pm\epsilon)\gamma_1\left(|s+b_{12}t|^{\alpha_1}+\theta|t|^{\alpha_1}\right), \ \ u\rw\IF$$
uniformly with respect to $(s,t)\in [-S,S]\times[-S_1,S_1]$.
It follows from Lemma \ref{PIPI1} that
\BQNY
\pi_{10}^{\pm}(u)\sim \Psi(u)\EE{\exp{\sup_{(s,t)\in [-S,S]\times
[-S_1,S_1]}\left[\sqrt{2}B_{\alpha_1}(s)-|s|^{\alpha_1}-(1\pm\epsilon)\gamma_1\left(|s+b_{12}t|^{\alpha_1}+\theta|t|^{\alpha_1}\right)\right]}}.
\EQNY
Since
$$\lim_{S_1\rw\IF}\EE{\exp{\sup_{(s,t)\in [-S,S]\times
[-S_1,S_1]}\left[\sqrt{2}B_{\alpha_1}(s)-|s|^{\alpha_1}-(1\pm\epsilon)\gamma_1\left(|s+b_{12}t|^{\alpha_1}+\theta|t|^{\alpha_1}\right)\right]}}
=\widehat{\mathcal{P}}_{\alpha_1}^{(1\pm \epsilon)\gamma_1  M_{\alpha_1} }(S),$$
we have
$$\pi_{10}^{\pm}(u)\sim \widehat{\mathcal{P}}_{\alpha_1}^{\gamma_1 M_{\alpha_1} }\Psi(u), \ \ u\rw\IF,  S_1\rw\IF, S\rw\IF,
\eMM{\ve \to 0}.$$
\KD{Using that
$\frac{Z(s,t)}{1+\frac{\kappa_1}{2}\left(v_1^2(|s|)+v_2^2(|t|)\right)}$
satisfies the conditions of Case iii) $\gamma_1, \gamma_2\in (0,\IF)$ of Theorem \ref{th.T1},
by the same argument as given in the proof of (\ref{PG4}), we obtain that}
$\pi_{11}(u)=o\left(\Psi(u)\right),  u\rw\IF, S_1\rw\IF, S\rw\IF.$
Thus the proof is completed. \\
\underline{Case iii) $\gamma_1=\gamma_2=\IF$.} \ \
It follows from (\ref{eq5}) and (\ref{Y}) with the specific $B$ in this case that
\BQNY
\mathbb{P}\left(X(0,0)>u\right)\leq \pi_1(u)\leq \mathbb{P}\left(\sup_{(s,t)\in
D_u}Y(s,t)>u\right),
\EQNY
where $\kappa_1$ is defined in Lemma \ref{simple1}.
The Gaussian random field involved in the right hand side of the above inequality satisfies the assumption of
Case iii) $\gamma_1=\gamma_2=\IF$ in Theorem \ref{th.T1}
and
therefore it follows that
$$\mathbb{P}\left(\sup_{(s,t)\in D_u}\frac{\overline{X}(s,t)}{1+\frac{\kappa_1}{2}\left(v_1^2(|s|)+v_2^2(|t|)\right)}>u\right)\sim \Psi(u), \ \ u\rw\IF.$$
This completes the proof. \QED

\prooftheo{th.T5}
\eMM{For $\ve>0$ sufficiently small, let}  $Z^{\pm \epsilon}$ be a stationary Gaussian process with continuous trajectories, unit variance and correlation function satisfying
$$1-r_{Z^{\pm \epsilon}}(t)\sim (1\mp \epsilon) \rho_1^2(|t|), \ \ t\rw 0.$$
By Slepain's inequality, we have
\BQNY
\pi_{12}^+(u)\leq \pi_1(u)\leq \pi_{12}^-(u),
\EQNY
where
$$\pi_{12}^\pm(u)=\mathbb{P}\left(\sup_{(s,t)\in D_u}\frac{Z^{\pm\epsilon}(s)}{1+(1\pm\epsilon)\left(v_1^2(|s|)+v_2^2(|t|)\right)}>u\right).$$
By the fact that for any $u>0$,
$$\sup_{(s,t)\in D_u}\frac{Z^{\pm\epsilon}(s)}{1+(1\pm\epsilon)\left(v_1^2(|s|)+v_2^2(|t|)\right)}=\sup_{|s|\leq \overleftarrow{\vv_1}(u^{-1}\ln
u)}\frac{Z^{\pm\epsilon}(s)}{1+(1\pm\epsilon)v_1^2(|s|)},$$
we have
$$\pi_{12}^\pm(u)=\mathbb{P}\left(\sup_{|s|\leq \overleftarrow{\vv_1}(u^{-1}\ln u)}\frac{Z^{\pm\epsilon}(s)}{1+(1\pm\epsilon)v_1^2(|s|)}>u\right).$$
Applying Theorem \ref{mainT}, we establish the claims.
\QED
\prooftheo{th.T6}
Set below for $u>0$
$$D_u=\{|s|\leq \overleftarrow{\vv_1}(u^{-1}\ln u), |t|\leq
2\mu\overleftarrow{\vv_1}(u^{-1}\ln u)\}.$$
Using the same $Z^{\pm \epsilon}$ as in the proof of Theorem \ref{th.T5}, by Slepian's inequality, we have
\BQNY
\pi_{13}^+(u)\leq \pi_1(u)\leq \pi_{13}^-(u),
\EQNY
where
$$\pi_{13}^\pm(u)=\mathbb{P}\left(\sup_{(s,t)\in D_u}\frac{Z^{\pm\epsilon}(s)}{1+(1\pm\epsilon)\left(v_1^2(|s+b_{12}t|)+v_2^2(|t|)\right)}>u\right).$$
\KD{The same analysis as given between (\ref{EQ14}) and (\ref{E11}) implies that,} for $u$ sufficiently large
$$(1-\epsilon) \Mbb v_1^2(|s|)\leq \inf_{|t|\leq 2\mu\overleftarrow{\vv_1}(u^{-1}\ln u)}v_1^2(|s+b_{12}t|)+v_2^2(|t|)\leq (1+\epsilon)
\Mbb v_1^2(|s|)$$
hold for $|s|\leq \overleftarrow{\vv_1}(u^{-1}\ln u)$.  Thus we have
$$
\pi_{13}^-(u)\leq \mathbb{P}\left(\sup_{|s|\leq \overleftarrow{\vv_1}(u^{-1}\ln u)}\frac{Z^{-\epsilon}(s)}{1+(1-\epsilon)^2
\Mbb v_1^2(|s|)}>u\right),$$ and $$\pi_{13}^+(u)\geq \mathbb{P}\left(\sup_{|s|\leq \overleftarrow{\vv_1}(u^{-1}\ln
u)}\frac{Z^{+\epsilon}(s)}{1+(1+\epsilon)^2 \Mbb v_1^2(|s|)}>u\right).$$
Hence the claim follows by Theorem \ref{mainT}. \QED

\section{Appendix A}
In this section \eM{we derive \eMM{some} key uniform expansions of the tail of maximum of Gaussian random fields
 over short intervals}.   \eMM{For any $\gamma \in (0,\IF), S>0$  we define }
 $$ \eMM{\mathcal{P}_{\alpha}^{\gamma} = \EE{ \sup_{[0,S]} e^{ \sqrt{2} B_\alpha(t)- (1+ \gamma) \abs{t}^\alpha}}}$$
 and we set
$$\mathcal{P}_{\alpha}^{\IF}[0,S]=1, \quad \mathcal{P}_{\alpha}^{0}[0,S]=\mathcal{H}_{\alpha}[0,S] , \quad
\alpha\in (0,2], S>0 .$$
\eM{The claim of the following three lemmas follows by Theorem 2.1 in \cite{KEP2016}};
 the detailed proofs are omitted here.

In the following  $h_k, k\in \mathcal{K}_u$ with $\mathcal{K}_u$ an index set are  positive functions such that $\limit{u}h_k(u)/u=1$ uniformly with respect to $k\in\mathcal{K}_u$.

\BEL\label{Pickands1} Let $X_{u,k}(t), t\in [0,T], k\in \mathcal{K}_u$ be a sequence of \eM{centered} Gaussian processes with continuous trajectories, variance 1  and correlation function $r(\cdot,\cdot)$ satisfying (\ref{eR}) uniformly with respect to $k\in \mathcal{K}_u$.
Suppose that $\rho\in\mathcal{R}_\alpha/2, \vv\in \mathcal{R}_{\beta/2}$ with $0<\alpha\leq 2,
\beta>0$.  If $\lim_{t\rw 0}\frac{\vv^2(t)}{\LL^2(t)} =\gamma\in [0,\IF]$, then
\BQNY
\lim_{u\rw\IF}\sup_{k\in \mathcal{K}_u}\left|  \frac{1} {\Psi(h_k(u))} \pk{ \sup_{t\in [0,\overleftarrow{\rho}(u^{-1})S]} \frac{X_{u,k}(t)}{1+\vv^2(t)}>h_k(u)}-\mathcal{P}_{\alpha}^\gamma[0,S]\right|=0.
\EQNY
\EEL

\bigskip
Let $\rho_i\in \mathcal{R}_{\alpha_i/2}$, $v_i\in \mathcal{R}_{\beta_i/2}, i=1,2$ be non-negative functions with $0<\alpha_i\leq 2, \beta_i>0, i=1,2.$
Let $\eM{X_{u,k}(s,t)}, k\in \mathcal{K}_u$ be  \CE{centered} Gaussian random fields over $\mathcal{E}(u):=\{(\overleftarrow{\LL_1}(u^{-1})\eM{s},\overleftarrow{\LL_2}(u^{-1}) \eM{t} ) , (s,t)\in \mathcal{E}\}$ with $\mathcal{E}$ an compact set containing $0$.
Suppose further that $X_{u,k}$ has unit variance, continuous trajectories  and correlation function $r_k(s,t,s_1,t_1)$ satisfying (\ref{Cor2}) uniformly with respect to $k\in \mathcal{K}_u$.
\BEL\label{PIPI}
Let  $d_u(s,t),u>0$ be continuous functions satisfying
 \BQN\label{eq1}
 \lim_{u\rw\IF}\sup_{(s,t)\in\mathcal{E}, k\in \mathbb{K}_u}|h_k^2(u)d_u(\overleftarrow{\LL_1}(u^{-1})s, \overleftarrow{\LL_2}(u^{-1})t)-d(s,t)|=0.
 \EQN
\eM{If further $d_u(s,t) > -1$ for any $s,t \in \mathcal{E}$ and all $u$ large},  then 
\BQNY
\lim_{u\rw\IF}\sup_{k\in \mathcal{K}_u}\left| \frac{1}{\Psi(h_k(u))) }\mathbb{P}\left(\sup_{(s,t)\in \mathcal{E}(u)}\frac{X_{u,k}(s,t)}{1+d_u(s,t)}>h_k(u)\right)-
\EE{e^{\sup_{(s,t)\in \mathcal{E}}\{W_{\alpha_1,\alpha_2}(s,t)-d(s,t)\}}}\right|=0,
\EQNY
where
$ W_{\alpha_1,\alpha_2}(s,t)=\sqrt{2}(B_{\alpha_1}(s)+\tilde{B}_{\alpha_2}(t))-|s|^{\alpha_1}-|t|^{\alpha_2}, s,t\in \mathbb{R}$,
with $B_{\alpha_1}$ and $B_{\alpha_2}$ being two independent fBm's with indices $\alpha_1, \alpha_2$, respectively.
\EEL
\BEL\label{PIPI1}
 Suppose that $d_u(s,t),u>0$ are continuous functions  satisfying
  $$\lim_{u\rw\IF}\sup_{(s,t)\in\mathcal{E}, k\in \mathcal{K}_u}|h_k^2(u)d_u(\overleftarrow{\LL_1}(u^{-1})s, \overleftarrow{\LL_1}(u^{-1})t)-d(s,t)|=0.$$
  If $\rho_2^2(t)=o(\rho_1^2(t))$ as $t\rw 0$ and \eM{$d_u(s,t)> -1$ for any $s,t\in $ and all $u$ large}, then
\BQNY
\lim_{u\rw\IF}\sup_{k\in\mathcal{K}_u}\left|\frac{1}{\Psi(h_k(u)))}\mathbb{P}\left(\sup_{(s,t)\in \widetilde{\mathcal{E}}(u)}\frac{\eM{X}_{u,k}(s,t)}{1+d_u(s,t)}>h_k(u)\right)-
\EE{\sup_{(s,t)\in \mathcal{E}}e^{\sqrt{2}B_{\alpha_1}(s)-|s|^{\alpha_1}-d(s,t)}}\right|=0,
\EQNY
with
 $B_{\alpha_1}$  an fBm with index $\alpha_1$ and $\widetilde{\mathcal{E}}(u):=\{(\overleftarrow{\LL_1}(u^{-1})\eM{s},\overleftarrow{\LL_1}(u^{-1})\eM{t}), (s,t)\in \mathcal{E}\}$.
\EEL
Assume now that $X(t), t=(t_1 \ldot t_d)\in \mathbb{R}^d$ is a Gaussian field with continuous trajectories, unit variance and covariance function satisfying
\BQN\label{C1}
1-Cov(X(s), X(t))\sim \sum_{i=1}^d\rho_i(|t_i-s_i|),  \quad s, t\rw 0,
\EQN
with $\rho_i $ positive regularly varying function with index $\alpha_i/2\in (0,1]$.
Denote by  $\overleftarrow{\rho}_{\vk{d}}(u^{-1})=(\overleftarrow{\rho}_1(u^{-1}), \dots, \overleftarrow{\rho}_d(u^{-1}))$ and define
$\overleftarrow{\rho}_{\vk{d}}(u^{-1})t=(\overleftarrow{\rho}_1(u^{-1})t_1, \dots, \overleftarrow{\rho}_d(u^{-1})t_d)$. Moreover,
set
$F(A,B)=\inf_{s\in A, t\in B}\sum_{i=1}^d|s_i-t_i|$ for any $A, B \subset \mathbb{R}^d$ and let
$$
D_u=\prod_{i=1}^d[-\frac{\delta_u}{\overleftarrow{\rho}_i(u^{-1})}, \quad
 \frac{\delta_u}{\overleftarrow{\rho}_i(u^{-1})}], \quad
\mathbb{K}=\{(\lambda_1, \lambda_2)\in D_u\times D_u, \lambda_1+\mathcal{E}_1, \lambda_2 +\mathcal{E}_2\subset D_u\}.$$
 Further, let $u_\lambda, \NE{\lambda\in D_u},$ with $\delta_u\rw 0, u\rw\IF$ satisfy
$$\lim_{u\rw\IF}\sup_{\lambda\in \NE{D_u}}\left|\frac{u_\lambda}{u}-1\right|=0.$$
\COM{Denote by
$$ D(\lambda_1, \lambda_2, \mathcal{E}_1, \mathcal{E}_2, u):=\pk{ \sup_{t\in \overleftarrow{\rho}_{\vk{d}}(u^{-1})(\lambda_1+\mathcal{E}_1)} X(t)> u_{\lambda_1}, \sup_{s\in \overleftarrow{\rho}_{\vk{d}}(u^{-1})(\lambda_2 +\mathcal{E}_2)} X(s)> u_{\lambda_2}}$$
with $\mathcal{E}_i\subset [0,S]^d, i=1,2$.
}
\HEH{We state next the result of Corollary 3.2 in \cite{KEP2016}, below  $\mathcal{E}_1, \mathcal{E}_2$ are assumed to be compact sets.}
 \BEL\label{uni} Suppose that $X(t), t=(t_1\ldot t_d)\in \mathbb{R}^d$ is a Gaussian field with continuous trajectories, unit variance and covariance function satisfying (\ref{C1}). Then there exists $\mathcal{C}, \mathcal{C}_1>0$ such that
 for
 any $S>1$ as $u$ sufficiently large,
 \BQNY
 \sup_{(\lambda_1,\lambda_2)\in \mathbb{K} , \mathcal{E}_1, \mathcal{E}_2\subset[0,S]^d}   e^{\frac{\mathcal{C}_1F^{\beta^*/2}(\lambda_1+\mathcal{E}_1, \lambda_2+\mathcal{E}_2)}{16}}  \frac{  \pk{ \sup_{t\in \overleftarrow{\rho}_{\vk{d}}(u^{-1})(\lambda_1+\mathcal{E}_1)} X(t)> u_{\lambda_1}, \sup_{s\in \overleftarrow{\rho}_{\vk{d}}(u^{-1})(\lambda_2 +\mathcal{E}_2)} X(s)> u_{\lambda_2}}}{ \HEH{S}^{2d} \Psi(\HEH{u}_{\lambda_1, \lambda_2}(u)) }
& \le & \mathcal{C}
 \EQNY
 with $u_{\lambda_1, \lambda_2}=\min(u_{\lambda_1}, u_{\lambda_2})$ and $\beta^*=\min_{i=1,\dots, d}\alpha_i$.
 \EEL
\section{Appendix B}
Let in this section $g$ be a positive function \KD{ such that
$$\lim_{u\rw\IF}g(u)=\IF, \quad \lim_{u\rw\IF}\frac{g(u)}{u}=0.$$
 Further, let $v \in \mathcal{R}_\beta, \beta>0$ be a non-negative function. }
We set throughout in the following
\BQN \label{tu}
t(u):= \overleftarrow{\vv}\left(\frac{g(u)}{u}\right), \quad u>0.
\EQN
We shall investigate first the asymptotic behaviour of an integral determined by $g$ and $v$.
\BEL\label{P1} i) For any $0<x\leq y<\IF$ and $c>0$, \KD{as $u\to\infty$}
\BQNY
\int_0^{xt(u)}e^{-c u^2v^2(t)}dt\sim \int_0^{yt(u)}e^{-c u^2v^2(t)}dt.
\EQNY
ii) If $a \in \mathcal{R}_\beta$ \KD{is} such that $a(t)\sim \vv(t)$ as $t\rw 0$, \KD{then as $u\to\infty$}
\BQNY
\int_0^{t(u)}e^{-u^2v^2(t)}dt \sim\int_0^{\overleftarrow{a}(u^{-1}g(u))}e^{-u^2a^2(t)}dt.
\EQNY
\EEL
\def\ttu{t(u)}
\prooflem{P1} i)
\KD{Using standard properties of regularly varying functions, see e.g., \cite{EKM97}, }
for $u$ sufficiently large and
and  $0<x<y<\IF$, we have
\BQNY
\int_{x \ttu }^{y\ttu }e^{-c u^2\vv^2(t)}dt
&\leq& e^{-cu^2\vv^2((x/2)\ttu )}(y-x)\ttu \\
&\leq& e^{-(x/3)^{2\beta}cu^2\vv^2(\ttu )}(y-x)\ttu \\
&\leq & e^{-(x/4)^{2\beta}c (g(u))^2}(y-x)\ttu
\EQNY
and
\BQNY
\int_{0}^{x\ttu }e^{-c u^2\vv^2(t)}dt&\geq& \int_{0}^{(x/8)\ttu }e^{-c u^2\vv^2(t)}dt\\
&\geq&e^{-c u^2\vv^2((x/7)\ttu )}(x/8)\ttu \\
&\geq& e^{-(x/6)^{2\beta}cu^2\vv^2(\ttu )}(x/8)\ttu \\
&\geq & e^{-(x/5)^{2\beta}c(g(u))^2}(x/8)\ttu ,
\EQNY
which imply that, \KD{as $u\to\infty$},
\BQNY
\int_{0}^{x\ttu }e^{-c u^2\vv^2(t)}dt\sim \int_{0}^{y \ttu }e^{-c u^2\vv^2(t)}dt.
\EQNY
ii)  For any $0<\epsilon<1/2$
$$(1-\epsilon) a(t)\leq \vv(t)\leq (1+\epsilon)a(t)$$
holds for $t$ sufficiently small. Consequently,  for $u$ sufficiently large
\BQNY
\int_{0}^{\ttu }e^{-u^2\vv^2(t)}dt&\leq& \int_{0}^{\ttu }e^{-(1-\epsilon)^2u^2 a^2(t)}dt\leq \int_{0}^{\ttu }e^{-u^2a^2((1-2\epsilon)^{1/\beta}t)}dt\\
&=& (1-2\epsilon)^{-1/\beta}\int_{0}^{(1-2\epsilon)^{1/\beta}\ttu }e^{-u^2 a^2(t)}dt\\
&\leq& (1-2\epsilon)^{-1/\beta}\int_{0}^{\ttu }e^{-u^2 a^2(t)}dt
\EQNY
and
\BQNY
\int_{0}^{\ttu }e^{-u^2\vv^2(t)}dt&\geq& \int_{0}^{\ttu }e^{-(1+\epsilon)^2u^2 a^2(t)}dt\geq \int_{0}^{t(u)}e^{-u^2
a^2((1+2\epsilon)^{1/\beta}t)}dt\\
&=& (1+2\epsilon)^{-1/\beta}\int_{0}^{(1+2\epsilon)^{1/\beta}\ttu }e^{-u^2 a^2(t)}dt\\
&\geq& (1+2\epsilon)^{-1/\beta}\int_{0}^{\ttu }e^{-u^2 a^2(t)}dt.
\EQNY
Letting $\epsilon\rw 0$ and by the fact that $\overleftarrow{a}(u^{-1}g(u))\sim t_u$, we  establish  the second claim. \QED

\BEL \label{integ} We have
\BQN
\int_0^{\ttu }e^{-cu^2v^2(t)}dt \sim c^{-1/(2\beta)}\Gamma(1+1/(2\beta)) \overleftarrow{v}(1/u), 
\quad u\to \IF.
\EQN
\EEL

\prooflem{integ} By \nelem{P1}, ii) we can assume that $v(x)= \ell(x) x^\beta$ with $\ell$ normalized slowly varying function at 0. It is well-known that
$\ell(x) x^\beta$ is ultimately monotone for any $\beta\not=0$, $\ell$ is continuously differentiable and
\BQN\label{ldif}
 \lim_{x\to0}  \frac{x \ell'(x)}{\ell(x)}=0.
 \EQN
Since $v$ is ultimately monotone, we have with \eM{$g(u)$ and $t(u)$  defined by \eqref{tu}}
\BQN
\int_0^{\ttu }e^{-cu^2v^2(t)}dt \sim u^{-1} \int_0^{g(u)} \frac{1}{v'(\overleftarrow{v}(y/u))} e^{- cy^2}\, dy, \quad u\to \IF.
\EQN
Further,  \eqref{ldif} implies
$$  \frac{1}{v'(\overleftarrow{v}(y/u))} \sim \frac{1}{\beta} \frac {\overleftarrow{v}(y/u)) } {v(\overleftarrow{v}(y/u))}
\sim \frac{1}{\beta}\frac{u}{y} \overleftarrow{v}(y/u)  $$
Consequently, as $u\to \IF$
\BQNY
\int_0^{\ttu}e^{-cu^2v^2(t)}dt
&\sim& \frac{1}{\beta} \int_0^{g(u)} \overleftarrow{v}(y/u)y^{-1}  e^{-c y^2}\, dy\notag\\
&\sim&  \frac{1}{\beta} \overleftarrow{v}(1/u)\int_0^{g(u)} \frac{\overleftarrow{v}(y/u)}{\overleftarrow{v}(1/u)}y^{-1}  e^{-c y^2}\, dy\notag.\\
\EQNY
Potter's theorem  shows that  there exists a constant $C$ such that for $u$ sufficiently large,
$$\frac{\overleftarrow{v}(y/u)}{\overleftarrow{v}(1/u)}\leq C(\max(1,y))^{2/\beta}, \ \ 0\leq y\leq g(u).$$
By the fact that for any $y>0$
$$ \limit{u} \frac{\overleftarrow{v}(y/u)}{\overleftarrow{v}(1/u)}= y^{1/\beta}$$
and the dominated convergence theorem, \eM{since $\limit{u} g(u)=\IF$}, we obtain
\BQNY
\int_0^{\ttu}e^{-u^2v^2(t)}dt
&\sim& \frac{1}{\beta}\overleftarrow{v}(1/u) \int_0^{\IF}  y^{1/\beta-1}  e^{-c y^2}\, dy \\ 
&\sim& c^{-1/(2\beta)}\Gamma(1+1/(2\beta)) \overleftarrow{v}(1/u)\notag. 
\EQNY
\eM{Note that alternatively, by \cite{Soulier}[Proposition 1.18] it follows that}
$$\int_0^{g(u)} \frac{\overleftarrow{v}(y/u)}{\overleftarrow{v}(1/u)}y^{-1}  e^{-c y^2}\, dy
\sim \int_0^{g(u)} y^{1/\beta-1}  e^{-c y^2}\, dy, \quad u\to \IF$$
and thus again the claim follows.} \QED

\BEL\label{simple} Suppose that $\rho^2_1\in \mathcal{R}_{\alpha_1}$ and $\rho_2^2\in \mathcal{R}_{\alpha_2}$ with $\alpha_1, \alpha_2 >0$. If
$\rho_1^2(|t|)=o(\rho_2^2(|t|))$ as $t\rw 0$, then for any $a, b \in \mathbb{R}$,
\BQNY
\rho_1^2(|as +b t|)+\rho_2^2(|t|)\sim \rho_1^2(|as|)+\rho_2^2(|t|),  \quad  s, t\rw 0.
\EQNY
\EEL
\prooflem{simple} The claim follows \eM{easily if $abst=0$}. Next we suppose that
$abst\neq 0$.
It suffices to prove that
$$\lim_{s,t\rw 0, st\neq 0}\frac{\abs{\rho_1^2(|as +b t|)-\rho_1^2(|as|)}}{\rho_1^2(|as|)+\rho_2^2(|t|)}=0.$$
For any $\epsilon\in (0,\alpha_1)$, if $\abs{\frac{as}{bt}}>\frac{4\alpha_1}{\epsilon}$, then
$$1-\frac{\epsilon}{4\alpha_1}\leq\frac{|as+bt|}{|as|}\leq 1+\frac{\epsilon}{4\alpha_1}.$$
Thus in light of UCT, we have, for $s,t$ sufficiently small
\BQNY
\frac{\abs{\rho_1^2(|as +b t|)-\rho_1^2(|as|)}}{\rho_1^2(|as|)+\rho_2^2(|t|)}&\leq& \frac{\rho_1^2(|as|)\abs{\frac{\rho_1^2(|as +b
t|)}{\rho_1^2(|as|)}-1}}{\rho_1^2(|as|)+\rho_2^2(|t|)}\leq \ABs{\frac{\rho_1^2(|as +b t|)}{\rho_1^2(|as|)}-1}\\
&\leq& \max\left(\left(1+\frac{\epsilon}{2\alpha_1}\right)^{\alpha_1}-1, 1-\left(1-\frac{\epsilon}{2\alpha_1}\right)^{\alpha_1}\right)
=: b_{\epsilon}.
\EQNY
For any $\epsilon \in (0,\alpha_1)$, if $\abs{\frac{as}{bt}}\leq\frac{4\alpha_1}{\epsilon}$, then
$$
\frac{|as+bt|}{|bt|}\leq 1+\frac{4\alpha_1}{\epsilon}.$$
\eM{Applying} again UCT we obtain
\BQNY
\frac{\abs{\rho_1^2(|as +b t|)-\rho_1^2(|as|)}}{\rho_1^2(|as|)+\rho_2^2(|t|)}&\leq&
\frac{\rho_1^2(|bt|)}{\rho_1^2(|as|)+\rho_2^2(|t|)}\ABs{\frac{\rho_1^2(|as +b t|)}{\rho_1^2(|bt|)}-\frac{\rho_1^2(|as|)}{\rho_1^2(|bt|)}}\\
&\leq&
\frac{\rho_1^2(|bt|)}{\rho_2^2(|t|)}\left[\left(1+\frac{8\alpha_1}{\epsilon}\right)^{\alpha_1}+\left(\frac{8\alpha_1}{\epsilon}\right)^{\alpha_1}\right]\\
&\rw& 0, \ \ \text{as}\ \  st\neq 0, s,t\rw 0.
\EQNY
Consequently,
$$\lim_{s,t\rw 0, st\neq 0}\frac{\abs{\rho_1^2(|as +b t|)-\rho_1^2(|as|)}}{\rho_1^2(|as|)+\rho_2^2(|t|)}\leq b_{\epsilon}\rw 0,\ \ \epsilon\rw 0.$$
This completes the proof. \QED

\BEL\label{simple1} Suppose that $v^2_1, v^2_2\in \mathcal{R}_{\beta}$ $\beta >0$. If $a_1 v^2_2(|t|)\leq v_1^2(|t|)\leq a_2 v_2^2(|t|))$ with $a_1,
a_2>0$ for $t$ sufficiently small, then for any reversible matrix
$B=\left(\begin{array}{cc}
 b_{11}& b_{12}\\
 b_{21}& b_{22}\\
 \end{array}
 \right),$
there exist two positive constants $\kappa_1$ and $\kappa_2$ such that
$$\kappa_1v^2_1(|s|)+\kappa_1v^2_2(|t|)\leq v^2_1(|b_{11}s+b_{12}t|)+v^2_2(|b_{21}s+b_{22}t|)\leq \kappa_2v^2_1(|s|)+\kappa_2v^2_2(|t|)$$
holds in a neighbourhood of \Kd{0}.
\EEL
\prooflem{simple1}  Without loss of generality, we assume that $|t|\geq |s|$ and $|t|>0$.
By UCT, we have
\BQNY
\frac{ v^2_1(|b_{11}s+b_{12}t|)+v^2_2(|b_{21}s+b_{22}t|)}{v^2_1(|s|)+v_2^2(|t|)}&\leq& \frac{
a_2v^2_2(|t|(|b_{12}|+|b_{11}\frac{s}{t}|))+v^2_2(|t|(|b_{22}|+|b_{21}\frac{s}{t}|))}{v^2_2(|t|)}\\
&\leq&2\left(a_2(|b_{11}|+|b_{12}|)^{\beta}+(|b_{21}|+|b_{22}|)^{\beta}\right),
\EQNY
for $t$ sufficiently small. Hence we get the upper bound. For the lower bound, making  a linear transformation
$$(s,t)^\top= B^{-1}(s',t')^\top=\left(\begin{array}{cc}
 b_{11}'& b_{12}'\\
 b_{21}'& b_{22}'\\
 \end{array}
 \right)(s',t')^\top,$$
  and then using the above conclusion, we have
 \BQNY
 v^2_1(|s|)+v^2_2(|t|)&=& v^2_1(|b_{11}'s'+b_{12}'t'|)+v^2_2(|b_{21}'s'+b_{22}'t'|)\\
 &\leq& 2\left(a_2(|b_{11}'|+|b_{12}'|)^{\beta}+(|b_{21}'|+|b_{22}'|)^{\beta}\right)(v^2_1(|s'|)+v^2_2(|t'|))\\
 &\leq& 2\left(a_2(|b_{11}'|+|b_{12}'|)^{\beta}+(|b_{21}'|+|b_{22}'|)^{\beta}\right)\left(v^2_1(|b_{11}s+b_{12}t|)+v^2_2(|b_{21}s+b_{22}t|)\right),
 \EQNY
provided $|t'|\geq |s'|$ and $|t'|>0$ for $t'$ sufficiently small. This completes the proof.
\QED

{\textbf {Proof of (\ref{EQ2})}}.
Note that
\BQN\label{EQ14}
v_1^2(|s+b_{12}t|)+v_2^2(|t|)=v_1^2(|s|)\frac{v_1^2(|s||1+b_{12}t/s|)+v_2^2(|s||t/s|)}{v_1^2(|s|)}, \ \ (s,t)\in D_u.
\EQN
If $|t/s|\leq M<\IF$, then by UCT
$$\sup_{(s,t)\in D_u, |t/s|\leq M }\left|\frac{v_1^2(|s||1+b_{12}t/s|)+v_2^2(|s||t/s|)}{v_1^2(|s|)}-|1+b_{12}t/s|^{\beta}-\theta|t/s|^{\beta}\right|\rw 0, \ \ u\rw
\IF.$$
If $|t/s|\geq M$, then using Potter's bound, for $u$ sufficiently large
$$\inf_{(s,t)\in D_u, |t/s|\geq M}\frac{v_1^2(|s||1+b_{12}t/s|)+v_2^2(|s||t/s|)}{v_1^2(|s|)}\geq 1/2\left(||b_{12}|M-1|^{\beta/2}+\theta M^{\beta/2}\right).$$
Therefore, the minimum of $v_1^2(|s+b_{12}t|)+v_2^2(|t|)$ \KD{is attained for}
$|t/s|\leq M$ with $M$ sufficiently large. Further, the minimum of
$|1+b_{12}t/s|^{\beta}+\theta|t/s|^{\beta}$ is attained at $\mu :=t/s\in [-1/|b_{12}|, 1/|b_{12}|]$.  Thus, for $(s,t)\in D_u$ and $u$ sufficiently large
\BQN\label{E11}
\frac{v_1^2(|s+b_{12}t|)+v_2^2(|t|)}{v_1^2(|(1+b_{12}\mu)s|)+v_2^2(|\mu
s|)}&=&\frac{v_1^2(|s+b_{12}t|)+v_2^2(|t|)}{v_1^2(|s|)}\frac{v_1^2(|s|)}{v_1^2(|(1+b_{12}\mu)s|)+v_2^2(|\mu s|)}\nonumber\\
&\geq& \frac{|1+b_{12}\mu|^{\beta}+\theta|\mu|^{\beta}}{2}\frac{1}{2(|1+b_{12}\mu|^{\beta}+\theta|\mu|^{\beta})}=1/4,
\EQN
Recall that $v(s,t)=v_1^2(|s+b_{12}t|)+v_2^2(|t|)-v_1^2(|(1+b_{12}\mu)s|)-v_2^2(|\mu s|), \ \
(s,t)\in D_u$.
Note that $v(s,t)$ may be negative at some point. It follows that
\BQNY
 &&\left[1+(1-2\epsilon)\left(v_1^2(|(1+b_{12}\mu)s|)+v_2^2(|\mu s|)\right)\right]\left[1+(1-2\epsilon)v(s,t)\right]\\
 && \ \ =1+(1-2\epsilon)\left(v_1^2(|s+b_{12}t|)+v_2^2(|t|)\right)+(1-2\epsilon)^2\left(v_1^2(|(1+b_{12}\mu)s|)+v_2^2(|\mu s|)\right)v(s,t).
\EQNY
\KD{Moreover} (\ref{E11}) yields that
$$\left(v_1^2(|(1+b_{12}\mu)s|)+v_2^2(|\mu s|)\right)v(s,t)=o\left(v_1^2(|s+b_{12}t|)+v_2^2(|t|)\right), \ \ (s,t)\in D_u, u\rw\IF.$$
Thus, we have for any $0<\epsilon<1/4$ \KD{and} sufficiently large $u$
\BQN\label{eq7}
\left[1+(1-2\epsilon)\left(v_1^2(|(1+b_{12}\mu)s|)+v_2^2(|\mu s|)\right)\right]\left[1+(1-2\epsilon)v(s,t)\right]\leq
1+(1-\epsilon)\left(v_1^2(|s+b_{12}t|)+v_2^2(|t|)\right),
\EQN
with $(s,t)\in D_u$. Since for
$|s|\in [\overleftarrow{\LL_2}(u^{-1})x/2, \overleftarrow{\LL_2}(u^{-1})2y], |t|\in [M\overleftarrow{\LL_2}(u^{-1}), \overleftarrow{\vv}_2(\ln
u/u)]$
\BQNY
v(s,t)&=&v_1^2(|t|)\frac{v_1^2(|s+b_{12}t|)+v_2^2(|t|)-v_1^2(|(1+b_{12}\mu)s|)-v_2^2(|\mu s|)}{v_1^2(|t|)}\\
&\sim& v_1^2(|t|)\left(\left|b_{12}+\frac{s}{t}\right|^{\beta}+\theta-|1+b_{12}\mu|^{\beta}\left|\frac{s}{t}\right|^{\beta}-\theta\left|\mu
\frac{s}{t}\right|^{\beta}\right), \ \ u\rw\IF
\EQNY
then for $M,u$ sufficiently large
\BQN\label{EQ1}
v(s,t)\geq \frac{1-3\epsilon}{1-2\epsilon}v(s_1,t_1),
\EQN
with $|s|, |s_1|\in [\overleftarrow{\LL_2}(u^{-1})x/2, \overleftarrow{\LL_2}(u^{-1})2y], |t|, |t_1|\in [M\overleftarrow{\LL_2}(u^{-1}),
\overleftarrow{\vv}_2(\ln u/u)]$.

Moreover, for any $\epsilon_1>0$, $|s|\in [\overleftarrow{\LL_2}(u^{-1})x/2, \overleftarrow{\LL_2}(u^{-1})2y]$ and $|t|\in [0,
M\overleftarrow{\LL_2}(u^{-1})]$, by UCT
\BQNY
v(s,t)&\geq&v_1^2(|s|)\left[(1-\epsilon_1)\left(\left|1+b_{12}t/s\right|^{\beta}+\theta\left|t/s\right|^{\beta}\right)-(1+\epsilon_1)
\left(|1+b_{12}\mu|^{\beta}+\theta|\mu|^{\beta}\right)\right]\nonumber\\
&\geq& -2\epsilon_1 \left(|1+b_{12}\mu|^{\beta}+\theta|\mu|^{\beta}\right)v_1^2(|s|),
\EQNY
and for any $|s|, |s_1|\in [\frac{i-1}{n}\overleftarrow{\LL_2}(u^{-1}), \frac{i+2}{n}\overleftarrow{\LL_2}(u^{-1})]$ with $x/2\leq \frac{i}{n}\leq 2y$
and $|t|\in [0, M\overleftarrow{\LL_2}(u^{-1})]$ and $u$ and $n$ sufficiently large
\BQNY
&&|v(s,t)-v(s_1,t)|
 \leq v_1^2(|s|) \sup_{d_1,d_2\in \{\pm \epsilon_1\}}\left|(1+d_1)\left(|1+b_{12}t/s|^{\beta}+|(1+b_{12}\mu)s_1/s|^{\beta}+\theta|\mu
s_1/s|^{\beta}\right)\right.\\
&& \ \ \left.-(1+d_2)\left(|1+b_{12}\mu|^{\beta}+\theta|\mu|^{\beta}+|s_1/s+b_{12}t/s|^{\beta}\right)\right|\\
&& \ \ \leq v_1^2(|s|)\mathbb{Q}\epsilon_1\left(|1+b_{12}\mu|^{\beta}+\theta|\mu|^{\beta}+|1+2|b_{12}|M/x|^{\beta}\right)+v_1^2(|s|)||s_1/s|^{\beta}-1|
\left(|1+b_{12}\mu|^{\beta}+\theta|\mu|^{\beta}\right)\\
&& \ \ \  \ +v_1^2(|s|)\sup_{|z|\in [0, 4M/x] }|h_{s_1/s}( z)-h_1(z)|,
\EQNY
where $h_s(z)=|s+b_{12}z|^{\beta}, s, z\in \mathbb{R}$. Therefore, for $|s|, |s_1|\in [\frac{i-1}{n}\overleftarrow{\LL_2}(u^{-1}), \frac{i+2}{n}\overleftarrow{\LL_2}(u^{-1})]$ with $x/2\leq \frac{i}{n}\leq 2y$
and $|t|\in [0, M\overleftarrow{\LL_2}(u^{-1})]$ when $\epsilon_1$ sufficiently small and $u$ and $n$ sufficiently large
$$v(s,t)\geq -\epsilon/4 \left(|1+b_{12}\mu|^{\beta}+\theta|\mu|^{\beta}\right)v_1^2(|s|), \ \ \ |v(s,t)-v(s_1,t)|\leq \epsilon/8
\left(|1+b_{12}\mu|^{\beta}+\theta|\mu|^{\beta}\right)v_1^2(|s|),$$
which implies that (\eM{recall that $\limit{u}\sup_{(s,t)\in D_u}|v(s,t)|= 0$)}
\BQNY
&&\epsilon \left(v_1^2(|(1+b_{12}\mu)s|)+v_2^2(|\mu s|)\right)+\epsilon v(i\overleftarrow{\LL_2}(u^{-1})/n,t)\\
&& \ \ \ \geq (1-2\epsilon)\left(v(i\overleftarrow{\LL_2}(u^{-1})/n,t)-v(s,t)\right)-(1-2\epsilon)^2\left(v_1^2(|(1+b_{12}\mu)s|)+v_2^2(|\mu
s|)\right)v(s,t)\\
&& \ \ \ \ \ +(1-3\epsilon)^2\left(v_1^2(|(1+b_{12}\mu)s|)+v_2^2(|\mu s|)\right)v(i\overleftarrow{\LL_2}(u^{-1})/n,t).
\EQNY
Hence, \KD{combining the above} with (\ref{eq7}) and (\ref{EQ1})  for any $0<\epsilon<1/4$, we have for $n,u$ sufficiently large,
\BQNY
 1+(1-\epsilon)\left(v_1^2(|s+b_{12}t|)+v_2^2(|t|)\right)&\geq& \left[1+(1-2\epsilon)\left(v_1^2(|(1+b_{12}\mu)s|)+v_2^2(|\mu
 s|)\right)\right]\left[1+(1-2\epsilon)v(s,t)\right]\nonumber\\
 &\geq&  \left[1+(1-3\epsilon)\left(v_1^2(|(1+b_{12}\mu)s|)+v_2^2(|\mu s|)\right)\right]\left[1+(1-3\epsilon)v(i\overleftarrow{\LL_2}(u^{-1})/n,t)\right],
\EQNY
holds for $|s|\in [\frac{i-1}{n}\overleftarrow{\LL_2}(u^{-1}), \frac{i+2}{n}\overleftarrow{\LL_2}(u^{-1})]$ with $x/2\leq \frac{i}{n}\leq 2y$ and  and
$|t|\in [0, \overleftarrow{\vv}_2(\ln u/u)]$, which completes the proof. \QED\\
{\bf Acknowledgement}: Thanks to Swiss National Science Foundation grant No. 200021-166274.  KD acknowledges partial support
by NCN Grant No 2015/17/B/ST1/01102 (2016-2019).

\bibliographystyle{plain}

 \bibliography{RRRS}

\end{document}